\documentclass[3p,computermodern,10pt]{elsarticle}
\usepackage{amsmath,amssymb,amsthm}
\usepackage{bm}
\usepackage{mathtools}
\usepackage{graphicx}
\usepackage{xcolor}
\usepackage{placeins}
\definecolor{bestkl}{RGB}{0,100,0}
\usepackage{subcaption}
\usepackage{hyperref}
\usepackage{enumitem}
\usepackage{algorithm}
\usepackage{algpseudocode}
\usepackage{comment}

\usepackage{tikz}
\usetikzlibrary{positioning, fit, arrows.meta, shapes}

\usetikzlibrary{
  arrows.meta, %
  fit, %
  positioning, %
}
\tikzset{
  neuron/.style={ %
    circle,draw,thick, %
    inner sep=0pt, %
    minimum size=3.5em, %
    node distance=1ex and 2em, %
  },
  group/.style={ %
    rectangle,draw,thick, %
    inner sep=0pt, %
  },
  io/.style={ %
    neuron, %
    fill=gray!15, %
  },
  conn/.style={ %
    -{Straight Barb[angle=60:2pt 3]}, %
    thick, %
  },
}

\newcommand{\HTwo}{\mathrm{H}_2}
\newcommand{\OTwo}{\mathrm{O}_2}
\newcommand{\HTwoO}{\mathrm{H}_2\mathrm{O}}
\newcommand{\defeq}{:=}

\newcommand{\residual}{\boldsymbol r}
\newcommand{\state}{\boldsymbol u}

\newcommand{\stateToObservableMap}{\boldsymbol s}

\newcommand{\obs}{\mathbf{d}}
\newcommand{\y}{\mathbf{y}}
\newcommand{\x}{\mathbf{x}}
\newcommand{\xDim}{d}
\newcommand{\nParams}{d}

\newcommand{\RR}[1]{\mathbb{R}^{#1}}

\newcommand{\variationalMean}{\boldsymbol \mu}
\newcommand{\variationalStd}{\boldsymbol \sigma}
\newcommand{\variationalLogStd}{\boldsymbol \eta}

\newcommand{\variationalStdScalar}{\sigma}

\newcommand{\fisherMatrix}{\mathbf{F}}
\newcommand{\nSamples}{\lambda}

\newcommand{\ELBOHF}{\ELBO_{\mathrm{HF}}}
\newcommand{\ELBOLF}{\ELBO_{\mathrm{LF}}}

\newcommand{\nSamplesHF}{\nSamples_{\mathrm{HF}}}
\newcommand{\nSamplesLF}{\nSamples_{\mathrm{LF}}}

\newcommand{\searchDistribution}{q}
\newcommand{\searchDistributionParameters}{\boldsymbol \zeta}
\newcommand{\elbo}{\mathcal{L}}
\newcommand{\ELBO}{\elbo}

\newcommand{\params}{\boldsymbol \theta}
\newcommand{\paramsSpace}{\mathcal{P}}

\newcommand{\gradient}{\boldsymbol g}
\newcommand{\monteCarloGradient}{\gradient^{\mathrm{MC}}}
\newcommand{\controlVariateGradient}{\gradient^{\mathrm{CV}}}
\newcommand{\approximateControlVariateGradient}{\widehat{\gradient}}%
\newcommand{\hessian}{\mathbf{H}}
\newcommand{\approximateControlVariateHessian}{\widehat{\hessian}}%
\newcommand{\laggedHessian}{\bar{\hessian}}
\newcommand{\regularizedHessian}{\widetilde{\hessian}}
\newcommand{\fallbackLearningRate}{\gamma}
\newcommand{\additiveRegularization}{\delta}

\newcommand{\acvGain}{\mathbf{A}}
\newcommand{\gradientAcvCorrection}{\boldsymbol{\chi}}
\newcommand{\biasCheck}{\textcolor{green!50!black}{\ensuremath{\checkmark}}}
\newcommand{\biasCross}{\textcolor{red!80!black}{\ensuremath{\times}}}

\newcommand{\controlVariate}{\mathbf{c}}

\newcommand{\score}[2]{s(#1;#2)}
\newcommand{\scoreH}[2]{S(#1;#2)}
\newcommand{\elboGradientIntegrand}{\boldsymbol \Xi}
\newcommand{\elboGradientIntegrandHF}{\elboGradientIntegrand_{\mathrm{HF}}}
\newcommand{\elboGradientIntegrandLF}{\elboGradientIntegrand_{\mathrm{LF}}}
\newcommand{\elboHessianIntegrand}{\boldsymbol{\Psi}}
\newcommand{\elboHessianIntegrandHF}{\elboHessianIntegrand_{\mathrm{HF}}}
\newcommand{\elboHessianIntegrandLF}{\elboHessianIntegrand_{\mathrm{LF}}}

\newcommand{\looGradientBaseline}[1]{\boldsymbol b_{g,#1}}
\newcommand{\looHessianBaseline}[1]{\boldsymbol b_{H,#1}}
\newcommand{\looElboGradientIntegrand}[1]{\elboGradientIntegrand_{#1}^{\mathrm{LOO}}}
\newcommand{\looElboHessianIntegrand}[1]{\elboHessianIntegrand_{#1}^{\mathrm{LOO}}}
\newcommand{\hessianAcvGain}[1]{\mathcal{B}_{H,#1}}

\newcommand{\rom}{\mathcal{M}}
\newcommand{\romParams}{\mathbf{w}}
\newcommand{\romTrainingMap}{\boldsymbol w}

\newcommand{\fomQoi}{\mathbf{y}}

\newcommand{\sampleSet}{\mathcal{S}}
\newcommand{\sampleSetHighFidelity}{\mathcal{S}^{\mathrm{HF}}}
\newcommand{\sampleSetLowFidelity}{\mathcal{S}^{\mathrm{LF}}}

\newcommand{\sampleSetGain}{\mathcal{S}^{\mathrm{gain}}}

\newcommand{\romError}[1]{\mathcal{E}_{#1}}
\newcommand{\romErrorFloor}{\varepsilon_{\mathrm{ROM}}}
\newcommand{\blockGain}[1]{\alpha_{#1}}
\newcommand{\componentGain}[1]{a_{#1}}
\newcommand{\podSingularValue}[1]{\varsigma_{#1}}
\newcommand{\podRank}[1]{n_{\mathrm{POD},#1}}
\newcommand{\podSingularMatrix}[1]{\mathbf{S}_{#1}}

\newcommand{\convectionAngle}{\varphi}
\newcommand{\reactionCoeff}{\kappa_{\mathrm{r}}}

\newtheoremstyle{smallitremark}%
  {6pt}%
  {6pt}%
  {\itshape}%
  {}%
  {\bfseries}%
  {.}%
  {0.5em}%
  {}%
\theoremstyle{smallitremark}
\newtheorem{remark}{Remark}
\theoremstyle{plain}

\begin{document}
\begin{frontmatter}
\title{Multifidelity variational inference for derivative-free PDE-constrained inverse problems}
\author[a]{Eric J. Parish}
\ead{ejparis@sandia.gov}
\address[a]{Sandia National Laboratories, Livermore, CA, USA}
\begin{abstract}
    We develop a multifidelity black-box variational inference (MF-BBVI) approach for PDE-constrained Bayesian inverse problems in which forward solves are expensive and PDE gradients are unavailable. The method is built on score-function estimators for the evidence lower bound (ELBO), which avoid differentiation through the PDE model but are known to suffer from high variance. To address this high variance, we employ multifidelity approximate control-variate estimators that combine a small number of high-fidelity samples with a larger number of low-fidelity samples to reduce estimator variance and the number of high-fidelity evaluations required for optimization. We further develop a PDE-derivative-free second-order optimization strategy based on both gradient and Hessian score-function estimators. For low-fidelity model generation, we propose an on-the-fly strategy that builds and adaptively updates a reduced-order model (ROM) as the variational distribution evolves throughout the optimization. We investigate strategies for adapting the ROM and estimating the approximate control-variate gains, examining tradeoffs between estimator bias, variance reduction, and computational cost. The resulting approach facilitates derivative-free VI for large-scale inverse problems through a reduction in overall computational cost. We demonstrate the method on four numerical examples: a linear benchmark with an exact posterior and three PDE-constrained inverse problems.
\end{abstract}
\begin{keyword}
variational inference \sep multifidelity methods \sep approximate control variates \sep score-function estimators \sep reduced-order models \sep Bayesian inverse problems
\end{keyword}
\end{frontmatter}

\section{Introduction}
The ability to efficiently solve large-scale, PDE-constrained Bayesian inverse problems is an enabling technology across a wide range of disciplines, with applications ranging from aero-thermal dynamics~\cite{RiPaBl26} to ice sheet modeling~\cite{Br22}. Despite major increases in computational power, these problems often remain prohibitively expensive in realistic settings because each likelihood evaluation may require one or more costly PDE solves, and inverse workflows typically demand many such evaluations. The challenge is amplified when the forward model is available only as a black box, as is common for, e.g., legacy codes, multiphysics models that chain together different codebases, and codes with non-differentiable components, in which case adjoint or automatic-differentiation formulations are unavailable.

\subsection{Related work}
In this work, we introduce a multifidelity, PDE-gradient-free variational inference approach to efficiently solve Bayesian inverse problems. To place the approach in context, we briefly review the main classes of methods used in these settings. A first distinction is between gradient-based and gradient-free inference. Gradient-based approaches, including Hamiltonian Monte Carlo~\cite{Brooks_2011}, adjoint-based Laplace approximations, Stein variational inference~\cite{liu2019steinvariationalgradientdescent}, and variational inference with reparameterization-based gradients~\cite{kingma2022autoencodingvariationalbayes}, can be highly effective in terms of both sample efficiency and scalability to high-dimensional parameter spaces. In many practical PDE-constrained settings, however, gradient information is unavailable or unreliable for the reasons mentioned in the previous paragraph. %
Gradient-free approaches provide an appealing alternative in these cases because they require only forward-model evaluations and therefore naturally accommodate black-box simulators, noisy objectives, and complex physics couplings. A primary tradeoff is computational cost: gradient-free methods typically require many more high-fidelity model evaluations than their gradient-based counterparts, and their performance often degrades as the parameter dimension grows. For large-scale PDE models, this simulation burden is a central obstacle.

Among gradient-free Bayesian methods, Markov chain Monte Carlo (MCMC) remains a gold standard because it targets the full posterior distribution. In large-scale PDE-constrained inverse problems, however, gradient-free (and gradient-based) MCMC can become impractical because of long burn-in periods, large sample requirements, and limited opportunities for parallelism. While reduced-order models and surrogates have been incorporated into MCMC in a number of ways~\cite{ChFo05,GaFiWi09,LiWiGh10,CaMaYo15,WaZa05}, these approaches can introduce bias and may still remain too expensive for demanding applications.
Ensemble Kalman methods provide another important class of popular approximate inference strategies for large-scale inverse problems. Emerging from data assimilation~\cite{evensen2009data,kalnay2003atmospheric}, ensemble Kalman filters and related variants propagate an ensemble of forward solves and update parameter estimates through empirical covariance information~\cite{chada2023review}. These methods are derivative-free, naturally parallel, and often effective at scale, but they do not directly target the exact posterior except under restrictive assumptions, and their performance can depend strongly on ensemble size, covariance regularization, and the quality of linear-Gaussian approximations. Relevant to the present work, a number of multilevel, multifidelity, and reduced-order variants of ensemble Kalman methods have been developed to improve computational efficiency~\cite{hoel2016multilevel,yan2019adaptive,popov2021multifidelity,gao2021bi,DoYa22}.

\subsubsection{Variational inference}
The focus of the present work is variational inference (VI). Emerging from the machine learning community, VI~\cite{Jo98,WaJo08} replaces posterior sampling with an optimization problem over a tractable family of distributions. Although this replacement introduces bias, it offers an attractive tradeoff between modeling fidelity and computational cost for large-scale, real-world inverse problems, where physically meaningful priors can be difficult to specify, observational noise is often poorly characterized, model-form errors are prevalent, and uncertainty quantification remains an important objective. While VI is now widely used in machine learning because of its efficiency and scalability~\cite{BlKuMc17}, its use for PDE-constrained inference remains an active and emerging research topic~\cite{rei2025blackboxvariationalinference,ZhNaZh21,AgKo24,ZhChZh24}.

Standard VI algorithms maximize the evidence lower bound (ELBO), typically using stochastic gradient-based optimization. These approaches typically employ pathwise, or reparameterization-based, estimators, which differentiate through the sampled parameters and the log-joint density~\cite{kingma2022autoencodingvariationalbayes,MoRoFiMn20}; in PDE-constrained inference, this generally requires forward-model derivatives. Here, however, we consider the setting in which PDE gradients are unavailable. In this case, the likelihood-ratio, or score-function, identity~\cite{Wi92,Fu06,Gl90} can be used to rewrite gradients of expectations using derivatives of the variational density alone. This idea underlies the original black-box variational inference (BBVI) framework~\cite{Ranganath2014BBVI}. %

The main drawback of score-function estimators is their notoriously high variance~\cite{kingma2022autoencodingvariationalbayes,paisley2012variationalbayesianinferencestochastic,schulman2016gradientestimationusingstochastic,rezende2014stochasticbackpropagationapproximateinference}. Without effective variance reduction, the resulting optimization can be slow, noisy, and unstable. This issue has motivated substantial work on variance reduction for VI to reduce estimator variance. In the initial BBVI reference, Ranganath et al.~\cite{Ranganath2014BBVI} proposed several control-variate strategies~\cite{rubinstein1985efficiency} for variance reduction. Many subsequent works have pursued related ideas for both BBVI and reparameterized VI~\cite{dayta2024variancecontrolblackbox,ng2024pathwisegradientvariancereduction,LiLiJi18,RiBoNuRuAk20}, and a general review of Monte Carlo gradient estimators and variance reduction techniques is provided in~\cite{MoRoFiMn20}.

In addition to variance reduction, various other approaches have been proposed for improving the performance of BBVI. Recently, Rei et al.~\cite{rei2025blackboxvariationalinference} proposed a score-function BBVI approach with adaptive importance-sampling reuse of previously evaluated samples. Their approach reduces repeated forward-model evaluations through sample reuse, which effectively enables larger sample sizes and more efficient inference.
Multilevel ideas have also been introduced in VI to improve performance. For instance, Fujisawa and Sato~\cite{FuSa21} use multilevel Monte Carlo to reuse reparameterized gradient information across successive optimization steps, while Li et al.~\cite{LiYaChHeSu25} combine multilevel control functionals with reparameterized VI gradients. In both cases, the levels are defined by successive optimization steps, and the methods rely on reparameterized gradients; in the context of PDE-constrained inference, we again emphasize that reparameterized gradients require gradients of the forward model. More closely related, Agrawal et al.~\cite{AgRaKo25} use multifidelity score-function estimators for stochastic optimization of design variables. Their focus, however, is optimization rather than Bayesian inference: they do not optimize an ELBO over a variational posterior or construct a low-fidelity model on the fly for PDE-constrained inference.

While not directly related, we additionally note that recent work has also examined derivative-free VI methods for large-scale Bayesian inverse problems that avoid stochastic score-function estimators. Che et al.~\cite{ChChHuHuWa25} formulate a derivative-free VI method for Gaussian-mixture variational families based on Fisher--Rao natural-gradient dynamics. The approach replaces the Monte Carlo estimators commonly used in black-box VI with specialized deterministic quadrature rules. Related work has incorporated reduced-order models directly within VI; Wu et al.~\cite{ZhChZh24}, for example, developed a ROM-based normalizing-flow VI approach for Bayesian elliptic inverse problems.

\subsubsection{Optimization}
Beyond estimator construction, one must also decide how to solve the resulting VI optimization problem. The basic approach is stochastic gradient ascent using Euclidean gradients, but this can be sensitive to scaling and variance collapse~\cite{Wierstra2008NES,Wierstra2014NES}. Natural-gradient methods~\cite{Amari1998NaturalGradient} address this through the use of the Fisher information of the variational family to define a metric on distributions. Score-function natural-gradient VI has been developed in several works~\cite{TrCoMe18,JiSuSu21,MaShIo26}. Relevant to the present work, Magris et al.~\cite{MaShIo26} propose a BBVI approach that exploits exponential-family dual coordinates to avoid explicitly forming and inverting the Fisher matrix. While the Fisher matrix generally differs from the Hessian of the ELBO, the update proposed in Magris et al. contains a likelihood-weighted score kernel and, as a result, incorporates curvature of the objective. %

Various works have directly considered the use of second-order methods, which compute (potentially an approximation to) the Hessian; the majority of these works require derivatives of the forward model~\cite{FaWaBeKwHe15,ReMc16,ReJoMc17,CoZhLiChCa21}. Related approaches have explored randomized quasi-Monte Carlo sampling with stochastic L-BFGS~\cite{LiOw21}. To the best of our knowledge, however, PDE-derivative-free second-order score-function estimators --- which directly approximate the Hessian using model evaluations and derivatives of the sampling density alone --- have not been explored within BBVI. However, these estimators have a long history outside of VI~\cite{Ru89,Gl04,FoFaAlRoXiWh18,FuLeBa16,MaMoPrBh23}. Collectively, these works have derived score-function Hessian identities, the potentially high variance of such an estimator, and established their use within optimization; however, these works have not targeted the setting of variational inference.
\subsection{Our contributions}
We develop a multifidelity black-box variational inference (MF-BBVI)
framework for PDE-constrained Bayesian inverse problems that combines
multifidelity score-function estimation, regularized second-order
optimization, and adaptive reduced-order modeling.
The proposed framework enables efficient VI without requiring derivatives of the underlying
PDE model. The main contributions are:

\begin{enumerate}[leftmargin=*]
    \item A multifidelity approach for estimating both the gradient and
    Hessian of the high-fidelity ELBO using score-function identities
    and approximate control variates (ACV). By combining a small number
    of high-fidelity evaluations with a larger number of inexpensive
    low-fidelity evaluations, the approach can reduce estimator variance, with unbiasedness depending on the sampling, gain-estimation, and ROM-adaptation strategies. To the best of our knowledge, this work is the first to employ multifidelity approximate control variates for the joint estimation of first- and second-order derivatives within VI. 

    \item A regularized stochastic Newton-type optimization strategy
    that leverages noisy estimates of the ELBO Hessian to construct
    efficient updates. The approach combines
    diagonal curvature approximations, lagged Hessian averaging,
    and regularization. %

    \item An on-the-fly reduced-order modeling strategy that
    adaptively constructs and updates a reduced-order model as the
    variational distribution evolves, enabling efficient integration
    of reduced-order models within the multifidelity ACV framework.

    \item An assessment of the proposed MF-BBVI approach on four inverse problems, ranging from a linear benchmark with an analytic posterior to nonlinear PDE-constrained inverse problems. We compare MF-BBVI to single-fidelity BBVI and examine the effects of sample allocation, ACV gain structure and estimation, and ROM adaptation on convergence, robustness, and high-fidelity evaluation cost. These experiments demonstrate the efficacy of MF-BBVI and the proposed second-order Newton-type optimization scheme. 

\end{enumerate}

The remainder of the paper is organized as follows. Section~\ref{sec:vi} reviews VI, including black-box VI, control variates, and optimization approaches. Section~\ref{sec:mf_vi} introduces our proposed multifidelity black-box VI formulation, including the approximate-control-variate estimators and the second-order optimization strategy. Section~\ref{sec:otf_rom} outlines the proposed on-the-fly adaptive reduced-order modeling approach and further discusses several approaches for computing the ACV gain and adapting the ROM. Section~\ref{sec:results} presents numerical experiments, while Section~\ref{sec:conclusions} provides conclusions and perspectives.

\section{Variational inference}
\label{sec:vi}
We consider the black-box PDE-constrained Bayesian inverse problem: given observations $\obs \in \RR{m}$, we aim to characterize the posterior distribution
\begin{equation}
\label{eq:posterior}
\begin{aligned}
\mathrm{find} & \; p(\params\mid \obs) \propto p(\obs\mid\params)p_0(\params),\\
\mathrm{subject} \; \mathrm{to} & \; \residual(\state;\params)=\mathbf{0},
\end{aligned}
\end{equation}
where $p(\params \mid \obs)$ is the posterior, $p(\obs \mid \params)$ is the likelihood, and $p_0(\params)$ is the prior. Additionally, $\params \in \paramsSpace \subseteq \RR{\xDim}$ denotes unknown, uncertain parameters, $\state \in \RR{M}$ denotes the (discretized) PDE state, $\fomQoi = \stateToObservableMap(\state,\params) \in \RR{m}$ denotes the observables obtained from the PDE simulation under the state-to-observable map $\stateToObservableMap$ such that $\fomQoi \approx \obs$ under appropriate $\params$, 
 and $\residual : \RR{M} \times \paramsSpace \rightarrow \RR{M}$ is the discretized PDE residual. We consider both steady and transient systems; in the latter setting the state and residual should be viewed in the space--time sense. In what follows, we focus on the setting where (i) evaluating the likelihood, which requires solving the forward PDE, is expensive and (ii) PDE gradients are unavailable.

The present work considers VI-based approaches for generating approximate solutions to the Bayesian inference problem~\eqref{eq:posterior}. 
VI seeks approximate solutions to~\eqref{eq:posterior} by approximating the posterior distribution $p(\params \mid \obs)$ with a \textit{search distribution} $\searchDistribution\left(\params;\searchDistributionParameters\right)$, where $\searchDistributionParameters$ denotes the variational parameters (e.g., mean and covariance in the Gaussian case).
VI then minimizes the Kullback-Leibler (KL) divergence,
\begin{equation}
\label{eq:vi_kl}
\searchDistributionParameters^* = \underset{\searchDistributionParameters}{\arg \, \min}  \; \mathrm{KL}\!\left(\searchDistribution(\params ; \searchDistributionParameters)\,\|\,p(\params\mid \obs)\right),
\end{equation}
where
\begin{equation}
\label{eq:kl_vi}
\mathrm{KL}\!\left(
\searchDistribution(\params;\searchDistributionParameters)
\,\|\,
p(\params \mid \obs)
\right)
=
\int
\searchDistribution(\params;\searchDistributionParameters)
\log\!\left(
\frac{
\searchDistribution(\params;\searchDistributionParameters)
}{
p(\params \mid \obs)
}
\right)
\, d\params,
\end{equation}
is the KL divergence. 
Minimizing the KL divergence is equivalent to maximizing the evidence lower bound (ELBO),
\begin{equation}
\label{eq:vi_elbo}
\searchDistributionParameters^* = \underset{\searchDistributionParameters}{\arg \, \max}\;  \mathbb{E}_{\params \sim \searchDistribution\left(\params;\searchDistributionParameters\right)} \left[ \ELBO\!\left(\params;\searchDistributionParameters \right) \right],
\end{equation}
where we refer to
\begin{equation}
\label{eq:elbo}
\ELBO(\params;\searchDistributionParameters)
:=\log p(\obs,\params)-\log \searchDistribution(\params;\searchDistributionParameters)
\end{equation}
as the ELBO integrand.

\subsection{Optimization approach}
\label{sec:vi_optimization}
A standard approach for maximizing the ELBO~\eqref{eq:vi_elbo} is gradient-based optimization. In a gradient-based optimization approach, one differentiates the ELBO with respect to the search distribution parameters,
\begin{equation}\label{eq:vi_gradient}
\nabla_{\searchDistributionParameters}  \; \;  \mathbb{E}_{\params \sim \searchDistribution\left(\params;\searchDistributionParameters\right)} \left[ \ELBO\left( \params; \searchDistributionParameters \right) \right].
\end{equation}
This gradient is used to compute a step. For example, the naive gradient ascent step would be
\begin{equation}\label{eq:elbo_update_gradient}
\searchDistributionParameters \leftarrow \searchDistributionParameters + \epsilon \nabla_{\searchDistributionParameters}   \mathbb{E}_{\params \sim \searchDistribution\left(\params;\searchDistributionParameters\right)} \left[ \ELBO\left( \params; \searchDistributionParameters \right) \right],
\end{equation}
where $\epsilon$ is the step size.

Before proceeding, we note that in practice, standard gradient ascent with a fixed step size can fail. To describe why, following Refs.~\cite{Wierstra2008NES,Wierstra2014NES}, we consider the canonical case of a univariate Gaussian variational family, $\searchDistributionParameters  = \left(\mu,\sigma\right)$, with $\mu$ being the mean and $\sigma$ the standard deviation. The gradients of the variational distribution with respect to these parameters are
\[
 \nabla_{\mu}  \log \searchDistribution \left( \params ; \searchDistributionParameters \right) =  \frac{\params - \mu}{\sigma^2},
\qquad
\nabla_{\sigma}  \log \searchDistribution \left( \params ; \searchDistributionParameters \right) =  \frac{ \left(\params - \mu\right)^2  - \sigma^2}{\sigma^3}.
\]
In a first-order update, the lack of scale-invariance of the gradient results in inefficient or unstable steps. In particular, if $\sigma$ is large the steps will be tiny, and as $\sigma \rightarrow 0$ the gradient explodes. Reparameterizing in terms of $\log\sigma$, as is done in Section~\ref{sec:mean_field}, removes the $\sigma^{-3}$ singularity in the second expression, but does not by itself remove the dependence of the step on the scale of $\sigma$.

Various mitigations addressing this issue have been proposed. For example, Ref.~\cite{Ranganath2014BBVI} proposes the use of AdaGrad, which utilizes adaptive step sizes. Another standard mitigation is the use of natural gradients~\cite{Amari1998NaturalGradient}. The key idea of the natural gradient is to remove dependence on the parameterization by measuring distance using the KL divergence between distributions. The natural gradient of the ELBO can be written as
\[
\widetilde{\nabla}_{\searchDistributionParameters}\,
\mathbb{E}_{\params \sim \searchDistribution(\params;\searchDistributionParameters)}
\left[ \ELBO(\params;\searchDistributionParameters) \right]
=
\fisherMatrix^{-1}\nabla_{\searchDistributionParameters}\,
\mathbb{E}_{\params \sim \searchDistribution(\params;\searchDistributionParameters)}
\left[ \ELBO(\params;\searchDistributionParameters) \right],
\]
where
\begin{equation}\label{eq:fisher_generic}
\fisherMatrix(\searchDistributionParameters)
:=
\mathbb{E}_{\params\sim q(\cdot;\searchDistributionParameters)}
\left[
\nabla_{\searchDistributionParameters}\log q(\params;\searchDistributionParameters)\nabla_{\searchDistributionParameters}\log q(\params;\searchDistributionParameters)^{\top}
\right]
\end{equation}
is the Fisher information matrix. Because natural gradients measure distance in terms of the KL divergence between distributions, they are invariant to reparameterization to first order (i.e., in the small-step limit), and they are widely used in variational inference and related optimization methods. The primary drawback of a natural gradient is that it requires the Fisher information matrix, which is not analytically tractable for all distributions.

\subsection{Black-box VI}
Evaluation of the update~\eqref{eq:elbo_update_gradient} requires computation of the gradient of the ELBO. Here, we summarize ``black-box'' VI estimators that do not require gradients of the underlying PDE. To this end, we note that the ELBO gradient~\eqref{eq:vi_gradient} is usually computed in one of two ways: through the use of the ``reparameterization'' trick~\cite{kingma2022autoencodingvariationalbayes} or through the use of ``score-function'' estimators~\cite{Wi92,Fu06,Gl90}. The reparameterization trick typically results in lower-variance estimators, but requires gradients of the PDE residual. Score-function estimators (also known as REINFORCE or the log-derivative trick) exploit the identity $\nabla_{\searchDistributionParameters}\searchDistribution = \searchDistribution\,\nabla_{\searchDistributionParameters}\log\searchDistribution$ to write~\eqref{eq:vi_gradient} as
\begin{equation}\label{eq:reinforce}
\begin{aligned}
\nabla_{\searchDistributionParameters}\,
\mathbb{E}_{\params \sim \searchDistribution(\params;\searchDistributionParameters)}                          
\left[ \ELBO(\params;\searchDistributionParameters) \right]
&=
\mathbb{E}_{\params \sim \searchDistribution(\params;\searchDistributionParameters)}
\left[ \ELBO(\params;\searchDistributionParameters)\,
\nabla_{\searchDistributionParameters}\log \searchDistribution(\params;\searchDistributionParameters) \right]
+
\mathbb{E}_{\params \sim \searchDistribution(\params;\searchDistributionParameters)}
\left[ \nabla_{\searchDistributionParameters}\ELBO(\params;\searchDistributionParameters) \right]\\
&=
\mathbb{E}_{\params \sim \searchDistribution(\params;\searchDistributionParameters)}
\left[ \ELBO(\params;\searchDistributionParameters)\,
\nabla_{\searchDistributionParameters}\log \searchDistribution(\params;\searchDistributionParameters) \right].
\end{aligned}
\end{equation}
The second term vanishes since
$\nabla_{\searchDistributionParameters}\ELBO(\params;\searchDistributionParameters)
=-\nabla_{\searchDistributionParameters}\log\searchDistribution(\params;\searchDistributionParameters)$
and the score has zero mean,
$\mathbb{E}_{\params \sim \searchDistribution(\params;\searchDistributionParameters)}
\left[\nabla_{\searchDistributionParameters}\log\searchDistribution(\params;\searchDistributionParameters)\right]=\mathbf 0$.
In practice, the expectation in~\eqref{eq:reinforce} is evaluated by Monte Carlo (MC) sampling.

Given i.i.d.\ samples $\params^{(i)} \sim \searchDistribution\left(\cdot;\searchDistributionParameters\right)$ for $i=1,\dots,\nSamples$, an (unbiased) MC estimator of the gradient is 
\begin{equation}\label{eq:reinforce_mc}
\monteCarloGradient
=
\left\langle
\ELBO(\params;\searchDistributionParameters)\nabla_{\searchDistributionParameters}\log \searchDistribution(\params;\searchDistributionParameters)
\right\rangle_{\{\params^{(i)}\}_{i=1}^{\nSamples}},
\end{equation}
where we use the notation
\[
\langle r\rangle_{\mathcal{S}}:=\frac{1}{|\mathcal{S}|}\sum_{\params\in\mathcal{S}} r(\params)
\]
to denote the average across the finite sample set \(\mathcal{S}\).

Critically, evaluating the gradient via~\eqref{eq:reinforce}--\eqref{eq:reinforce_mc} only requires gradients of the search distribution (which are often available analytically). No gradients of the PDE are required, and these estimators are thus amenable for PDE-gradient-free optimization. 

\subsection{Mean-field variational family}
\label{sec:mean_field}
The so-called ``mean-field'' variational family is a simple and commonly used choice for the variational distribution $q$. In the mean-field approximation, the variational distribution $\searchDistribution$ is assumed to factorize into independent componentwise distributions. A common choice is to represent each component with a Gaussian distribution, 
\begin{equation}
q(\params;\searchDistributionParameters)
=\mathcal{N}\!\left(
\variationalMean,\operatorname{diag}(\variationalStdScalar_1^2,\ldots,\variationalStdScalar_{\nParams}^2)
\right).
\end{equation}
To enforce a positive variance, it is common to parameterize each standard deviation by its logarithm,
\(\eta_i=\log\sigma_i\), and optimize
\(\searchDistributionParameters =(\variationalMean,\boldsymbol \eta)\). %
The gradients $\nabla_{\searchDistributionParameters} \log \searchDistribution(\params;\searchDistributionParameters )$ required by the score-function estimator are available analytically and are given by
\begin{equation}
\nabla_{\mu_i}\log q(\params;\searchDistributionParameters)
=
\frac{\theta_i-\mu_i}{\sigma_i^2},
\qquad
\nabla_{\eta_i}\log q(\params;\searchDistributionParameters)
=
\frac{(\theta_i-\mu_i)^2}{\sigma_i^2}-1.
\end{equation}

\subsection{Variance reduction through control variates}\label{sec:control_variates}
Unfortunately, score-function MC estimators of the gradient~\eqref{eq:reinforce_mc} have a notoriously high variance~\cite{kingma2022autoencodingvariationalbayes,paisley2012variationalbayesianinferencestochastic,schulman2016gradientestimationusingstochastic,rezende2014stochasticbackpropagationapproximateinference}, which results in inaccurate or unstable optimization unless a large number of samples are employed. As sampling $\elbo$ requires solving the high-fidelity model, the computational cost of the required number of samples for low-variance estimators is typically prohibitively expensive, and variance reduction schemes are essential for a practical technique. %
\begin{remark}[Variance of Monte Carlo estimators]
Define
\[
\elboGradientIntegrand
:=
\ELBO(\params;\searchDistributionParameters)\,\nabla_{\searchDistributionParameters}\log q(\params;\searchDistributionParameters),
\]
with $\params\sim q(\cdot;\searchDistributionParameters)$. Since~\eqref{eq:reinforce_mc} involves sample means of i.i.d.\ draws, the covariance of the mean estimator scales as~\cite{wasserman2004all}
\begin{equation}\label{eq:mc_variance_scaling}
\mathrm{Cov}\!\left[\monteCarloGradient\right]
=
\frac{1}{\nSamples}\,\mathrm{Cov}\!\left[\elboGradientIntegrand \right].
\end{equation}
Since the estimator is unbiased, the trace of this covariance is its mean-squared error, which decreases at the standard Monte Carlo rate $\mathcal{O}(\nSamples^{-1})$. Note that, while this rate is independent of the parameter dimension $d$, the magnitude of the covariance $\mathrm{Cov}\!\left[ \elboGradientIntegrand \right]$ of the $2d$-dimensional gradient integrand can display a strong dependence on dimension.
\end{remark}

A standard approach for reducing the variance of Monte Carlo estimators is the use of control variates~\cite{Ranganath2014BBVI,rubinstein1985efficiency}. The basic idea is to augment a target random variable by a correlated auxiliary quantity whose expectation is known. To describe the idea, let \(X\in\mathbb{R}^{d}\) denote a random vector of interest and let \(Y\in\mathbb{R}^{m}\) denote a control variate with known mean \(\mu_Y=\mathbb{E}[Y]\). Then, for any matrix \( \mathbf{A} \in\mathbb{R}^{d\times m}\),
\begin{equation}\label{eq:generic_control_variate}
X_{\mathrm{cv}}(\mathbf{A}):=X+\mathbf{A}\left(Y-\mu_Y\right)
\end{equation}
is unbiased, since
\[
\mathbb{E}\!\left[X_{\mathrm{cv}}(\mathbf{A})\right]=\mathbb{E}[X].
\]
Its covariance is
\begin{equation}\label{eq:generic_control_variate_covariance}
\mathrm{Cov}\!\left[X_{\mathrm{cv}}(\mathbf{A})\right]
=
\mathrm{Cov}[X]
+\mathbf{A}\,\mathrm{Cov}[Y]\,\mathbf{A}^{\top}
+\mathbf{A}\,\mathrm{Cov}[Y,X]
+\mathrm{Cov}[X,Y]\,\mathbf{A}^{\top}.
\end{equation}
When \(Y\) is strongly correlated with \(X\), an appropriate choice of \(\mathbf{A}\) can significantly reduce the estimator variance without introducing bias.
\begin{remark}[Optimal control-variate gains]
Provided \(\mathrm{Cov}[Y]\) is nonsingular, minimizing the trace of the covariance~\eqref{eq:generic_control_variate_covariance} gives~\cite{rubinstein1985efficiency,popov2021multifidelity} (see~\ref{sec:appendix_optimal_gains} for the derivation)
\begin{equation}\label{eq:generic_control_variate_optimal_matrix}
\mathbf{A}^{\star}=-\,\mathrm{Cov}[X,Y]\mathrm{Cov}[Y]^{-1}.
\end{equation}
Alternatively, when \(m=d\), one can restrict \(\mathbf{A}\) to be diagonal and solve for its entries by a componentwise minimization of the corresponding scalar variances, which results in
\begin{equation}\label{eq:compwise}
\mathbf{A}^{\star}_{\mathrm{diag}}
=
\mathrm{diag}(a_1^{\star},\ldots,a_d^{\star}),
\qquad
a_i^{\star}
=
-\,\frac{\mathrm{Cov}[X_i,Y_i]}{\mathrm{Var}[Y_i]},
\end{equation}
where \(X_i\) and \(Y_i\) denote the \(i\)th components of \(X\) and \(Y\), respectively.
\end{remark}

In black-box VI, Ranganath et al.~\cite{Ranganath2014BBVI} propose the use of the score function of the variational approximation as a natural zero-mean control variate. Using the notation introduced above, we define
\[
    \controlVariate(\params;\searchDistributionParameters):=\nabla_{\searchDistributionParameters}\log q(\params;\searchDistributionParameters),
\]
so that \(\mathbb{E}[\controlVariate]=\mathbf{0}\). This leads to the unbiased control-variate estimator
\begin{equation}\label{eq:cv_gradient_estimator}
\controlVariateGradient
(\mathbf{A}_g)
:=
\left\langle
    \elboGradientIntegrand(\params;\searchDistributionParameters)+\mathbf{A}_g \controlVariate(\params;\searchDistributionParameters)
\right\rangle_{\{\params^{(i)}\}_{i=1}^{\nSamples}},
\end{equation}
for a suitable, fixed coefficient matrix \(\mathbf{A}_g\). Ref.~\cite{Ranganath2014BBVI} employs the componentwise formulation given in Eq.~\eqref{eq:compwise}, wherein the model covariances are estimated from the same Monte Carlo samples used to compute the score functions; note that this estimation technique can result in a biased estimator, as will be discussed in Section~\ref{sec:acv_sampling}.  %

\section{Multifidelity black-box variational inference}
\label{sec:mf_vi}
While control variates can be effective, black-box VI still suffers from high-variance estimators and slow convergence. In practice, these issues result in high computational costs that can make black-box VI impractical for large-scale applications. To address these issues, here we propose a novel \textit{multifidelity} black-box VI approach utilizing approximate control-variate estimators, regularized stochastic second-order updates, and on-the-fly adaptive surrogates.

\subsection{Approximate control variates}
The control-variate formulation described in Section~\ref{sec:control_variates} assumes access to an auxiliary random variable with known expectation. In multifidelity settings, the control variate takes the form of a low-cost, lower-fidelity model. Typically, the expectation of the low-fidelity quantity is not available in closed form. A standard workaround is to replace the exact control variate by an \emph{approximate} control variate (ACV)~\cite{GORODETSKY2020109257}: one uses a correlated low-fidelity model whose expectation can be estimated from additional low-fidelity samples. The resulting estimator remains unbiased provided an appropriate sampling scheme is used (e.g., the correction is written so that the same low-fidelity quantity appears both with and without the high-fidelity sample subset).

To describe the idea, let \(X\) again denote a high-fidelity random quantity and let \(Y\) again denote a correlated low-fidelity quantity. If \(\mathbb{E}[Y]\) were known, the additive control-variate estimator would take the form
\[
X_{\mathrm{cv}}(\mathbf{A})=X+\mathbf{A}\left(Y-\mathbb{E}[Y]\right).
\]
In the multifidelity setting, \(\mathbb{E}[Y]\) is replaced by a cheap Monte Carlo estimate formed from a larger low-fidelity sample set, and the resulting two-model estimator, which coincides with the multifidelity Monte Carlo (MFMC) estimator of~\cite{PeWiGu16} and is a special case of the ACV framework of~\cite{GORODETSKY2020109257}, takes the form
\begin{equation}\label{eq:approximate_control_variate_generic}
\widehat X_{\mathrm{acv}}(\mathbf{A})
:=
\left\langle X\right\rangle_{\mathcal{S}^{\mathrm{HF}}}
+
\mathbf{A}
\left(
\left\langle Y\right\rangle_{\mathcal{S}^{\mathrm{HF}}}
-\left\langle Y\right\rangle_{\mathcal{S}^{\mathrm{LF}}}
\right),
\end{equation}
where $\sampleSetHighFidelity$ and $\sampleSetLowFidelity$ with $\sampleSetHighFidelity \subseteq \sampleSetLowFidelity$ are nested high-fidelity and low-fidelity sample sets. 
Both low-fidelity averages have the same expectation, and hence for any fixed $\mathbf{A}$ the multifidelity estimator~\eqref{eq:approximate_control_variate_generic} is unbiased, i.e.,
\[\mathbb{E}[\widehat X_{\mathrm{acv}}(\mathbf{A})] = \mathbb{E}[X].\]
The same result holds for any random $\mathbf{A}$, provided its construction is independent of the samples used to construct the ACV correction; this point will be revisited later in this work.

Relative to standard Monte Carlo, the benefit of the multifidelity estimator~\eqref{eq:approximate_control_variate_generic} is that one can take \(|\sampleSetLowFidelity|\gg |\sampleSetHighFidelity|\) because the low-fidelity model is much cheaper to evaluate. When the LF model is correlated with the HF model, the ACV estimator can substantially reduce the estimator variance.

\begin{remark}[Covariance and practical choices for the control-variate matrix]
Define
\[
\mathbf{\Sigma}_{XX}:=\mathrm{Cov}(X),\qquad
\mathbf{\Sigma}_{YY}:=\mathrm{Cov}(Y),\qquad
\mathbf{\Sigma}_{XY}:=\mathrm{Cov}(X,Y),\qquad
\mathbf{\Sigma}_{YX}:=\mathbf{\Sigma}_{XY}^{\top},
\]
and let \(\nSamplesHF:=|\mathcal{S}^{\mathrm{HF}}|\), \(\nSamplesLF:=|\mathcal{S}^{\mathrm{LF}}|\), with \(\mathcal{S}^{\mathrm{HF}}\subseteq\mathcal{S}^{\mathrm{LF}}\). Then the covariance of the generic approximate-control-variate estimator~\eqref{eq:approximate_control_variate_generic} is
\begin{equation}\label{eq:mf_covariance}
\mathrm{Cov}\!\left[\widehat X_{\mathrm{acv}}(\mathbf{A})\right]
=
\frac{1}{\nSamplesHF}\mathbf{\Sigma}_{XX}
+
\left(\frac{1}{\nSamplesHF}-\frac{1}{\nSamplesLF}\right)
\left(
\mathbf{A}\mathbf{\Sigma}_{YY}\mathbf{A}^{\top}
+\mathbf{A}\mathbf{\Sigma}_{YX}
+\mathbf{\Sigma}_{XY}\mathbf{A}^{\top}
\right).
\end{equation}
The gain matrix $\mathbf{A}$ is chosen to minimize the estimator variance. In practice, the required covariances are estimated from samples, and small HF batches can compromise the reliability of the fitted gain. We therefore consider two structured formulations—a componentwise gain and a scalar gain—and refer the reader to~\ref{sec:appendix_optimal_gains} for their derivations. An unrestricted matrix gain could additionally exploit correlations across components, but its estimation requires more coefficients and can involve a singular or poorly conditioned empirical covariance matrix. We consequently restrict the present work to scalar and componentwise gains.
\begin{itemize}[leftmargin=*]

\item \textit{Diagonal (component-wise) control variate.}
Restricting \(\mathbf{A}=\mathbf{D}=\mathrm{diag}(a_1,\ldots,a_d)\) and minimizing the trace of the covariance, the optimal gains are
\[
a_i^\star=-\,\frac{(\mathbf{\Sigma}_{XY})_{ii}}{(\mathbf{\Sigma}_{YY})_{ii}},\qquad i=1,\ldots,d.
\]
The total variance is
\begin{equation}\label{eq:mf_diag_trace_opt}
\mathrm{tr}\!\left(\mathrm{Cov}[\widehat X_{\mathrm{acv}}(\mathbf{D}^\star)]\right)
=
\frac{1}{\nSamplesHF}\mathrm{tr}(\mathbf{\Sigma}_{XX})
-
\left(\frac{1}{\nSamplesHF}-\frac{1}{\nSamplesLF}\right)
\sum_{i=1}^{d}\frac{(\mathbf{\Sigma}_{XY})_{ii}^2}{(\mathbf{\Sigma}_{YY})_{ii}}.
\end{equation}

\item \textit{Scalar control variate.}
Lastly, restricting \(\mathbf{A}=\alpha \mathbf{I}\) and minimizing the trace of the covariance, the optimal scalar gain parameter is
\[
\alpha^\star=-\,\frac{\mathrm{tr}(\mathbf{\Sigma}_{XY})}{\mathrm{tr}(\mathbf{\Sigma}_{YY})},
\]
with
\begin{equation}\label{eq:mf_scalar_trace_opt}
\mathrm{tr}\!\left(\mathrm{Cov}[\widehat X_{\mathrm{acv}}(\alpha^\star \mathbf{I})]\right)
=
\frac{1}{\nSamplesHF}\mathrm{tr}(\mathbf{\Sigma}_{XX})
-
\left(\frac{1}{\nSamplesHF}-\frac{1}{\nSamplesLF}\right)
\frac{\mathrm{tr}(\mathbf{\Sigma}_{XY})^2}{\mathrm{tr}(\mathbf{\Sigma}_{YY})}.
\end{equation}
\end{itemize}
\end{remark}

\subsubsection{Multifidelity VI estimators}
With the above discussion in mind, we now describe our multifidelity estimators. For simplicity, we describe our proposed approach in the context of a two-model hierarchy: one low-fidelity model and one high-fidelity model. We emphasize, however, that the approach can be extended to handle an arbitrary number of low-fidelity models (see, e.g.,~\cite{GORODETSKY2020109257} for this in the context of UQ). %

At iteration \(k\), we draw \(\nSamplesLF\ge\nSamplesHF\) independent parameter samples \(\params^{(i)}\sim q(\cdot;\searchDistributionParameters_k)\) and define the nested sets
\[
\sampleSetHighFidelity_k
=
\{\params^{(i)}\}_{i=1}^{\nSamplesHF},
\qquad
\sampleSetLowFidelity_k
=
\{\params^{(i)}\}_{i=1}^{\nSamplesLF}.
\]
By construction, \(\sampleSetHighFidelity_k\subseteq\sampleSetLowFidelity_k\).
Both models are evaluated at the first \(\nSamplesHF\) parameter samples;
the remaining samples are evaluated only with the low-fidelity model and,
since all samples are drawn independently, are independent of the shared samples.
Let \(p_{\mathrm{HF}}(\obs,\params)\) and
\(p_{\mathrm{LF}}(\obs,\params)\) denote the high- and low-fidelity joint
densities, respectively. These differ only through the likelihood: writing
\(p_{\mathrm{HF}}(\obs,\params)=p(\obs\mid\fomQoi(\params))\,p_0(\params)\), the low-fidelity joint density is
\(p_{\mathrm{LF}}(\obs,\params)=p(\obs\mid\rom(\params))\,p_0(\params)\), where \(\rom\) is the low-fidelity
approximation to the parameter-to-observable map (Section~\ref{sec:otf_rom}). The prior, the
\(-\log q\) term, and (for bounded parameters) the Jacobian term of Section~\ref{sec:bounds} are
shared by both fidelities, so the HF--LF difference in the ELBO integrands below is
purely the difference in log-likelihoods.
We define the high- and low-fidelity ELBO integrands as
    \begin{equation}\label{eq:joint_elbo_integrands}
\begin{aligned}
\ELBOHF(\params;\searchDistributionParameters)
&:=\log p_{\mathrm{HF}}(\obs,\params)-\log q(\params;\searchDistributionParameters),\\
\ELBOLF(\params;\searchDistributionParameters)
&:=\log p_{\mathrm{LF}}(\obs,\params)-\log q(\params;\searchDistributionParameters).
\end{aligned}
\end{equation}
The corresponding high- and low-fidelity score-function gradient contributions are
\[
\elboGradientIntegrandHF (\params;\searchDistributionParameters)
:=\ELBOHF(\params;\searchDistributionParameters)\nabla_{\searchDistributionParameters}\log q(\params;\searchDistributionParameters),
\qquad
\elboGradientIntegrandLF(\params;\searchDistributionParameters)
:=\ELBOLF(\params;\searchDistributionParameters)\nabla_{\searchDistributionParameters}\log q(\params;\searchDistributionParameters).
\]
The corresponding ELBO gradients are their expectations,
\begin{equation}\label{eq:joint_entropy_gradient_expectations}
\begin{aligned}
    \nabla_{\searchDistributionParameters}\mathbb{E}_{\params\sim q(\cdot;\searchDistributionParameters)}
\left[\ELBOHF(\params;\searchDistributionParameters)\right]
&=\mathbb{E}_{\params\sim q(\cdot;\searchDistributionParameters)}
\left[\elboGradientIntegrandHF(\params;\searchDistributionParameters)\right],\\
    \nabla_{\searchDistributionParameters}\mathbb{E}_{\params\sim q(\cdot;\searchDistributionParameters)}
\left[\ELBOLF(\params;\searchDistributionParameters)\right]
&=\mathbb{E}_{\params\sim q(\cdot;\searchDistributionParameters)}
\left[\elboGradientIntegrandLF(\params;\searchDistributionParameters)\right].
\end{aligned}
\end{equation}
A multifidelity estimator for the high-fidelity ELBO gradient at iteration $k$ is thus
\begin{equation}\label{eq:mfmc_gradient_joint_entropy}
\approximateControlVariateGradient_k
:=
\left\langle \elboGradientIntegrandHF(\params;\searchDistributionParameters_k)\right\rangle_{\mathcal{S}^{\mathrm{HF}}_k}
+
\acvGain_k
\left(
\left\langle \elboGradientIntegrandLF(\params;\searchDistributionParameters_k)\right\rangle_{\mathcal{S}^{\mathrm{HF}}_k}
-\left\langle \elboGradientIntegrandLF(\params;\searchDistributionParameters_k)\right\rangle_{\mathcal{S}^{\mathrm{LF}}_k}
\right),
\end{equation}
where $\acvGain_k$ is the ACV gain at iteration $k$, with the same
HF-minus-LF correction convention as Eq.~\eqref{eq:approximate_control_variate_generic}. 
Under the restriction that the gain matrix $\acvGain_k$ and LF model are fixed, the multifidelity approximation is unbiased and can have a substantially lower variance
than its single-fidelity counterpart and therefore requires fewer
high-fidelity samples for a target accuracy.

\subsection{Regularized stochastic diagonal-Newton optimization}
Once the variance of the score-function gradient has been reduced, it becomes practical to exploit second-order information. Here, we propose the use of the same PDE-gradient-free machinery used for the gradient estimator to construct an estimator for the Hessian such that we can perform a Newton-type update. In a Newton-type update, the gradient is preconditioned by the inverse of the Hessian and the step is
\begin{equation}\label{eq:elbo_update_gradient_second_order}
    \searchDistributionParameters \leftarrow \searchDistributionParameters - \epsilon \left[ \nabla_{\searchDistributionParameters}^2\mathbb{E}_{q_{\searchDistributionParameters}}[\ELBO(\params;\searchDistributionParameters)] \right]^{-1}  \nabla_{\searchDistributionParameters}   \mathbb{E}_{\params \sim \searchDistribution\left(\params;\searchDistributionParameters\right)} \left[ \ELBO\left( \params; \searchDistributionParameters \right) \right],
\end{equation}
where $0 < \epsilon \le 1$ is a step size. Equation~\eqref{eq:elbo_update_gradient_second_order} is the idealized Newton update; the update used in practice, Eq.~\eqref{eq:mfmc_newton}, replaces the exact Hessian by a lagged, regularized diagonal estimate.
The remainder of this subsection describes the use of multifidelity estimators for computing the Hessian and the proposed optimization approach.

\subsubsection{Hessian estimator}
To incorporate second-order information into the inference process, we twice apply the score-function (log-derivative) identity to obtain
\begin{equation}\label{eq:reinforce_second_order}
\nabla_{\searchDistributionParameters}^2\mathbb{E}_{q_{\searchDistributionParameters}}[\ELBO(\params;\searchDistributionParameters)]
=
\mathbb{E}_{\params\sim q(\cdot;\searchDistributionParameters)}\!\left[
\log p(\obs,\params)\left(\nabla_{\searchDistributionParameters}\log q \, \nabla_{\searchDistributionParameters}\log q^{\top}
+\nabla_{\searchDistributionParameters}^2\log q\right)
\right]
+\nabla_{\searchDistributionParameters}^2\mathcal{H}(q_{\searchDistributionParameters}),
\end{equation}
where
\(\mathcal{H}(q_{\searchDistributionParameters}):=-\mathbb{E}_{q_{\searchDistributionParameters}}[\log q(\params;\searchDistributionParameters)]\)
is the entropy of the variational distribution.
Equation~\eqref{eq:reinforce_second_order} follows by applying the second-order score identity to \(\log p(\obs,\params)\), which does not depend explicitly on \(\searchDistributionParameters\).
In writing~\eqref{eq:reinforce_second_order}, we note that we have separated out the entropy term; for the mean-field Gaussian family parameterized in \((\variationalMean,\log\variationalStd)\), it is straightforward to show $\nabla_{\searchDistributionParameters}^2\mathcal{H}(q_{\searchDistributionParameters}) = 0$. We note that a score-function-type construction could equivalently be formulated for variational families whose entropy contribution is not analytically tractable.

As with the gradient estimator in~\eqref{eq:reinforce}, evaluating~\eqref{eq:reinforce_second_order} only requires first- and second-order derivatives of the search distribution with respect to the variational parameters. No derivatives of the underlying PDE model are required. To this end, we define the high and low-fidelity Hessian score-function quantities as
\[
\begin{aligned}
\elboHessianIntegrandHF(\params;\searchDistributionParameters)
&:=\log p_{\mathrm{HF}}(\obs,\params)\left(
\nabla_{\searchDistributionParameters}\log q(\params;\searchDistributionParameters)\nabla_{\searchDistributionParameters}\log q(\params;\searchDistributionParameters)^\top
+\nabla_{\searchDistributionParameters}^2\log q(\params;\searchDistributionParameters)
\right),\\
\elboHessianIntegrandLF(\params;\searchDistributionParameters)
&:=\log p_{\mathrm{LF}}(\obs,\params)\left(
\nabla_{\searchDistributionParameters}\log q(\params;\searchDistributionParameters)\nabla_{\searchDistributionParameters}\log q(\params;\searchDistributionParameters)^\top
+\nabla_{\searchDistributionParameters}^2\log q(\params;\searchDistributionParameters)
\right).
\end{aligned}
\]

\subsubsection{Lagged Hessian estimator}
Using the same parameter samples as the gradient estimator, the multifidelity Hessian estimate at iteration $k$ is
\begin{equation}\label{eq:mfmc_hessian}
{\approximateControlVariateHessian}_k
:=
\left\langle \elboHessianIntegrandHF(\params;\searchDistributionParameters_k)\right\rangle_{\mathcal{S}^{\mathrm{HF}}_k}
+
\hessianAcvGain{k}
\left(
\left\langle \elboHessianIntegrandLF(\params;\searchDistributionParameters_k)\right\rangle_{\mathcal{S}^{\mathrm{HF}}_k}
-\left\langle \elboHessianIntegrandLF(\params;\searchDistributionParameters_k)\right\rangle_{\mathcal{S}^{\mathrm{LF}}_k}
\right)
+\nabla_{\searchDistributionParameters}^2\mathcal{H}(q_{\searchDistributionParameters})\big|_{\searchDistributionParameters=\searchDistributionParameters_k}.
\end{equation}
Here \(\hessianAcvGain{k}\) is the Hessian ACV gain, the counterpart of the
gradient gain \(\acvGain_k\) in Eq.~\eqref{eq:mfmc_gradient_joint_entropy}.
Because the Hessian correction is a matrix, the gain acts on its entries:
for the scalar gains used in this work, \(\hessianAcvGain{k}\) multiplies
the mean and log-standard-deviation blocks of the correction by separate
scalar coefficients (see Eq.~\eqref{eq:scalar_hessian_gain_blocks}).

As in the gradient case, the two low-fidelity averages have the same
expectation so long as the LF model and gain are fixed. As a result, the estimator~\eqref{eq:mfmc_hessian} is unbiased for the
high-fidelity Hessian. %

Using this same-sample Hessian estimate directly to precondition the gradient, however, does not preserve a zero expected step at a fixed point.
Specifically,
$
\mathbb E[{\approximateControlVariateHessian}_k^{-1} \approximateControlVariateGradient_k]
\ne
\mathbb E[{\approximateControlVariateHessian}_k^{-1}]\,
\mathbb E[\approximateControlVariateGradient_k]
$
since ${\approximateControlVariateHessian}_k$ and $\approximateControlVariateGradient_k$ are not independent.
As a result, a zero expected gradient \(\mathbb E[ \approximateControlVariateGradient_k]=0\) does not imply that the expected
update vanishes.

A robust and rigorous solution to this issue is to use independent sample sets for the Hessian and gradient estimators. This approach, while mathematically correct, incurs a
substantial computational cost as it involves evaluation of an entirely new sample set. As a practical alternative, here we instead propose the use of a \textit{lagged} Hessian update that makes use
of an exponential moving average of preceding Hessian estimates. Let
\(\laggedHessian_k\) denote the curvature used at iteration \(k\), with
\begin{equation}\label{eq:lagged_hessian_average}
\laggedHessian_k=\beta\laggedHessian_{k-1}+(1-\beta){\approximateControlVariateHessian}_{k-1},
\qquad k\ge1,\quad 0\le\beta<1.
\end{equation}
We initialize \(\laggedHessian_0={\approximateControlVariateHessian}_0\),
so the first step uses same-sample curvature. Conceptually, the lagged update is
\begin{equation}\label{eq:elbo_update_hessian_lagged}
    \searchDistributionParameters_{k+1}=\searchDistributionParameters_k-\epsilon_k\laggedHessian_k^{-1} \approximateControlVariateGradient_k.
\end{equation}
In practice, $\laggedHessian_k$ is further regularized before use, as described in Section~\ref{sec:fisher_regularization}.
After initialization, the current estimate
\({\approximateControlVariateHessian}_k\) does not enter the curvature
used for the current direction; it contributes to
\(\laggedHessian_{k+1}\). As a result, at a given iteration,\footnote{The expectation is understood to be taken given a fixed $\searchDistributionParameters_k$, not over the optimization trajectory.} 
\begin{equation}\label{eq:lagged_hessian_conditional_factorization}
\mathbb E[\laggedHessian_k^{-1}\approximateControlVariateGradient_k
]
=\mathbb E[\laggedHessian_k^{-1}]
\mathbb E[\approximateControlVariateGradient_k
]
\end{equation}
Thus, if $\mathbb{E}[\approximateControlVariateGradient_k]$ vanishes, so too does  $\mathbb{E}[\laggedHessian_k^{-1}]
\mathbb E[\approximateControlVariateGradient_k
]$. 

\subsubsection{Fisher-scaled Hessian regularization}
\label{sec:fisher_regularization}
Monte Carlo noise can produce inaccurate or even wrong-sign curvature in the Hessian estimates, which can in turn degrade the performance of the optimizer. 
To address this, we employ a regularized diagonal Hessian estimate. Specifically, we retain only the diagonal entries of the Hessian to reduce sensitivity to noisy
cross-coordinate estimates, and further regularize  the diagonal using the Fisher metric.
For the mean-field Gaussian family with variational parameters
\(\searchDistributionParameters=[\variationalMean^\top,
\variationalLogStd^\top]^\top\), where
\(\variationalLogStd=\log\variationalStd\), the Fisher information matrix is
\begin{equation}\label{eq:natural_metric_scaling}
\fisherMatrix_k=
\begin{bmatrix}
\operatorname{diag}(\variationalStd_k^{-2})&\mathbf{0}\\
\mathbf{0}&2\mathbf{I}
\end{bmatrix}.
\end{equation}
At iteration \(k\), we regularize the lagged curvature
\(\laggedHessian_k\) using the Fisher matrix \(\fisherMatrix_k\)
of the current search distribution. Define
\begin{equation}\label{eq:lagged_hessian_diagonal}
h_{k,i}:=(\laggedHessian_k)_{ii},
\qquad i=1,\ldots,2d.
\end{equation}
Note that, since we maximize the
ELBO, suitable curvature is negative. We thus construct a diagonal
regularized Hessian with entries
\begin{equation}\label{eq:absolute_hessian}
\regularizedHessian_{k,ii}=
\begin{cases}
h_{k,i}-\additiveRegularization\fisherMatrix_{k,ii},
&-h_{k,i}/\fisherMatrix_{k,ii}\ge\tau,\\
-(\fallbackLearningRate^{-1}+\additiveRegularization)\fisherMatrix_{k,ii},
&-h_{k,i}/\fisherMatrix_{k,ii}<\tau,
\end{cases}
\end{equation}
where \(\tau>0\) is a curvature tolerance (we use \(\tau=1 \times 10^{-9}\)), \(\fallbackLearningRate>0\) is the fallback
learning rate, and \(\additiveRegularization\ge0\) is the hyperparameter controlling additive damping.
All off-diagonal entries of \(\regularizedHessian_k\) are zero. This rule
retains negative curvature and replaces near-zero  
or wrong-sign curvature with a natural-gradient update: the fallback step in Eq.~\eqref{eq:mfmc_newton} is \((\Delta\searchDistributionParameters_k)_i=(\fallbackLearningRate^{-1}+\additiveRegularization)^{-1}(\fisherMatrix_k^{-1}\approximateControlVariateGradient_k)_i\), i.e., a natural-gradient step with learning rate approximately \(\fallbackLearningRate\) when \(\additiveRegularization\ll\fallbackLearningRate^{-1}\).\footnote{We note that equation~\eqref{eq:absolute_hessian} is equivalent to scaling by
\(\fisherMatrix_k^{-1/2}\), replacing scaled negative-Hessian entries strictly
below \(\tau\) by \(\fallbackLearningRate^{-1}\), adding \(\additiveRegularization\mathbf I\), and mapping
back to the original coordinates.}  
The resulting update is
\begin{equation}\label{eq:mfmc_newton}
\begin{aligned}
\Delta\searchDistributionParameters_k
&=-\regularizedHessian_k^{-1}\approximateControlVariateGradient_k,\\
\searchDistributionParameters_{k+1}
&=\searchDistributionParameters_k
+\epsilon_k\Delta\searchDistributionParameters_k.
\end{aligned}
\end{equation}
In the present work, we employ a fixed unit step size, $\epsilon_k=1.0$. 

\subsection{Additional ingredients}
\label{sec:additional_ingredients}
Lastly, we supplement the multifidelity estimators with a leave-one-out baseline for additional variance reduction and a transformation
for bounded parameters.

\subsubsection{Leave-one-out baseline for variance reduction}
A leave-one-out (LOO) baseline~\cite{kool2019buy} reduces variance in 
the score-function estimators by subtracting the average ELBO integrand over
the other samples in the same batch. We have found LOO baselines to be effective for small sample sizes, and as such we employ them here.

For fidelity
\(r\in\{\mathrm{HF},\mathrm{LF}\}\), the sample set \(\sampleSet\) with $| \sampleSet | > 1 $,
and \(\params^{(i)}\in\sampleSet\), we define the baseline
\begin{equation}\label{eq:loo_baseline}
\looGradientBaseline{r}^{(-i)}(\sampleSet)
:=\frac{1}{|\sampleSet|-1}
\sum_{\substack{\params^{(j)}\in\sampleSet\\j\ne i}}
\ELBO_r(\params^{(j)};\searchDistributionParameters_k).
\end{equation}
The centered gradient terms are
\[
\looElboGradientIntegrand{r}(\params^{(i)};\searchDistributionParameters_k)
:=\elboGradientIntegrand_{r}(\params^{(i)};\searchDistributionParameters_k)
-\looGradientBaseline{r}^{(-i)}(\sampleSet)\,
\nabla_{\searchDistributionParameters}\log\searchDistribution(\params^{(i)};\searchDistributionParameters_k).
\]
For HF evaluations, \(\sampleSet=\sampleSetHighFidelity_k\).
For LF evaluations, we center separately within the sample sets
\(\sampleSetHighFidelity_k\) and
\(\sampleSetLowFidelity_k\setminus\sampleSetHighFidelity_k\), i.e., each LF evaluation uses
\(\sampleSet=\sampleSetHighFidelity_k\) if \(\params^{(i)}\in\sampleSetHighFidelity_k\) and
\(\sampleSet=\sampleSetLowFidelity_k\setminus\sampleSetHighFidelity_k\) otherwise.
Centering the shared LF evaluations with the same baseline set as the HF evaluations
keeps the paired HF and LF terms on \(\sampleSetHighFidelity_k\) centered consistently, which
preserves their correlation. In Eq.~\eqref{eq:loo_mfmc_gradient}, the average over
\(\sampleSetLowFidelity_k\) is taken over these centered terms, each with its own baseline.
The multifidelity gradient estimator in
Eq.~\eqref{eq:mfmc_gradient_joint_entropy} with a LOO baseline becomes
\begin{equation}\label{eq:loo_mfmc_gradient}
\approximateControlVariateGradient_k
=
\left\langle\looElboGradientIntegrand{\mathrm{HF}}\right\rangle_{\sampleSetHighFidelity_k}
+\acvGain_k
\left(
\left\langle\looElboGradientIntegrand{\mathrm{LF}}\right\rangle_{\sampleSetHighFidelity_k}
-\left\langle\looElboGradientIntegrand{\mathrm{LF}}\right\rangle_{\sampleSetLowFidelity_k}
\right).
\end{equation}
For independent parameter samples,
the LOO baseline is independent of the excluded sample. Since the score term (i.e., 
$\nabla_{\searchDistributionParameters}\log\searchDistribution(\params^{(i)};\searchDistributionParameters)$)
has zero expectation, centering preserves the expectation of each
score-function term and hence of the multifidelity estimator for a fixed
gain and a fixed LF model. 

For the Hessian estimator, the baseline
\(\looHessianBaseline{r}^{(-i)}(\sampleSet)\) is defined by
Eq.~\eqref{eq:loo_baseline} with \(\ELBO_r\) replaced by
\(\log p_r(\obs,\params)\). Writing
\(\score{\params}{\searchDistributionParameters}:=\nabla_{\searchDistributionParameters}\log\searchDistribution(\params;\searchDistributionParameters)\)
for the score and using the second-order score
\[
\scoreH{\params}{\searchDistributionParameters}
:=\score{\params}{\searchDistributionParameters}
\score{\params}{\searchDistributionParameters}^{\top}
+\nabla_{\searchDistributionParameters}^2
\log\searchDistribution(\params;\searchDistributionParameters),
\]
which also has zero expectation under \(\searchDistribution(\cdot;\searchDistributionParameters)\) (the information identity), the centered terms are
\[
\looElboHessianIntegrand{r}(\params^{(i)};\searchDistributionParameters_k)
:=\elboHessianIntegrand_{r}(\params^{(i)};\searchDistributionParameters_k)
-\looHessianBaseline{r}^{(-i)}(\sampleSet)\,\scoreH{\params^{(i)}}{\searchDistributionParameters_k}.
\]
The corresponding LOO-centered multifidelity Hessian estimator is
\begin{equation}\label{eq:loo_mfmc_hessian}
{\approximateControlVariateHessian}_k
=
\left\langle\looElboHessianIntegrand{\mathrm{HF}}\right\rangle_{\sampleSetHighFidelity_k}
+\hessianAcvGain{k}
\left(
\left\langle\looElboHessianIntegrand{\mathrm{LF}}\right\rangle_{\sampleSetHighFidelity_k}
-\left\langle\looElboHessianIntegrand{\mathrm{LF}}\right\rangle_{\sampleSetLowFidelity_k}
\right)
+\nabla_{\searchDistributionParameters}^2\mathcal H(\searchDistribution_{\searchDistributionParameters}).
\end{equation}
This retains the Hessian decomposition in
Eq.~\eqref{eq:reinforce_second_order}, where entropy is differentiated
analytically. The centered gradient and Hessian estimates are used in
the Newton update~\eqref{eq:mfmc_newton} and the lagged Hessian
average~\eqref{eq:lagged_hessian_average}, respectively.

\subsubsection{Variable bounds}
\label{sec:bounds}
Physical parameters often have meaningful bounds. To enforce finite bounds
\(a_i<b_i\), we map physical coordinates \(\params\) to unconstrained,
transformed coordinates \(\params^*=\boldsymbol g(\params)\). The map
\(\boldsymbol g\) is the component-wise inverse of a sigmoid:
\begin{equation}\label{eq:vi_sigmoid_bounds}
\theta_i^*=g_i(\theta_i)
=\log\!\left(\frac{u_i}{1-u_i}\right),\qquad
u_i=\frac{\theta_i-a_i}{b_i-a_i},\qquad i=1,\ldots,\nParams.
\end{equation}
Equivalently, we map from transformed coordinates to physical space as
\(\theta_i=g_i^{-1}(\theta_i^*)
=a_i+(b_i-a_i)/(1+e^{-\theta_i^*})\).

The Jacobian of the physical-to-transformed map is the diagonal matrix
\begin{equation}\label{eq:vi_bound_jacobian}
\boldsymbol J_{\boldsymbol g}(\params)
:=\frac{\partial\boldsymbol g(\params)}{\partial\params^\top}
=\operatorname{diag}\!\left(g_i'(\theta_i)\right)_{i=1}^{\nParams},\qquad
g_i'(\theta_i)
=\frac{1}{(b_i-a_i)u_i(1-u_i)}.
\end{equation}
For problems with bounds, we take the mean-field variational distribution to be Gaussian in transformed
coordinates, with mean \(\boldsymbol\mu^*\) and diagonal covariance
\(\boldsymbol\Sigma^*\), and denote its density by
\(\searchDistribution^*(\params^*;\searchDistributionParameters)\).
The joint density must be expressed in the same coordinates. For
\(r\in\{\mathrm{HF},\mathrm{LF}\}\), the change-of-variables formula gives
\begin{equation}\label{eq:vi_transformed_joint}
p_r^*(\obs,\params^*)
=p_r\!\left(\obs,\boldsymbol g^{-1}(\params^*)\right)
\left|\det\boldsymbol J_{\boldsymbol g^{-1}}(\params^*)\right|,
\end{equation}
where
\(\boldsymbol J_{\boldsymbol g^{-1}}(\params^*)
=\boldsymbol J_{\boldsymbol g}(\boldsymbol g^{-1}(\params^*))^{-1}\).
Thus, the ELBO in transformed coordinates is the expectation under
\(\searchDistribution^*\) of
\begin{equation}\label{eq:vi_transformed_elbo_integrand}
\ELBO_r^*(\params^*;\searchDistributionParameters)
={}\log p_r\!\left(\obs,\boldsymbol g^{-1}(\params^*)\right)
+\log\left|\det\boldsymbol J_{\boldsymbol g^{-1}}(\params^*)\right|
-\log\searchDistribution^*(\params^*;\searchDistributionParameters).
\end{equation}
The likelihood and prior are evaluated at the mapped physical parameters,
while the variational density and entropy are
evaluated in transformed coordinates. The same Jacobian contribution is
included in the HF and LF log-joint terms used by the gradient and Hessian
estimators.

Given physical-coordinate initialization settings
\(\boldsymbol\mu_{\mathrm{init}}\) and diagonal
\(\boldsymbol\Sigma_{\mathrm{init}}\), we initialize the transformed
Gaussian using the local linearization
\begin{equation}\label{eq:vi_initial_bound_covariance}
\begin{aligned}
\boldsymbol\mu_0^*&=\boldsymbol g(\boldsymbol\mu_{\mathrm{init}}),\\
\boldsymbol\Sigma_0^*&=
\boldsymbol J_{\boldsymbol g}(\boldsymbol\mu_{\mathrm{init}})
\boldsymbol\Sigma_{\mathrm{init}}
\boldsymbol J_{\boldsymbol g}(\boldsymbol\mu_{\mathrm{init}})^\top.
\end{aligned}
\end{equation}
We then optimize \(\boldsymbol\mu^*\) and the logarithms
of the transformed standard deviations.
At iteration \(k\), the physical mean and covariance are
\begin{equation}\label{eq:vi_forward_bound_covariance}
\begin{aligned}
\boldsymbol\mu^{\mathrm{phys}}_k
&:=\mathbb E_{\params^*\sim\searchDistribution^*(\cdot;\searchDistributionParameters_k)}
\!\left[\boldsymbol g^{-1}(\params^*)\right],\\
\boldsymbol\Sigma^{\mathrm{phys}}_k
&:=\operatorname{Cov}_{\params^*\sim\searchDistribution^*(\cdot;\searchDistributionParameters_k)}
\!\left[\boldsymbol g^{-1}(\params^*)\right].
\end{aligned}
\end{equation}
We note that the mapped center \(\boldsymbol g^{-1}(\boldsymbol\mu_k^*)\)
generally differs from \(\boldsymbol\mu^{\mathrm{phys}}_k\) because the map is nonlinear.

\section{On-the-fly adaptive low-fidelity modeling}
\label{sec:otf_rom}
The preceding multifidelity estimators and optimization methods apply to
any low-fidelity model that approximates the parameter-to-observable map. Its construction may be independent of
the inference process, for example through a coarse discretization,
a simplified physical model, or an offline surrogate. The computational
benefit depends on its cost and on the correlation between its
score-function terms and their HF counterparts.

Here, we construct a data-driven reduced-order model (ROM) from HF evaluations generated during inference. Throughout, we use the term ROM broadly to denote any reduced, inexpensive approximation of the parameter-to-observable map, including projection-based ROMs and non-intrusive regression models such as the POD--GP model of Section~\ref{sec:otf_rom_gp}. This avoids a separate
offline training stage and concentrates training data in the region
explored by the current variational distribution. As that distribution
evolves, the ROM is updated to maintain local accuracy. The approach
is similar to the adaptive multifidelity strategy of Donoghue and
Yano~\cite{DoYa22} for ensemble Kalman filtering. We first describe the
training and adaptation procedure, then discuss the consequences of reusing HF samples for training and
correlation estimation, and lastly specify an automated POD--GP ROM formulation.

\subsection{Training and adaptation from high-fidelity samples}
Let \(\fomQoi(\params)\in\mathbb R^m\) denote the HF QoI and let
\(\rom(\params;\romParams)\) be a trainable approximation with model parameters $\romParams$. Neural
networks, Gaussian processes, and projection-based ROMs provide possible
constructions. At iteration \(k\), the shared parameter samples
\(\params_k^{(i)}\in\sampleSetHighFidelity_k\), with
\(\sampleSetHighFidelity_k=\{\params_k^{(i)}\}_{i=1}^{\nSamplesHF}\),
yield the training pairs
\[
\mathcal D_k
=\{(\params_k^{(i)},\fomQoi(\params_k^{(i)}))\}_{i=1}^{\nSamplesHF}.
\]
Projection-based ROMs may additionally require solution snapshots and,
when hyper-reduction is used, data for the nonlinear operators. The training 
data at each iteration $k$ can be collected and used to fit the ROM.
To limit training cost and retain data relevant to the current search
distribution, we use a rolling window of the most recent \(n_{\mathrm{hist}}\) HF sample sets. Specifically, defining the training dataset at iteration $k$ as $\mathcal{T}_k$, we employ
\[
\mathcal T_k
:=\bigcup_{j=\max(0,k-n_{\mathrm{hist}}+1)}^k\mathcal D_j,
\qquad |\mathcal T_k|\le n_{\mathrm{hist}}\nSamplesHF.
\]
The ROM is fitted to these pairs using a procedure appropriate to
its construction; a least-squares fit, for example, minimizes
\[
\underset{\romParams}{\operatorname{minimize}}\;
\sum_{(\params,\fomQoi)\in\mathcal T_k}
\|\rom(\params;\romParams)-\fomQoi\|_2^2.
\]

\begin{figure}
\begin{center}
\begin{minipage}[t]{0.46\textwidth}
\centering
\scriptsize\parbox[t]{0.97\linewidth}{\centering Step 1\\Draw \(\mathcal{S}^{\mathrm{HF}}_k\) and sample HF model}
\begin{tikzpicture}[x=0.76cm,y=0.76cm]
  \path[use as bounding box] (0,0) rectangle (7.2,4.2);
  \draw[thick] (0.35,0.3) rectangle (6.85,3.7);
  \draw[blue!70!black,thick] (3.3,2.0) ellipse (2.0 and 1.15);
  \node[blue!50!black] at (1.25,3.05) {\scriptsize \(q(\cdot;\zeta_k)\)};
  \foreach \x/\y in {2.1/2.55,2.7/1.55,3.4/2.15,4.0/1.75}
    \fill[red!80!black] (\x,\y) circle (1.8pt);
  \fill[red!80!black] (5.0,0.85) circle (1.8pt);
  \node[anchor=west] at (5.05,0.95) {\scriptsize HF};
  \fill[green!50!black] (5.0,0.5) circle (1.5pt);
  \node[anchor=west] at (5.05,0.5) {\scriptsize LF};
\end{tikzpicture}
\end{minipage}
\hspace{0.004\textwidth}
\begin{minipage}[t]{0.46\textwidth}
\centering
    \scriptsize\parbox[t]{0.97\linewidth}{\centering Step 2\\Potentially adapt ROM, draw additional samples, and evaluate LF model on \(\mathcal{S}^{\mathrm{LF}}_k\).}
\begin{tikzpicture}[x=0.76cm,y=0.76cm]
  \path[use as bounding box] (0,0) rectangle (7.2,4.2);
  \draw[thick] (0.35,0.3) rectangle (6.85,3.7);
  \draw[blue!70!black,thick] (3.3,2.0) ellipse (2.0 and 1.15);
  \node[blue!50!black] at (1.25,3.05) {\scriptsize \(q(\cdot;\zeta_k)\)};
  \foreach \x/\y in {1.7/1.45,2.05/2.6,2.4/1.95,2.9/2.75,3.35/1.45,3.75/2.2,4.15/1.8,4.45/2.55,2.55/1.35,3.85/1.4}
    \fill[green!50!black] (\x,\y) circle (1.5pt);
  \foreach \x/\y in {2.1/2.55,2.7/1.55,3.4/2.15,4.0/1.75}
    \fill[red!80!black] (\x,\y) circle (1.8pt);
\end{tikzpicture}
\end{minipage}

\vspace{0.5em}

\begin{minipage}[t]{0.46\textwidth}
\centering
\scriptsize\parbox[t]{0.97\linewidth}{\centering Step 3\\Form MF gradient/Hessian estimates}
\begin{tikzpicture}[x=0.76cm,y=0.76cm]
  \path[use as bounding box] (0,0) rectangle (7.2,4.2);
  \draw[thick] (0.35,0.3) rectangle (6.85,3.7);
  \draw[blue!70!black,thick] (3.3,2.0) ellipse (2.0 and 1.15);
  \node[blue!50!black] at (1.25,3.05) {\scriptsize \(q(\cdot;\zeta_k)\)};

  \foreach \x/\y in {1.7/1.45,2.05/2.6,2.4/1.95,2.9/2.75,3.35/1.45,3.75/2.2,4.15/1.8,4.45/2.55,2.55/1.35,3.85/1.4}
    \fill[green!50!black] (\x,\y) circle (1.5pt);
  \foreach \x/\y in {2.1/2.55,2.7/1.55,3.4/2.15,4.0/1.75}
    \fill[red!80!black] (\x,\y) circle (1.8pt);
\end{tikzpicture}
\end{minipage}
\hspace{0.004\textwidth}
\begin{minipage}[t]{0.46\textwidth}
\centering
\scriptsize\parbox[t]{0.97\linewidth}{\centering Step 4\\Update to \(q(\cdot;\zeta_{k+1})\) and repeat}
\begin{tikzpicture}[x=0.76cm,y=0.76cm]
  \path[use as bounding box] (0,0) rectangle (7.2,4.2);
  \draw[thick] (0.35,0.3) rectangle (6.85,3.7);
  \draw[blue!70!black,dashed,thick] (3.3,2.0) ellipse (2.0 and 1.15);
  \node[blue!50!black] at (1.25,3.05) {\scriptsize \(q(\cdot;\zeta_k)\)};
  \draw[blue!85!black,thick] (4.8,1.95) ellipse (1.4 and 0.8);
  \node[blue!85!black] at (5.45,2.9) {\scriptsize \(q(\cdot;\zeta_{k+1})\)};
  \draw[->,very thick] (3.4,2.0) -- (4.0,1.95);
  \foreach \x/\y in {4.2/1.45,4.55/2.35,4.9/1.7,5.2/2.45,5.55/1.55,4.65/2.7}
    \fill[green!50!black] (\x,\y) circle (1.5pt);
  \foreach \x/\y in {4.55/1.95,5.05/2.05,4.85/2.45}
    \fill[red!80!black] (\x,\y) circle (1.8pt);
\end{tikzpicture}
\end{minipage}
\caption{Schematic of the on-the-fly multifidelity VI loop used to construct and refresh the low-fidelity ROM. Reading left-to-right and then top-to-bottom: draw the high-fidelity sample set \(\mathcal{S}^{\mathrm{HF}}_k\), potentially adapt the ROM and sample at the low-fidelity sample set \(\mathcal{S}^{\mathrm{LF}}_k\), form the multifidelity score-function estimators, and update the search distribution to \(q(\cdot;\zeta_{k+1})\). Red points denote the shared samples evaluated with both models; green points denote additional LF-only samples.}
\label{fig:otf_mfvi_sequence}
\end{center}
\end{figure}
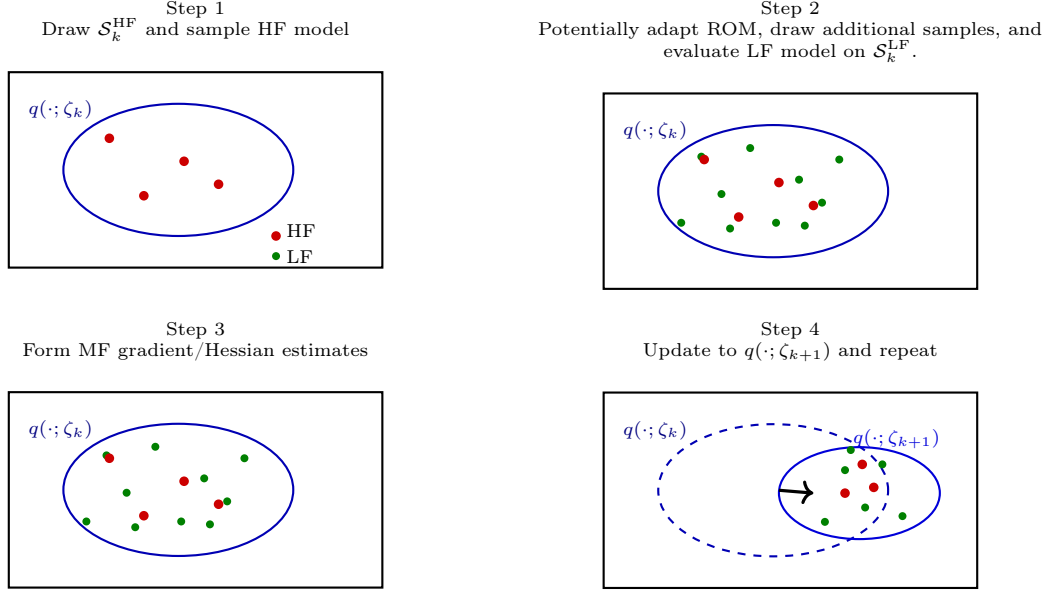

At each iteration, additional HF data becomes available, and as such the ROM could be adapted.  
Here, we test the existing LF model against the new HF
responses before deciding whether to adapt the LF model. Writing
\(\rom_k(\params):=\rom(\params;\romParams_k)\) for the ROM retained
at iteration \(k\), we test the model retained from iteration
\(k-1\) on the current HF sample set \(\sampleSetHighFidelity_k\),
before deciding whether to retrain:
\[
\romError{k}
:=\frac{\displaystyle\sum_{\params\in\sampleSetHighFidelity_k}
\|\rom_{k-1}(\params)-\fomQoi(\params)\|_2}
{\displaystyle\sum_{\params\in\sampleSetHighFidelity_k}
\|\fomQoi(\params)\|_2+\romErrorFloor},
\qquad k>0.
\]
The index \(k-1\) identifies the preceding retained model, which may
have been trained at an earlier iteration, while the test samples belong to
iteration \(k\).
Here \(\romErrorFloor>0\) prevents division by zero. If
\(\romError{k}>\tau_{\mathrm{ROM}}\), we retrain on \(\mathcal T_k\);
otherwise, \(\rom_k=\rom_{k-1}\). An initial ROM is trained when
none is available. Figure~\ref{fig:otf_mfvi_sequence} summarizes this procedure.

\subsection{Adaptation and estimation of correlations}
\label{sec:adaptation_correlations}
The gradient and Hessian gains \(\acvGain_k\) and
\(\hessianAcvGain{k}\) depend on correlations between the HF
and LF score-function quantities. These correlations vary with the search
distribution and can change as the ROM is adapted, and must therefore be estimated
online. We emphasize that estimating the gain requires paired HF--LF evaluations, rather than
LF evaluations alone. 

A natural choice to estimate the gains at iteration $k$ is to utilize the HF samples at iteration $k$. This choice, 
while efficient and simple, introduces several subtle challenges. First, if the ROM is re-trained on the HF data at step $k$,
the paired LF-HF data at step $k$ can result in optimistically biased estimates of the correlation since these data were directly used for training the ROM. 
Second, even if the ROM is not adapted, the use of samples from iteration $k$ to estimate the gain can result in a biased estimator. We now describe this in more detail and provide several mitigations. Table~\ref{tab:adaptive_correlation_strategies} summarizes the six
combinations that will be discussed. 

\begin{table}[!htbp]
\centering
\small
\caption{Removal of bias sources by ROM adaptation and gradient-gain
estimation: \biasCheck\ indicates that the condition is satisfied;
\biasCross\ indicates that bias is not excluded. The classification
assumes fresh independent samples after initialization.}
\label{tab:adaptive_correlation_strategies}
\begin{tabular}{@{}llcc@{}}
\hline
ROM adaptation & Gradient gain
& \shortstack{No ROM-induced\\bias}
& \shortstack{No ACV-gain-induced\\bias}\\
\hline
Same sample & Same sample & \biasCross & \biasCross\\
Same sample & LOO & \biasCross & \biasCheck\\
Same sample & Independent & \biasCross & \biasCheck\\
Delayed & Same sample & \biasCheck & \biasCross\\
Delayed & LOO & \biasCheck & \biasCheck\\
Delayed & Independent & \biasCheck & \biasCheck\\
\hline
\end{tabular}
\end{table}

\subsubsection{Same-sample estimation of ACV gain introduces bias}
\label{sec:acv_sampling}
We first illustrate how estimating the gradient gain using the current
paired HF--LF samples can introduce bias. For brevity, write
\(\elboGradientIntegrand_{r,k}(\params)
:=\elboGradientIntegrand_r(\params;\searchDistributionParameters_k)\),
for \(r\in\{\mathrm{HF},\mathrm{LF}\}\), and let
\[
\overline{\elboGradientIntegrand}_{\mathrm{HF},k}
:=\left\langle\elboGradientIntegrand_{\mathrm{HF},k}\right\rangle_{\sampleSetHighFidelity_k},
\qquad
\gradient_k:=\mathbb E[\elboGradientIntegrand_{\mathrm{HF},k}].
\]
Here \(\gradient_k\) is the HF ELBO gradient in
Eq.~\eqref{eq:joint_entropy_gradient_expectations}. All expectations
below are taken at fixed \(\searchDistributionParameters_k\),
conditional on the previous sampling history. Recall that the
multifidelity gradient estimator can be written as
\[
\approximateControlVariateGradient_k
=\overline{\elboGradientIntegrand}_{\mathrm{HF},k}
+\acvGain_k\gradientAcvCorrection_k,
\]
where, for \(\nSamplesLF>\nSamplesHF>1\),
\[
\gradientAcvCorrection_k
:=\left\langle\elboGradientIntegrand_{\mathrm{LF},k}\right\rangle_{\sampleSetHighFidelity_k}
-\left\langle\elboGradientIntegrand_{\mathrm{LF},k}\right\rangle_{\sampleSetLowFidelity_k}.
\]
The bias therefore satisfies
\begin{equation}\label{eq:adaptive_acv_bias}
\mathbb E[\approximateControlVariateGradient_k]-\gradient_k
=\mathbb E[\acvGain_k]\mathbb E[\gradientAcvCorrection_k]
+\mathbb E\!\left[
\bigl(\acvGain_k-\mathbb E[\acvGain_k]\bigr)
\bigl(\gradientAcvCorrection_k-\mathbb E[\gradientAcvCorrection_k]\bigr)
\right].
\end{equation}
For a fixed ROM, \(\mathbb E[\gradientAcvCorrection_k]=\boldsymbol 0\).
However, fitting \(\acvGain_k\) on the same paired samples used in
\(\gradientAcvCorrection_k\) can make the second term nonzero (since $\acvGain_k$ and $\gradientAcvCorrection_k$ are no longer independent).
Thus, sample sharing can introduce bias even when the ROM is fixed.

We consider three approaches for computing \(\acvGain_k\), each with
scalar and componentwise gains. For the scalar gain, 
we fit a separate scalar gain to the mean block and log standard deviation block 
of the gradient. We further note that, for the
Hessian gain \(\hessianAcvGain{k}\), we employ ``same-sample'' estimation
in all cases, which is described below. This approach results in a biased estimator; we permit bias in the Hessian estimate as it is only used to precondition the gradient. 

For notational simplicity, define the empirical gradient
cross-covariance computed on the set $\sampleSet$ as
\[
\widehat{\boldsymbol\Sigma}_{r,s}(\sampleSet)
:=\Big\langle
\bigl(\elboGradientIntegrand_{r,k}
-\langle\elboGradientIntegrand_{r,k}\rangle_{\sampleSet}\bigr)
\bigl(\elboGradientIntegrand_{s,k}
-\langle\elboGradientIntegrand_{s,k}\rangle_{\sampleSet}\bigr)^\top
\Big\rangle_{\sampleSet}.
\]
When score baselines are used, the covariance is formed from the
corresponding baseline-centered fitting contributions. For the scalar gain, we fit separate coefficients for the mean and
log-standard-deviation blocks. Writing
\[
\elboGradientIntegrand_{r,k}
=\begin{bmatrix}
\elboGradientIntegrand_{r,k,\mu}\\
\elboGradientIntegrand_{r,k,\eta}
\end{bmatrix},
\]
let $\widehat{\boldsymbol\Sigma}_{r,s}^{(b)}(\sampleSet)$ denote the
corresponding $(b,b)$ diagonal block of size $d\times d$, for
$b\in\{\mu,\eta\}$. The two gain
restrictions are
\begin{equation}\label{eq:gradient_gain_families}
\begin{aligned}
\blockGain{b}(\sampleSet)
&:=-\,\frac{\operatorname{tr}\widehat{\boldsymbol\Sigma}_{\mathrm{HF},\mathrm{LF}}^{(b)}(\sampleSet)}
{\operatorname{tr}\widehat{\boldsymbol\Sigma}_{\mathrm{LF},\mathrm{LF}}^{(b)}(\sampleSet)},
\qquad b\in\{\mu,\eta\},\\
\componentGain{j}(\sampleSet)
&:=-\,\frac{\bigl(\widehat{\boldsymbol\Sigma}_{\mathrm{HF},\mathrm{LF}}(\sampleSet)\bigr)_{jj}}
{\bigl(\widehat{\boldsymbol\Sigma}_{\mathrm{LF},\mathrm{LF}}(\sampleSet)\bigr)_{jj}}.
\end{aligned}
\end{equation}
The scalar restriction therefore gives the block-diagonal gain
\begin{equation}\label{eq:scalar_gradient_gain_blocks}
\acvGain^{\mathrm{scalar}}(\sampleSet)
=\operatorname{diag}\!\left(
\blockGain{\mu}(\sampleSet)\mathbf I_d,
\blockGain{\eta}(\sampleSet)\mathbf I_d\right).
\end{equation}
Throughout our numerical experiments, ``scalar gain''
refers to this blockwise restriction rather than one coefficient for the
entire $2d$-component gradient. 

For the diagonal Hessian used in optimization, we fit separate scalar
gains $\blockGain{H,B}(\sampleSet)$ to the vectorized mean--mean and
log-standard-deviation--log-standard-deviation blocks,
$B\in\{\mu\mu,\eta\eta\}$.
The resulting action on the retained diagonal entries of a matrix $\mathbf{Z}$ is
\begin{equation}\label{eq:scalar_hessian_gain_blocks}
    \bigl(\hessianAcvGain{k}(\mathbf{Z})\bigr)_{ii}=
\begin{cases}
\blockGain{H,\mu\mu}Z_{ii},&1\le i\le d,\\
\blockGain{H,\eta\eta}Z_{ii},&d<i\le2d,
\end{cases}
\end{equation}
where the fitting-set argument is suppressed. The resulting matrix is
diagonal. We consider the following approaches for estimating the gains:

\begin{enumerate}
\item \textbf{Same-sample.}
The gains are fitted using all paired evaluations on
\(\sampleSetHighFidelity_k\). Thus, the scalar and componentwise estimators are given by Eq.~\eqref{eq:gradient_gain_families} with $\sampleSet \leftarrow \sampleSetHighFidelity_k$.
This is the simplest approach and uses all available local training pairs without additional model
 evaluations. However, as the gain and correction depend on the same
samples, the second term in Eq.~\eqref{eq:adaptive_acv_bias} does not, in general, vanish, and the estimator may be biased. 

\item \textbf{Leave-one-out.}
For each \(i\), define
\(\sampleSetHighFidelity_{k,-i}
:=\sampleSetHighFidelity_k\setminus\{\params_k^{(i)}\}\).
In an LOO approach, we employ a separate gain for each $i$. The coefficients \(\blockGain{b}(\sampleSetHighFidelity_{k,-i})\), for $b\in\{\mu,\eta\}$,
and \(\componentGain{j}(\sampleSetHighFidelity_{k,-i})\) are fitted without pair
\(i\); the additional LF samples do not enter these fits either. When baselines are used, the baseline-centered contributions used to fit the gain for pair \(i\) are formed within \(\sampleSetHighFidelity_{k,-i}\), so that pair \(i\) does not influence the fitted gain.
Here the baseline-centered contribution is
\[
\elboGradientIntegrand_{r,k}^{(-i)}(\params)
:=\elboGradientIntegrand_{r,k}(\params)
-\looGradientBaseline{r}^{(-i)}(\sampleSetHighFidelity_k)
\score{\params}{\searchDistributionParameters_k}.
\]
The same excluded-pair LF
baseline is applied to the held-out term and the additional LF average.
Defining
\[
\gradientAcvCorrection_{k,i}^{\mathrm{LOO}}
:=\left(1-\frac{\nSamplesHF}{\nSamplesLF}\right)
\left(
\elboGradientIntegrand_{\mathrm{LF},k}^{(-i)}(\params_k^{(i)})
-\left\langle\elboGradientIntegrand_{\mathrm{LF},k}^{(-i)}\right\rangle_{\sampleSetLowFidelity_k\setminus\sampleSetHighFidelity_k}
\right),
\]
the two LOO estimators are
\begin{equation}\label{eq:gain_loo_mfmc_gradient}
\begin{aligned}
\approximateControlVariateGradient_k^{\mathrm{LOO,scalar}}
    &=\frac1{\nSamplesHF}\sum_{i=1}^{\nSamplesHF}\Big[
\elboGradientIntegrand_{\mathrm{HF},k}^{(-i)}(\params_k^{(i)})
+\acvGain^{\mathrm{scalar}}(\sampleSetHighFidelity_{k,-i})\gradientAcvCorrection_{k,i}^{\mathrm{LOO}}\Big],\\
\bigl(\approximateControlVariateGradient_k^{\mathrm{LOO,comp}}\bigr)_j
    &=\frac1{\nSamplesHF}\sum_{i=1}^{\nSamplesHF}\Big[
\bigl(\elboGradientIntegrand_{\mathrm{HF},k}^{(-i)}(\params_k^{(i)})\bigr)_j
+\componentGain{j}(\sampleSetHighFidelity_{k,-i})
\bigl(\gradientAcvCorrection_{k,i}^{\mathrm{LOO}}\bigr)_j\Big],
\end{aligned}
\end{equation}
The LOO approaches require additional gain fits, but no additional model evaluations.

It is straightforward to show that the LOO gain results in an unbiased estimator for a fixed (or delayed) ROM, as each pair
is corrected using a gain fitted without that pair. Holding the ROM and
the other paired samples fixed, the held-out LF contribution and the
additional LF average have the same expected value. Their difference
therefore has zero mean, and multiplying it by the already fitted gain
preserves that zero mean. The same reasoning applies to scalar and componentwise gains. We emphasize, however, that although the LOO gain results in an unbiased estimator, it does not necessarily result 
in a lower variance estimator than the same-sample approach.
\item \textbf{Independent.}
Lastly, we consider an ``independent'' approach where an additional paired sample set
\(\sampleSetGain_k=\{\params_{\mathrm{gain},k}^{(i)}\}_{i=1}^{|\sampleSetGain_k|}\)
is drawn independently from
\(\searchDistribution(\cdot;\searchDistributionParameters_k)\) and
 evaluated with both the HF and LF models. We emphasize that these samples are excluded from ROM training
and from the estimator's HF and LF averages; they are only used to compute the gain. 
The scalar and componentwise estimators are given by Eq.~\eqref{eq:gradient_gain_families} with $\sampleSet \leftarrow \sampleSetGain_k$.
Provided the ROM is fixed (or adapted with a delay), the resulting estimator is unbiased, as the gain is independent of the samples used to construct the MF estimator.
We again emphasize that, although this independent estimator is unbiased, it does not necessarily result in a lower variance estimator.
\end{enumerate}

\subsubsection{Same-sample LF adaptation introduces bias}
Adapting the LF ROM on \(\sampleSetHighFidelity_k\) can also
introduce bias into the gradient estimator, regardless of how the gain
is fitted. Consider a ROM \(\rom_k\) trained on \(\mathcal T_k\supseteq\mathcal D_k\).
Holding this ROM fixed, we define its mean LF gradient contribution as
\[
\boldsymbol\mu_{\mathrm{LF}}(\rom_k)
:=\mathbb E_{\params\sim\searchDistribution(\cdot;\searchDistributionParameters_k)}
[\elboGradientIntegrand_{\mathrm{LF},k}(\params)].
\]
The additional LF samples are independent of the training data and
hence of the fitted ROM. Therefore,
\[
\mathbb E\!\left[
\left\langle\elboGradientIntegrand_{\mathrm{LF},k}\right\rangle_{\sampleSetLowFidelity_k\setminus\sampleSetHighFidelity_k}
\middle|\rom_k\right]=\boldsymbol\mu_{\mathrm{LF}}(\rom_k).
\]
The shared HF samples do not have this independence, since they helped
determine \(\rom_k\). Noting that we can write $\gradientAcvCorrection_k$ in the form 
\[
\gradientAcvCorrection_k
=\left(1-\frac{\nSamplesHF}{\nSamplesLF}\right)
\left(
\left\langle\elboGradientIntegrand_{\mathrm{LF},k}\right\rangle_{\sampleSetHighFidelity_k}
-\left\langle\elboGradientIntegrand_{\mathrm{LF},k}\right\rangle_{\sampleSetLowFidelity_k\setminus\sampleSetHighFidelity_k}
\right),
\]
we have 
\begin{equation}\label{eq:adaptive_acv_correction_mean}
\mathbb E[\gradientAcvCorrection_k\mid\rom_k]
=\left(1-\frac{\nSamplesHF}{\nSamplesLF}\right)
\left(
\mathbb E\!\left[
\left\langle\elboGradientIntegrand_{\mathrm{LF},k}\right\rangle_{\sampleSetHighFidelity_k}
\middle|\rom_k\right]
-\boldsymbol\mu_{\mathrm{LF}}(\rom_k)
\right).
\end{equation}
Taking the expectation again with respect to \(\rom_k\) gives
\begin{equation}\label{eq:adaptive_acv_unconditional_correction_mean}
\begin{aligned}
\mathbb E[\gradientAcvCorrection_k]
&=\mathbb E\!\left[\mathbb E[\gradientAcvCorrection_k\mid\rom_k]\right]\\
&=\left(1-\frac{\nSamplesHF}{\nSamplesLF}\right)
\left(
\mathbb E\!\left[
\left\langle\elboGradientIntegrand_{\mathrm{LF},k}\right\rangle_{\sampleSetHighFidelity_k}
\right]
-\mathbb E[\boldsymbol\mu_{\mathrm{LF}}(\rom_k)]\right).
\end{aligned}
\end{equation}
Thus, the correction has zero mean only if the 
expected mean of the LF model evaluated on $\sampleSetHighFidelity_k$ equals the
expected mean of the LF model evaluated on the fresh samples $\sampleSetLowFidelity_k \setminus \sampleSetHighFidelity_k$. If the LF model is constructed using $\sampleSetHighFidelity_k$, this equality does not hold, in general.

With the above in mind, we consider two approaches for ROM adaptation. To this end, let
\(\romTrainingMap(\mathcal T_k)\) denote the parameters obtained by
fitting the ROM on the retained training window \(\mathcal T_k\). We consider the following approaches for adapting the ROM:
\begin{enumerate}
\item \textbf{Same-sample.} The current HF batch is added to the
retained training window and, when adaptation is required, the ROM is
retrained before constructing the multifidelity estimates. The active ROM at iteration $k$ is thus 
\[
\rom_k(\params)=
\begin{cases}
\rom\bigl(\params;\romTrainingMap(\mathcal T_k)\bigr),
&\romError{k}>\tau_{\mathrm{ROM}},\\
\rom_{k-1}(\params),&\romError{k}\le\tau_{\mathrm{ROM}}.
\end{cases}
\]
Here \(\romError{k}\) is evaluated using the preceding ROM, before
retraining. The LF quantities in the estimators then use \(\rom_k\).
This immediately incorporates the new
responses, but introduces the training dependence in
Eq.~\eqref{eq:adaptive_acv_correction_mean} whenever the current batch
contributes to fitting the active model.

\item \textbf{Delayed adaptation.} The pre-adapted LF model is retained throughout
the current iteration. When adaptation is
        required, the new model is trained but not used until the subsequent iteration. In this case, we evaluate \(\romError{k}\) \textit{and} the LF quantities in the estimators using
the active model \(\rom_k\), and then define the model for the next
iteration as
\[
\rom_{k+1}(\params)=
\begin{cases}
\rom\bigl(\params;\romTrainingMap(\mathcal T_k)\bigr),
&\romError{k}>\tau_{\mathrm{ROM}},\\
\rom_k(\params),&\romError{k}\le\tau_{\mathrm{ROM}}.
\end{cases}
\]
This separates training from the current correction
without additional HF evaluations, at the cost of a delay in responding
to changes in the search distribution.
\end{enumerate}

At initialization, if no ROM is available, both strategies train
\(\rom_0(\params)=\rom(\params;\romTrainingMap(\mathcal T_0))\)
immediately and use the fixed gain $\acvGain_0=-\mathbf I$. The initial iteration therefore always uses a ROM trained on its own HF batch and is subject to the same-sample ROM bias described above, regardless of the adaptation strategy.

\subsubsection{Discussion}
The preceding strategies expose a tradeoff between estimator bias, computational cost, and variance reduction. For instance, same-sample ROM adaptation incorporates the latest HF evaluations immediately. This potentially improves ROM accuracy and HF--LF correlation without additional forward-model evaluations. However, using these samples both to train the ROM and to construct the multifidelity correction can introduce bias. Delayed adaptation removes direct training-sample dependence from fresh correction batches, but the lag in updating the ROM may reduce its effectiveness as the search distribution evolves.
A similar tradeoff arises in estimating the ACV gain. Same-sample fitting uses all available paired evaluations without additional model evaluations, but can introduce dependence between the gain and the correction. With a fixed ROM and independent samples, LOO fitting removes direct dependence between each held-out correction and its fitted gain, but fitting each gain with fewer pairs can increase its variability. Independent fitting separates these quantities explicitly, but comes at the cost of additional HF evaluations. Neither approach guarantees greater variance reduction than same-sample fitting; an unbiased gradient estimator can still have substantial variance when the fitted gains are inaccurate or the HF--LF correlation is weak.
Consequently, the optimal approach depends on the analysis at hand and the desired properties of the inference algorithm. We empirically assess these strategies in our numerical experiments to quantify these tradeoffs. Unless otherwise noted, we employ delayed ROM adaptation with LOO gain estimation as the default configuration, since it yields an unbiased gradient estimator (after initialization) at no additional HF cost. In our experiments, we find that same-sample approaches perform comparably, suggesting that the incurred bias has a minimal impact on the overall optimization problem.

\subsection{Automated model reduction via Gaussian processes}
\label{sec:otf_rom_gp}
The proposed MF-BBVI approach requires the construction of a ROM. The construction of such a ROM, however, 
can be a barrier to entry for use of the algorithm.
To automate ROM construction for low-to-moderate dimensional
parameter spaces and mitigate this issue, we employ a non-intrusive ROM that utilizes proper orthogonal decomposition (POD) of
the QoI responses and Gaussian-process (GP) regression. The construction of the ROM 
requires only parameter--response pairs and can therefore be applied
without access to the governing equations or their derivatives. We briefly describe the construction of the GP.

At a retraining iteration \(k\), we use the training data
\(\mathcal T_k\) defined above and write \(N_k:=|\mathcal T_k|\).
The training inputs, response mean, and centered response matrix are
\[
\begin{aligned}
\boldsymbol\Theta_k
&=[\params]_{(\params,\fomQoi)\in\mathcal T_k},\\
\bar{\fomQoi}_k
&=\frac{1}{N_k}\sum_{(\params,\fomQoi)\in\mathcal T_k}\fomQoi,\\
\mathbf Y_k
&=[\fomQoi-\bar{\fomQoi}_k]_{(\params,\fomQoi)\in\mathcal T_k}
\in\mathbb R^{m\times N_k}.
\end{aligned}
\]
Each column corresponds to one training pair, with the same ordering in
\(\boldsymbol\Theta_k\) and \(\mathbf Y_k\).
To account for differences in parameter scales, we normalize each physical
parameter using its range over the retained training data:
\[
\widetilde\theta_{i,k}(\params)
=\frac{\theta_i-\theta_{i,k}^{\min}}{\Delta\theta_{i,k}},
\qquad
\theta_{i,k}^{\min}
:=\min_{(\params,\fomQoi)\in\mathcal T_k}\theta_i,
\qquad
\Delta\theta_{i,k}
:=\max_{(\params,\fomQoi)\in\mathcal T_k}\theta_i-\theta_{i,k}^{\min}.
\]
We write the normalized input vector as
$\widetilde{\params}_k(\params)$.

From the thin singular value decomposition
\(\mathbf Y_k=\mathbf U_k\podSingularMatrix{k}\mathbf V_k^\top\), we take
\(\mathbf\Phi_k=(\mathbf U_k)_{:,1:\podRank{k}}\). The rank \(\podRank{k}\) is the
smallest positive integer satisfying
\[
\frac{\sum_{j=1}^{\podRank{k}}\podSingularValue{j}^2}
{\sum_{j=1}^{\min(m,N_k)}\podSingularValue{j}^2}\ge\rho_{\mathrm{POD}},
\]
where \(\podSingularValue{j}\) are the singular values of \(\mathbf Y_k\).
The reduced training responses are
\[
\mathbf C_k=\mathbf\Phi_k^\top\mathbf Y_k
\in\mathbb R^{\podRank{k}\times N_k}.
\]
The POD is formed from the centered, unscaled responses. Before GP
regression, we also normalize each retained coefficient using its
training range:
\[
\widetilde C_{k,jn}
=\frac{C_{k,jn}-c_{j,k}^{\min}}{\Delta c_{j,k}},
\qquad
c_{j,k}^{\min}:=\min_{1\le n\le N_k}C_{k,jn},
\qquad
\Delta c_{j,k}:=\max_{1\le n\le N_k}C_{k,jn}-c_{j,k}^{\min}.
\]
We replace a zero range by one. 
For each mode $j$, a scalar GP is
trained on the normalized parameter inputs and the $j$th row of
$\widetilde{\mathbf C}_k$.

Letting $m_{j,k}$ denote the predicted mean,
the dimensional coefficient prediction is
\[
\hat c_{j,k}(\params)
=c_{j,k}^{\min}+\Delta c_{j,k}\,
m_{j,k}\!\left(\widetilde{\params}_k(\params)\right).
\]
Collecting the predicted POD coefficient
means in \(\hat{\boldsymbol c}_k(\params)\in\mathbb R^{\podRank{k}}\), the
vector-valued ROM response is
\[
\rom_k(\params)
=\bar{\fomQoi}_k+\mathbf\Phi_k\hat{\boldsymbol c}_k(\params).
\]

Algorithm~\ref{alg:otf_mfvi} summarizes the method for the default configuration: delayed
ROM adaptation and LOO scalar gain estimation for the gradient, with same-sample scalar gains for the Hessian. Separate scalar
gains are fitted for the two gradient blocks in
Eq.~\eqref{eq:scalar_gradient_gain_blocks} and the two Hessian blocks in
Eq.~\eqref{eq:scalar_hessian_gain_blocks}. Training sets retain the latest
$n_{\mathrm{hist}}$ HF batches.

\begin{algorithm}
\caption{Adaptive multifidelity VI (default configuration): delayed ROM updates, LOO scalar gradient gains, and same-sample scalar Hessian gains}
\label{alg:otf_mfvi}
\small
\begin{algorithmic}[1]
\Require $\searchDistributionParameters_0$, sample sizes, training and optimization controls
\State Draw $\sampleSetHighFidelity_0\subseteq\sampleSetLowFidelity_0$ from $\searchDistribution(\cdot;\searchDistributionParameters_0)$; evaluate HF on $\sampleSetHighFidelity_0$
\State Form $\mathcal T_0$ from these HF pairs and train $\rom_0=\rom(\cdot;\romTrainingMap(\mathcal T_0))$; set $\rom_1\gets\rom_0$
\State Evaluate LF on $\sampleSetLowFidelity_0$; set all initial blockwise scalar gains to $-1$
\State Compute $\approximateControlVariateGradient_0$, $\approximateControlVariateHessian_0$ by Eqs.~\eqref{eq:loo_mfmc_gradient}--\eqref{eq:loo_mfmc_hessian}
\State $\laggedHessian_0\gets\approximateControlVariateHessian_0$; $k\gets0$
\While{budget and stopping criteria permit}
  \State Form $\regularizedHessian_k$ from $\laggedHessian_k$ using Eq.~\eqref{eq:absolute_hessian}
  \State $\searchDistributionParameters_{k+1}\gets\searchDistributionParameters_k-\epsilon_k\regularizedHessian_k^{-1}\approximateControlVariateGradient_k$
  \State $\laggedHessian_{k+1}\gets\beta\laggedHessian_k+(1-\beta)\approximateControlVariateHessian_k$
  \State Draw $\sampleSetHighFidelity_{k+1}\subseteq\sampleSetLowFidelity_{k+1}$ from $\searchDistribution(\cdot;\searchDistributionParameters_{k+1})$; evaluate HF on $\sampleSetHighFidelity_{k+1}$
  \State Form $\mathcal T_{k+1}$ from the latest $n_{\mathrm{hist}}$ HF batches
  \State Evaluate the active ROM $\rom_{k+1}$ on $\sampleSetLowFidelity_{k+1}$ and form baseline-centered score contributions
  \State For each $i$, fit the two gradient scalar gains on $\sampleSetHighFidelity_{k+1,-i}$ (LOO); fit the two Hessian scalar gains on $\sampleSetHighFidelity_{k+1}$
  \State Compute $\approximateControlVariateGradient_{k+1}$ by Eq.~\eqref{eq:gain_loo_mfmc_gradient} and $\approximateControlVariateHessian_{k+1}$ by Eq.~\eqref{eq:loo_mfmc_hessian}
  \State Compute $\romError{k+1}$ using $\rom_{k+1}$ and the current HF responses
  \If{$\romError{k+1}>\tau_{\mathrm{ROM}}$}
    \State $\rom_{k+2}\gets\rom(\cdot;\romTrainingMap(\mathcal T_{k+1}))$ \Comment{used from the next iteration}
  \Else
    \State $\rom_{k+2}\gets\rom_{k+1}$
  \EndIf
  \State $k\gets k+1$
\EndWhile
\State \Return $\searchDistribution(\cdot;\searchDistributionParameters_k)$
\end{algorithmic}
\end{algorithm}

\subsection{Software implementation}
The MF-BBVI algorithms discussed above are implemented in the \texttt{rom-tools-and-workflows} Python package through a small set of model and builder abstractions. The package comprises Python drivers that manage the outer VI algorithms.  The high-fidelity forward model is supplied through the \texttt{QoiModel} interface, which encapsulates run-directory setup, model execution, and QoI extraction. In practice, the model execution layer involves an external call to a PDE solver. 

A simplified pseudocode view of the model API is
\begin{center}
\fcolorbox{black!35}{black!6}{%
\parbox{0.94\linewidth}{%
\ttfamily
class QoiModel:\par
\hspace*{1.5em}populate\_run\_directory(run\_directory, parameter\_sample)\par
\hspace*{1.5em}run\_model(run\_directory, parameter\_sample) -> int\par
\hspace*{1.5em}compute\_qoi(run\_directory, parameter\_sample) -> ndarray\par
}}
\end{center}
A user can gain access to the algorithms by creating a concrete instance of the QoiModel for their application.  

The low-fidelity models are additionally managed through a \texttt{QoiModel} API, and the on-the-fly construction of these low-fidelity models is handled through the \texttt{QoiModelBuilderWithTrainingData} abstraction. This builder accepts the current training directories, sampled parameters, and QoIs, and returns a new low-fidelity \texttt{QoiModel}. 
A corresponding pseudocode view of the builder API is
\begin{center}
\fcolorbox{black!35}{black!6}{%
\parbox{0.94\linewidth}{%
\ttfamily
class QoiModelBuilderWithTrainingData:\par
\hspace*{1.5em}build\_from\_training\_dirs(\par
\hspace*{3em}offline\_data\_dir,\par
\hspace*{3em}training\_data\_dirs,\par
\hspace*{3em}training\_parameters,\par
\hspace*{3em}training\_qois\par
\hspace*{1.5em}) -> QoiModel%
}}
\end{center}
We implement a concrete instance of the QoiModelBuilderWithTrainingData for the POD--GP ROM discussed above, but any arbitrary model can be substituted. 

Finally, and of practical importance, the implementation supports concurrent evaluation of both the HF and low-fidelity model ensembles. Since the evaluation of the HF and LF batches at each iteration can be performed in parallel, concurrent evaluation enables substantial gains in terms of wall time.

\section{Numerical experiments}
\label{sec:results}
We now assess the proposed multifidelity BBVI approaches on four inverse problems: a simple sine-series benchmark problem with an analytic solution, a steady convection--diffusion--reaction problem, a premixed H$_2$-air flame example, and a shallow-water bathymetry problem. In each case, we compare single-fidelity, natural gradient BBVI with the defacto first order optimizer Adam, single-fidelity BBVI with the stochastic Newton optimizer, and their multifidelity counterparts. We refer to these methods as BBVI (Adam), BBVI (Newton), MF-BBVI (Adam), and MF-BBVI (Newton), respectively. Table~\ref{tab:experiment_hyperparameters} summarizes the optimizer hyperparameters used across all examples.

\begin{table}[!htbp]
\centering
\small
\renewcommand{\arraystretch}{1.12}
\setlength{\tabcolsep}{8pt}
\caption{Shared optimizer hyperparameters for the sine-series, CDR, and H$_2$-air method comparisons. Settings apply to both single-fidelity and multifidelity variants.}
\label{tab:experiment_hyperparameters}
\begin{tabular}{@{}lr@{}}
\hline
\textbf{Parameter} & \textbf{Value} \\
\hline
\multicolumn{2}{@{}l}{\textit{Newton optimizer}} \\[2pt]
Diagonal threshold, $\rho$ & $1.00\times10^{-5}$ \\
Additive regularization, $\additiveRegularization_0$ & $0.1$ \\
Additive cap & $20$ \\
Additive growth factor & $1.02$ \\
Fallback learning rate, $\fallbackLearningRate$ & $0.05$ \\
Mean update cap & $0.5$ \\
\hline
\multicolumn{2}{@{}l}{\textit{Adam optimizer}} \\[2pt]
Learning-rate scale & $0.2$ \\
\hline
\end{tabular}
\end{table}

Unless otherwise noted, all MF methods use delayed ROM updates and a scalar gain computed using the LOO approach; this configuration results in an unbiased estimator. The Newton methods use lagged curvature with an averaging factor of $\beta=0.1$. Following Ref.~\cite{rei2025blackboxvariationalinference}, Adam uses natural gradients obtained by preconditioning the ELBO gradient with $(\fisherMatrix+\delta_{\mathrm{F},k}\mathbf{I})^{-1}$ before the Adam update. The positive Fisher regularization coefficient $\delta_{\mathrm{F},k}$ limits gradient amplification when Fisher eigenvalues are small; it starts at $10^{-2}$, decays exponentially after iteration 50, and has a floor of $10^{-6}$. The effective learning rate equals the listed scale divided by the parameter dimension.

The multifidelity methods use the Gaussian-process surrogate described in Section~\ref{sec:otf_rom_gp}, constructed from the four most recent HF batches. Parameters and responses are normalized as described in Section~\ref{sec:otf_rom_gp}. We use a ROM-update tolerance of $\tau_{\mathrm{ROM}}=1\times10^{-3}$ and a retained POD energy fraction of $0.999999$. The GP employs a squared-exponential covariance kernel
\[
\mathcal K\!\left(\widetilde{\params}_k(\params),\widetilde{\params}_k(\params')\right)
=\sigma_f^2\exp\!\left(-\frac{\|\widetilde{\params}_k(\params)-\widetilde{\params}_k(\params')\|_2^2}{2\ell^2}\right),
\]
where $\widetilde{\params}_k(\params)$ and $\widetilde{\params}_k(\params')$ are the normalized parameter inputs defined in Section~\ref{sec:otf_rom_gp},
$\ell>0$ is an isotropic length scale, and $\sigma_f^2>0$ is the signal variance.
We select these hyperparameters independently for each retained POD mode
by maximizing the GP log marginal likelihood over a tensor-product grid.
We use nine logarithmically spaced length scales in $[10^{-2},10^2]$ and
nine signal variances in $[10^{-6},10^2]$, giving 81 candidate pairs per mode.

In practice, evaluating a GP is trivial relative to a HF evaluation, and thus we measure computational work by the cumulative number of HF evaluations. We choose the initial HF batch size using the number of parameters as a guide. Our default is to round the number of parameters to the nearest power of two on a logarithmic scale, giving a batch size of $2^{\operatorname{round}(\log_2(\nParams))}$.
To measure the accuracy of the various approaches, we examine both the KL divergence and variational parameters against a high-fidelity reference solution. 
For the sine-series problem, this reference solution is an exact Gaussian posterior. For the remaining problems, we employ a separately computed, high-fidelity numerical VI reference in the same variational family. Descriptions of these reference solutions are given in the experiment descriptions. We emphasize that these latter values measure agreement with the reference distribution rather than the exact Bayesian posterior; the bias incurred by the variational approximation is outside the scope of the present contribution.
In all cases, the KL divergence is evaluated analytically from the variational parameters. Writing the Gaussian distributions in the coordinates used for this calculation as $\mathcal N(\variationalMean_k,\boldsymbol\Sigma_k)$ and $\mathcal N(\variationalMean_{\mathrm{ref}},\boldsymbol\Sigma_{\mathrm{ref}})$, we use
\begin{equation}\label{eq:numerical_gaussian_kl}
D_{\mathrm{KL}}(q_k\|q_{\mathrm{ref}})
=\frac12\Big[\operatorname{tr}(\boldsymbol\Sigma_{\mathrm{ref}}^{-1}\boldsymbol\Sigma_k)
+(\variationalMean_k-\variationalMean_{\mathrm{ref}})^\top
\boldsymbol\Sigma_{\mathrm{ref}}^{-1}
(\variationalMean_k-\variationalMean_{\mathrm{ref}})
-\nParams+\log\frac{\det\boldsymbol\Sigma_{\mathrm{ref}}}{\det\boldsymbol\Sigma_k}\Big].
\end{equation}
For the CDR, air-flame, and shallow-water bathymetry problems, the Gaussian parameters in the above are the latent-coordinate moments $\boldsymbol\mu^*$ and $\boldsymbol\Sigma^*$ introduced in Section~\ref{sec:bounds}. Physical means and variances are evaluated by one-dimensional quadrature.

\subsection{Gaussian sine-series inverse problem}
\label{sec:analytic_sine_benchmark}
\begin{figure}
\centering
\includegraphics[width=\linewidth]{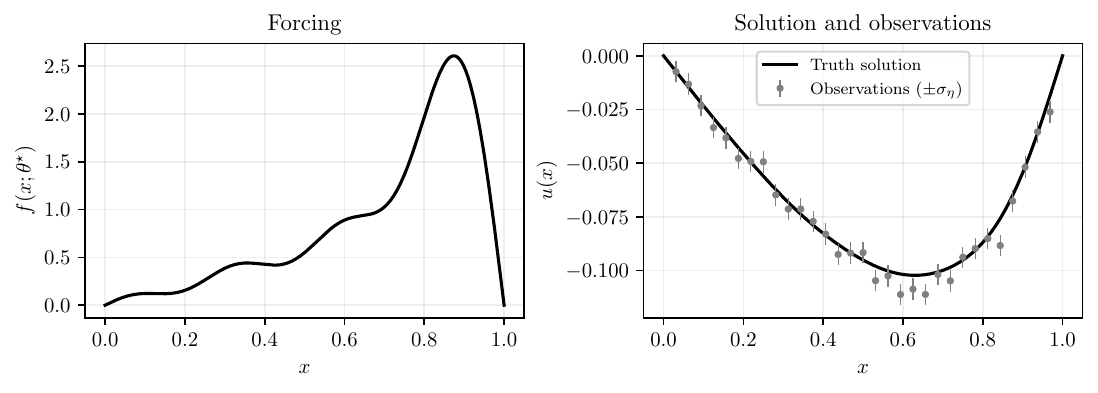}
\caption{Sine-series problem: forcing and solution at the generating parameters. Observation error bars indicate one noise standard deviation.}
\label{fig:sine_solution}
\end{figure}

We first consider a linear-Gaussian inverse problem for which the exact posterior is available. The forward problem on $x \in (0,1)$ is
\begin{equation}
 u''(x)=f(x;\params),\qquad u(0)=u(1)=0,
 \qquad f(x;\params)=\sum_{j=1}^{7}\theta_j\sin(j\pi x).
\end{equation}
The solution is
\begin{equation}
 u(x;\params)=-\sum_{j=1}^{7}\frac{\theta_j}{(j\pi)^2}\sin(j\pi x).
\end{equation}
We generate observations at 31 equispaced interior points, $x_i=i/32$, using
\begin{equation}
    y_i=u(x_i;\params^{\mathrm{truth}})+\eta_i,\qquad
 \eta_i\overset{\mathrm{iid}}{\sim}\mathcal N(0,(5\times10^{-3})^2),\qquad
    \params^{\mathrm{truth}}=(1,-0.8,0.65,-0.5,0.35,-0.25,0.15)^\top.
\end{equation}
The prior and initial variational distribution are both $\mathcal N(0,0.75^2\mathbf I)$. Figure~\ref{fig:sine_solution} shows the forcing, solution, and noisy observations. 

Writing the observation model as $\obs=\mathbf G\params+\eta$, the observation matrix has entries $G_{ij}=-(j\pi)^{-2}\sin(ij\pi/32)$. The exact posterior mean and covariance satisfy
\begin{equation}
 \boldsymbol\Sigma_{\mathrm{post}}^{-1}=\mathbf C_0^{-1}+\mathbf G^\top\boldsymbol\Gamma^{-1}\mathbf G,
 \qquad \boldsymbol\mu_{\mathrm{post}}=\boldsymbol\Sigma_{\mathrm{post}}\mathbf G^\top\boldsymbol\Gamma^{-1}\obs,
\end{equation}
where $\mathbf C_0=0.75^2\mathbf I$ and $\boldsymbol\Gamma=(5\times10^{-3})^2\mathbf I$. Noting that $\mathbf G^\top\boldsymbol\Gamma^{-1}\mathbf G$ is diagonal, the posterior belongs to the mean-field Gaussian variational family.

\subsubsection{Method comparison}
\begin{figure}
\centering
\includegraphics[width=0.92\linewidth]{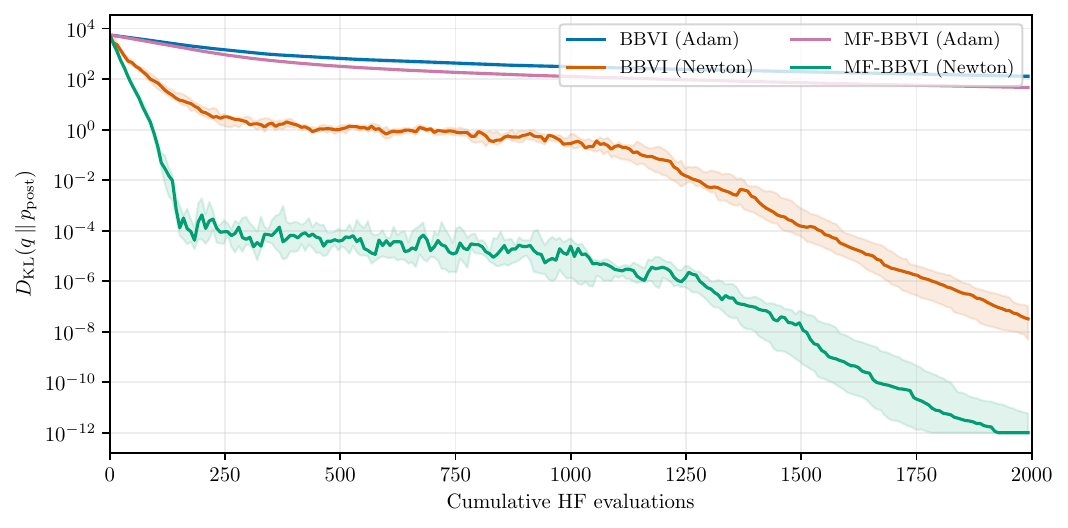}
\caption{Sine-series method comparison: KL versus cumulative HF evaluations. Curves show medians over ten realizations and bands show interquartile ranges. Values below $10^{-12}$ are clipped for plotting.}
\label{fig:analytic_sine_benchmark}
\end{figure}

\begin{figure}
\centering
\includegraphics[width=0.95\linewidth]{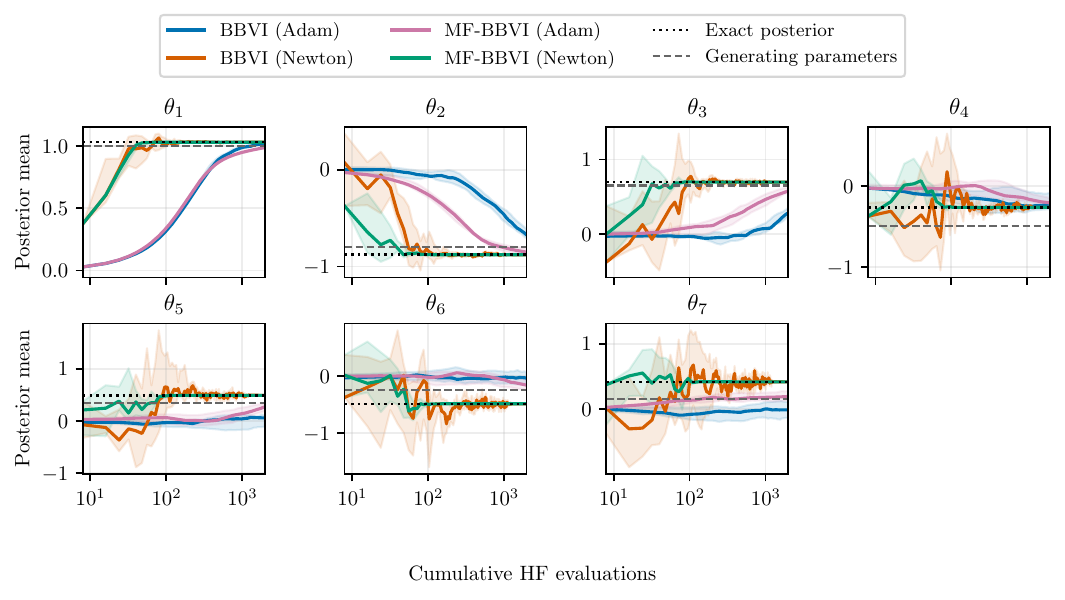}
\caption{Sine-series posterior means. Black dotted lines mark the exact posterior means; gray dashed lines mark the generating parameters.}
\label{fig:sine_moments}
\end{figure}

\begin{figure}
\centering
\includegraphics[width=0.95\linewidth]{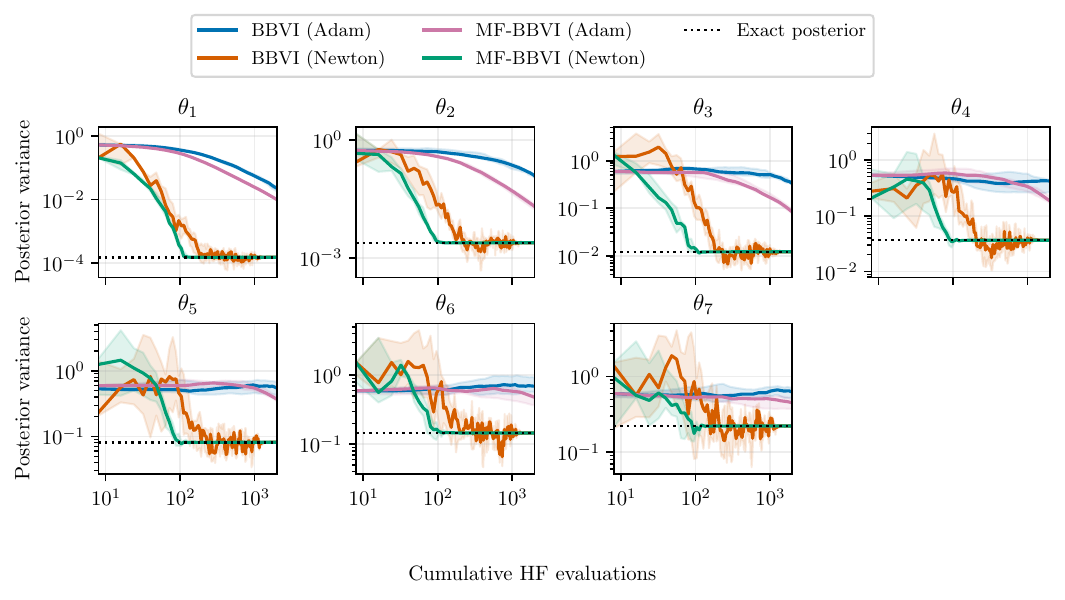}
\caption{Sine-series marginal variances. Black dotted lines mark the exact posterior variances.}
\label{fig:sine_moments_variance}
\end{figure}

\begin{figure}
\centering
\includegraphics[width=\linewidth]{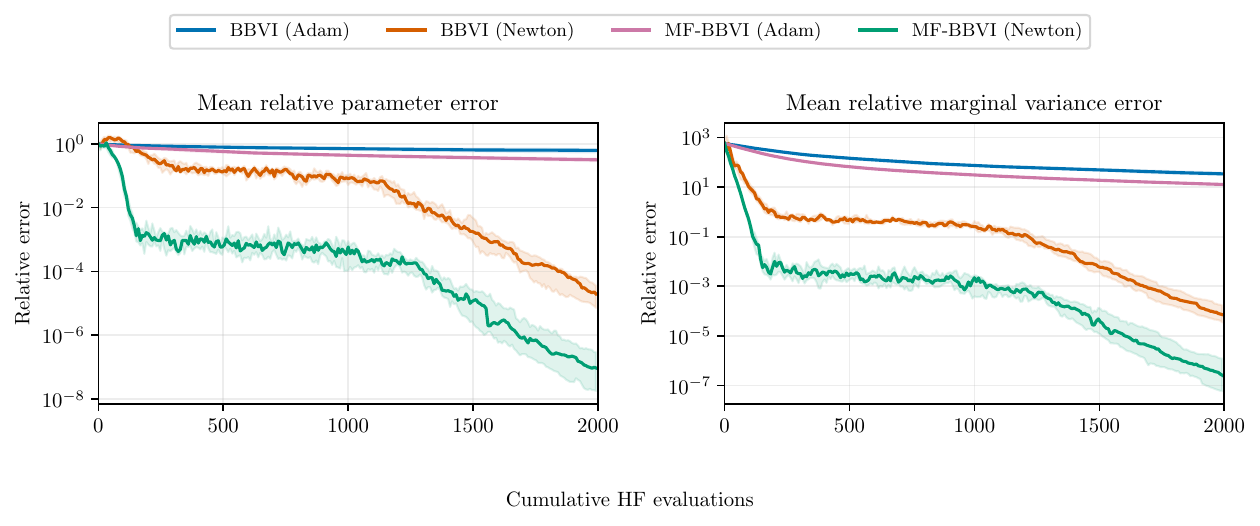}
\caption{Sine-series mean relative parameter error (left) and marginal variance error (right) with respect to the exact posterior.}
\label{fig:sine_relative_moment_errors}
\end{figure}

We first assess each method given a fixed computational budget of 2000 HF evaluations. We employ 8 HF samples per batch, and 512 additional LF samples for the MF methods. As optimization is stochastic, each method is repeated 10 times. Both BBVI (Newton) and MF-BBVI (Newton) approach the exact posterior (Figure~\ref{fig:analytic_sine_benchmark}, Table~\ref{tab:sine_method_ablation}). The terminal median KL value is below $10^{-12}$ for MF-BBVI (Newton) and below $10^{-7}$ for BBVI (Newton). BBVI (Adam) has the highest median KL of $130$, while MF-BBVI (Adam) substantially improves upon this with a terminal median KL of $46.9$. Multifidelity sampling improves both optimizers. 
Figures~\ref{fig:sine_moments}--\ref{fig:sine_relative_moment_errors}
show convergence of the mean and variances. We observe that MF-BBVI (Newton) closely approximates the correct posterior within a hundred iterations for all quantities. BBVI (Newton) is the next best performing method. Both Adam-based methods never reach the true posterior for the considered budget.

\begin{table}[!htbp]
\centering
\caption{Sine method comparison: median KL at the last recorded iteration across 10 seeds. Green and red mark the smallest and largest median KL, respectively.}
\label{tab:sine_method_ablation}
\begin{tabular}{lr}
\hline
Method & Median KL \\
\hline
BBVI (Adam) & \textcolor{red}{$130$} \\
BBVI (Newton) & $3.23\times10^{-8}$ \\
MF-BBVI (Adam) & $46.9$ \\
MF-BBVI (Newton) & \textcolor{bestkl}{$<10^{-12}$} \\
\hline
\end{tabular}
\end{table}


\subsubsection{Sample-size study}
\begin{figure}
\centering
\includegraphics[width=\linewidth]{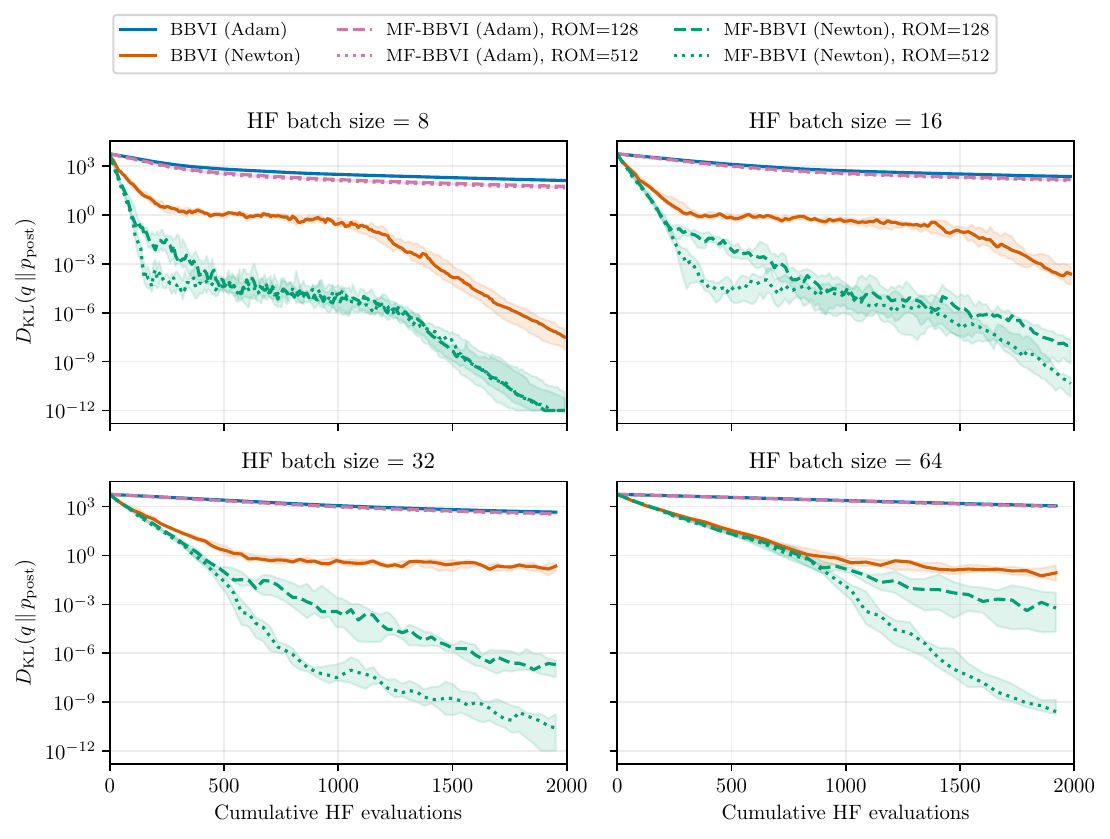}
\caption{Sine-series sample-size study: exact-posterior KL for each HF batch size. Solid and dashed multifidelity curves use 128 and 512 additional LF samples, respectively. Curves and bands show medians and interquartile ranges. Values below $10^{-12}$ are clipped for plotting.}
\label{fig:analytic_sine_sample_sweep}
\end{figure}

To assess the effect of sample allocation, we vary $\nSamplesHF\in\{8,16,32,64\}$ and use either 128 or 512 additional LF samples, retaining the 2000-HF budget. Figure~\ref{fig:analytic_sine_sample_sweep} and Table~\ref{tab:sine_sample_size} show the resulting KL histories and terminal medians. MF-BBVI (Newton) has the smallest terminal median for all configurations. In all cases, increasing the additional number of LF samples from 128 to 512 results in improved performance. Increasing the number of HF samples results in more stable optimization, but comes at the cost of more expensive iterations. As a result, for a fixed budget, increasing the HF sample size at each iteration does not result in uniform improvement. Across all configurations, we observe the same hierarchy in terms of model performance: MF-BBVI (Newton) is the best performing method, followed by BBVI (Newton). MF-BBVI (Adam) outperforms BBVI (Adam) in every case, but yields substantially worse performance than the Newton methods. 
\begin{table}[!htbp]
\centering
\small
\caption{Sine sample-size study: median terminal KL across 10 seeds. ROM denotes additional low-fidelity samples for the multifidelity methods. Single-fidelity values are repeated for comparison. Green and red mark the smallest and largest KL in each row.}
\label{tab:sine_sample_size}
\begin{tabular}{rr|rrrr}
\hline
 & & BBVI & BBVI & MF-BBVI & MF-BBVI \\
HF batch & ROM & (Adam) & (Newton) & (Adam) & (Newton) \\
\hline
8 & 128 & \textcolor{red}{$130$} & $3.23\times10^{-8}$ & $57.7$ & \textcolor{bestkl}{$1.78\times10^{-13}$} \\
8 & 512 & \textcolor{red}{$130$} & $3.23\times10^{-8}$ & $46.9$ & \textcolor{bestkl}{$4.04\times10^{-13}$} \\
16 & 128 & \textcolor{red}{$229$} & $2.43\times10^{-4}$ & $150$ & \textcolor{bestkl}{$1.08\times10^{-8}$} \\
16 & 512 & \textcolor{red}{$229$} & $2.43\times10^{-4}$ & $132$ & \textcolor{bestkl}{$4.55\times10^{-11}$} \\
32 & 128 & \textcolor{red}{$456$} & $0.22$ & $373$ & \textcolor{bestkl}{$2.01\times10^{-7}$} \\
32 & 512 & \textcolor{red}{$456$} & $0.22$ & $344$ & \textcolor{bestkl}{$2.54\times10^{-11}$} \\
64 & 128 & \textcolor{red}{$1.11\times10^{3}$} & $0.0855$ & $1.05\times10^{3}$ & \textcolor{bestkl}{$5.91\times10^{-4}$} \\
64 & 512 & \textcolor{red}{$1.11\times10^{3}$} & $0.0855$ & $997$ & \textcolor{bestkl}{$2.60\times10^{-10}$} \\
\hline
\end{tabular}
\end{table}

\subsubsection{Multifidelity studies}
\begin{figure}
\centering
\includegraphics[width=\linewidth]{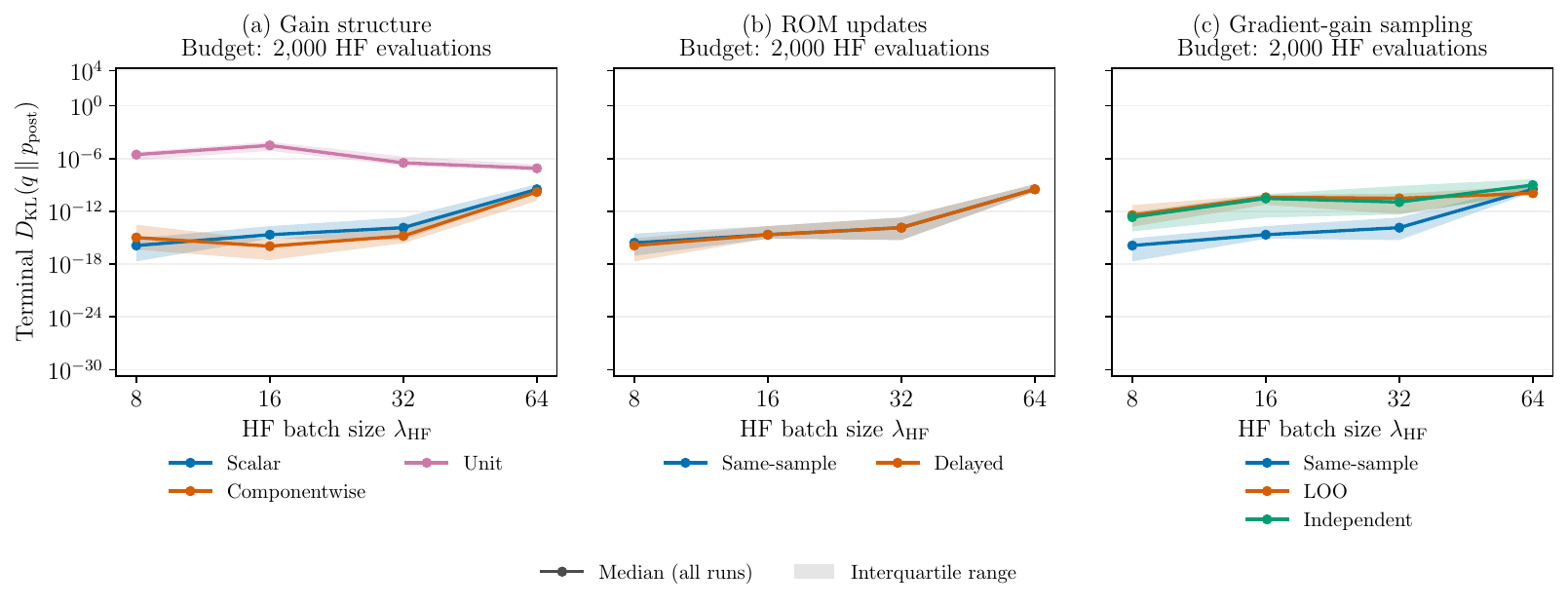}
\caption{Sine-series multifidelity studies: terminal exact-posterior KL versus HF batch size for (a) gain structure, (b) LF-update timing, and (c) gradient-gain sampling. Lines and markers show medians across ten runs, and shaded bands span the 25th--75th percentiles.}
\label{fig:sine_mf_ablation_summary}
\end{figure}

We now examine the impact of the gain structure, the manner in which the LF model is updated, and the manner in which the gain is estimated. We restrict our study to MF-BBVI (Newton), using ten realizations per configuration and 512 additional LF samples. Figure~\ref{fig:sine_mf_ablation_summary}(a--c) summarizes terminal KL versus HF batch size.%

For gain structure, we compare fitted scalar and componentwise gains with fixed unit-magnitude gains in both gradient and Hessian estimators. Figure~\ref{fig:sine_mf_ablation_summary}(a) shows that the componentwise and scalar gains clearly outperform the unit gains for all HF batch sizes. The scalar and componentwise gains perform similarly. These results suggest that increasing the complexity of the gain structure does not necessarily lead to improved estimators due to the finite sampling error associated with computing the gain.

For LF model adaptation, we compare same-sample and delayed updates with same-sample scalar gradient and Hessian gains, using $\nSamplesHF\in\{8,16,32,64\}$ and a 2000-HF budget (Figure~\ref{fig:sine_mf_ablation_summary}(b)). The two methods perform similarly, suggesting that the bias induced from immediately using a retrained LF model has a minimal impact on the optimization. 

Lastly, for gradient-gain sampling, we compare same-sample, LOO, and independent scalar gains with delayed LF updates and same-sample Hessian gains (Figure~\ref{fig:sine_mf_ablation_summary}(c)). Despite being biased, the same-sample approach slightly outperforms LOO and independent sampling. We again emphasize that independent gain-fitting evaluations are not included in the budget; in practice, these runs use approximately twice the total HF work. Thus, there is not a practical advantage from separating gain-fitting samples in this example. 

\subsubsection{High-dimensional configuration}
\label{sec:sine_highdim}
To examine performance in a higher-dimensional configuration, we briefly consider a configuration where we infer a 30-term sine-series expansion with generating coefficients
\[
\theta_j^{\mathrm{truth}}=0.6\sin(0.73j)+0.3\cos(1.37j),\qquad j=1,\ldots,30.
\]
We increase the number of observations to 127, located at $x_i=i/128$, while
retaining the prior, initialization, and observation-noise level of the previous configuration. The exact posterior remains a diagonal Gaussian. 
We compare the four methods using 32 HF samples
per batch, 512 additional LF samples for the multifidelity methods, and a nominal
budget of 2000 HF evaluations, with five realizations per method. 

\begin{figure}
\centering
\includegraphics[width=0.92\linewidth]{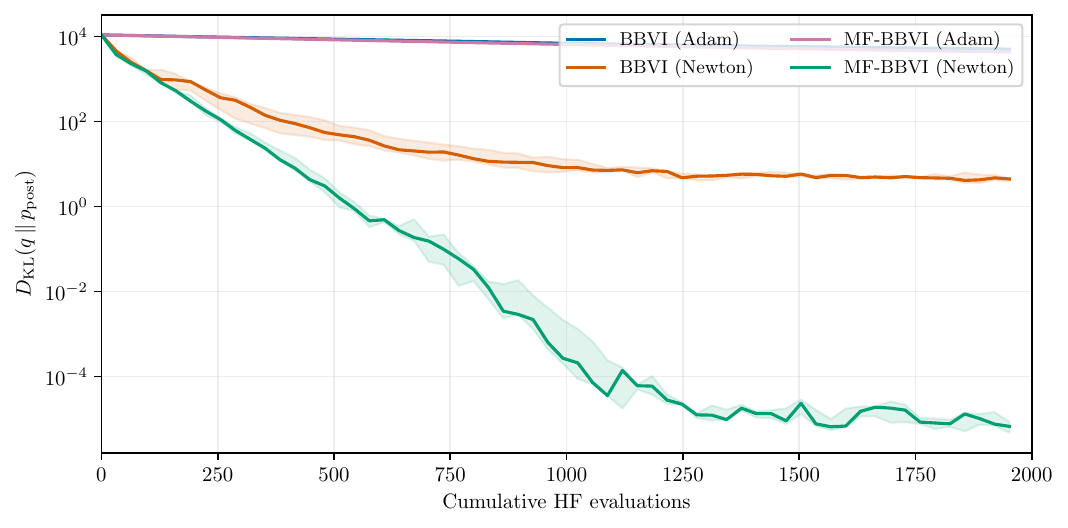}
\caption{Thirty-term sine-series method comparison: KL to the exact posterior versus cumulative HF evaluations. Curves and bands show medians and interquartile ranges across five realizations. The final recorded iterates use 1952 HF evaluations within the nominal 2000-evaluation budget.}
\label{fig:sine_highdim_method_ablation}
\end{figure}

\begin{figure}
\centering
\includegraphics[width=\linewidth]{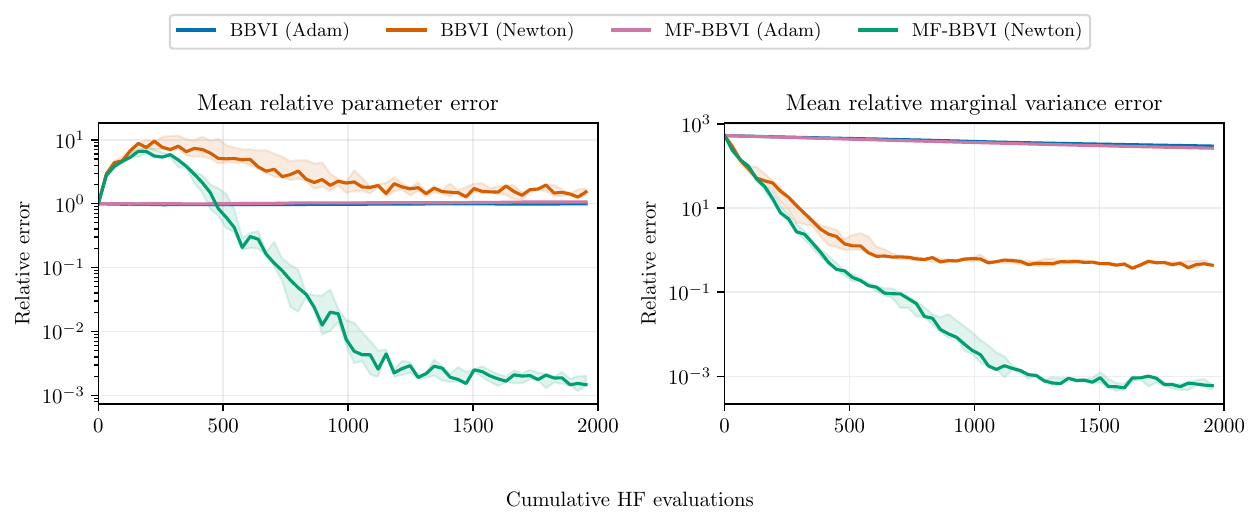}
\caption{Thirty-term sine-series method comparison: mean componentwise relative parameter error (left) and mean componentwise relative marginal variance error (right) with respect to the exact posterior. Curves and bands show medians and interquartile ranges across five realizations.}
\label{fig:sine_highdim_relative_moment_errors}
\end{figure}

Figure~\ref{fig:sine_highdim_method_ablation} shows KL convergence, and
Table~\ref{tab:sine_30_method_ablation} reports the terminal medians.
In this higher-dimensional configuration, MF-BBVI (Newton) still closely approximates the exact posterior, yielding a median
terminal KL of $6.73\times10^{-6}$, compared with $4.48$ for BBVI (Newton).
The corresponding medians for MF-BBVI (Adam) and BBVI (Adam) are
$4.38\times10^3$ and $5.08\times10^3$, respectively. Multifidelity sampling
improves both optimizers, and Newton again outperforms Adam. 
Figure~\ref{fig:sine_highdim_relative_moment_errors} shows the average componentwise
relative mean and marginal variance errors. MF-BBVI (Newton) attains
terminal median mean and variance errors of $0.147\%$ and $0.0601\%$, respectively,
compared with $153\%$ and $44.0\%$ for BBVI (Newton). The corresponding values
are $106\%$ and $2.67\times10^4\%$ for MF-BBVI (Adam), and $99.3\%$ and
$3.00\times10^4\%$ for BBVI (Adam). 

\begin{table}[!htbp]
\centering
\caption{30-term sine-series method comparison: median KL at the last recorded iteration across 5 seeds. Green and red mark the smallest and largest median KL, respectively.}
\label{tab:sine_30_method_ablation}
\begin{tabular}{lr}
\hline
Method & Median KL \\
\hline
BBVI (Adam) & \textcolor{red}{$5.08\times10^{3}$} \\
BBVI (Newton) & $4.48$ \\
MF-BBVI (Adam) & $4.38\times10^{3}$ \\
MF-BBVI (Newton) & \textcolor{bestkl}{$6.73\times10^{-6}$} \\
\hline
\end{tabular}
\end{table}

\FloatBarrier
\subsection{Steady convection--diffusion--reaction inverse problem}
\label{sec:cdr_vi_benchmark}
\begin{figure}
\centering
\includegraphics[width=0.66\linewidth]{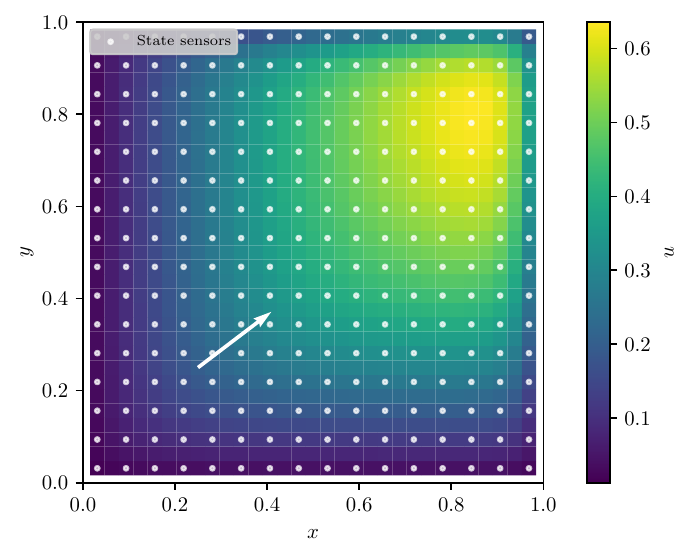}
\caption{CDR truth solution. White points mark the 256 state sensors. The white arrow indicates the convection direction.}
\label{fig:cdr_solution}
\end{figure}

We next consider a steady convection--diffusion--reaction (CDR) problem on $\Omega=[0,1]^2$,
\begin{equation}
 \nu\nabla^2u-\mathbf b\cdot\nabla u-\reactionCoeff u+1=0,
 \qquad \mathbf b=b_{\mathrm{mag}}[\cos\convectionAngle,\sin\convectionAngle]^\top,
\end{equation}
with homogeneous Dirichlet boundary conditions. Here $b_{\mathrm{mag}}$ is the convection speed, $\convectionAngle$ the convection angle (in radians), $\nu$ the diffusivity, and $\reactionCoeff$ the reaction coefficient. We infer $\params=(b_{\mathrm{mag}},\convectionAngle,\nu,\reactionCoeff)$ using synthetic observations generated at $\params^\star=(1.2,0.65,0.035,0.45)$. The finite-difference discretization has $31\times31$ interior unknowns. We employ second-order central differences for the diffusion term and second-order upwind differences for the convection term. Figure~\ref{fig:cdr_solution} shows the truth solution and observation locations. 

We observe the state at every other interior point in each coordinate direction, giving 256 observations. Independent Gaussian noise is added with standard deviation
\begin{equation}
 \sigma_u=0.05\frac{\|\obs_u^\star\|_2}{\sqrt{256}},
 \qquad \Gamma=\sigma_u^2I_{256},
\end{equation}
where the starred vector denotes the noiseless observations.

We restrict the convection angle to $[0,2\pi]$ and the remaining parameters to $[0,10]$, and use a Gaussian prior,
\begin{equation*}
 p(\params)\propto
 \mathcal N(\params;\boldsymbol\mu_0,\boldsymbol\Sigma_0),
 \qquad
 \boldsymbol\mu_0=(0.9,0.8,0.05,0.6)^\top,
\end{equation*}
where $\boldsymbol\Sigma_0=\operatorname{diag}(0.45^2,0.40^2,0.025^2,0.30^2)$.
We use a diagonal-Gaussian variational distribution in the transformed coordinates $\params^*=\boldsymbol g(\params)$ defined in~\eqref{eq:vi_sigmoid_bounds}. Its initial distribution is
\begin{equation*}
 \params^*\sim q_0=\mathcal N(\boldsymbol\mu_0^*,\boldsymbol\Sigma_0^*),
 \qquad
 \boldsymbol\mu_0^*=\boldsymbol g(\boldsymbol\mu_0),
 \qquad
 \boldsymbol\Sigma_0^*=
 \boldsymbol J_{\boldsymbol g}(\boldsymbol\mu_0)
 \boldsymbol\Sigma_0
 \boldsymbol J_{\boldsymbol g}(\boldsymbol\mu_0)^\top,
\end{equation*}
The prior and initialization use the same Gaussian mean and covariance.

This benchmark does not admit an analytic posterior. As a result, we compare our results to a numerically obtained ``truth solution'' computed using 2000 VI (Newton) updates with 1024 HF samples per update. The reference is initialized at the generating parameters, with standard deviations equal to 5\% of the corresponding parameter values. We use this saved Newton reference for all CDR studies. Its physical mean is approximately $(1.2039,0.6563,0.03518,0.4441)$.

\subsubsection{Method comparison}
\begin{figure}
\centering
\includegraphics[width=0.92\linewidth]{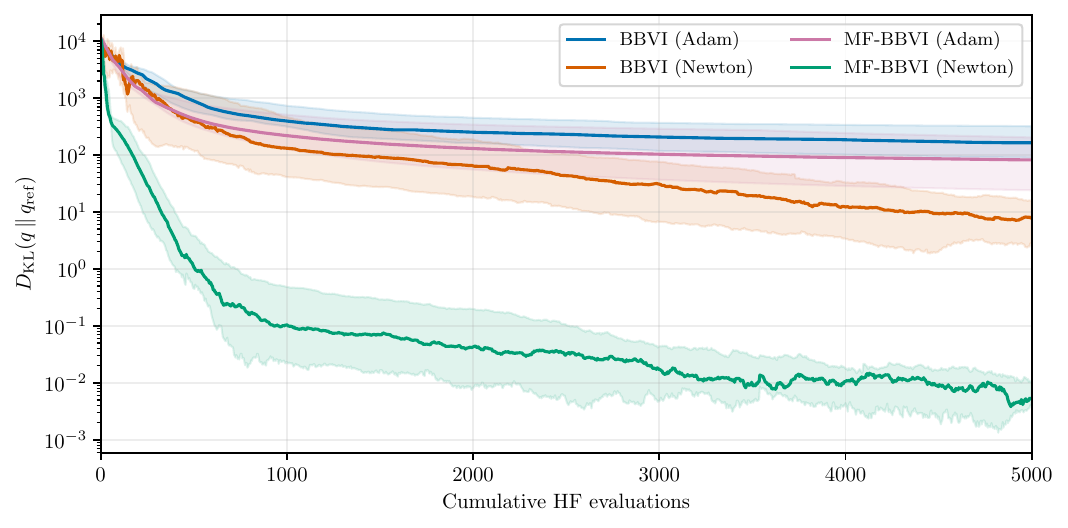}
\caption{CDR method comparison: KL to the numerical VI reference versus cumulative HF evaluations. Curves and bands show medians and interquartile ranges.}
\label{fig:cdr_method_ablation}
\end{figure}

We again compare the performance of the various methods given a fixed budget.  
All methods use four HF samples per batch and a nominal budget of 5000 HF evaluations. The multifidelity cases use 512 additional ROM samples. As optimization is stochastic, each method is again repeated 10 times. Figure~\ref{fig:cdr_method_ablation} shows the KL histories, and Table~\ref{tab:cdr_method_ablation} reports the terminal medians. MF-BBVI (Newton) yields the smallest median terminal KL of $5.36\times10^{-3}$. BBVI (Newton) is the next best performing method, with a median terminal KL of $7.79$. The corresponding values for MF-BBVI (Adam) and BBVI (Adam) are $82$ and $164$, respectively. Multifidelity sampling again improves both optimizers, and both Newton methods outperform their Adam counterparts.

\begin{table}[!htbp]
\centering
\caption{CDR method comparison: median KL at the last recorded iteration across 10 seeds. Green and red mark the smallest and largest median KL, respectively.}
\label{tab:cdr_method_ablation}
\begin{tabular}{lr}
\hline
Method & Median KL \\
\hline
BBVI (Adam) & \textcolor{red}{$164$} \\
BBVI (Newton) & $7.79$ \\
MF-BBVI (Adam) & $82$ \\
MF-BBVI (Newton) & \textcolor{bestkl}{$5.36\times10^{-3}$} \\
\hline
\end{tabular}
\end{table}

\begin{figure}
\centering
\includegraphics[width=0.92\linewidth]{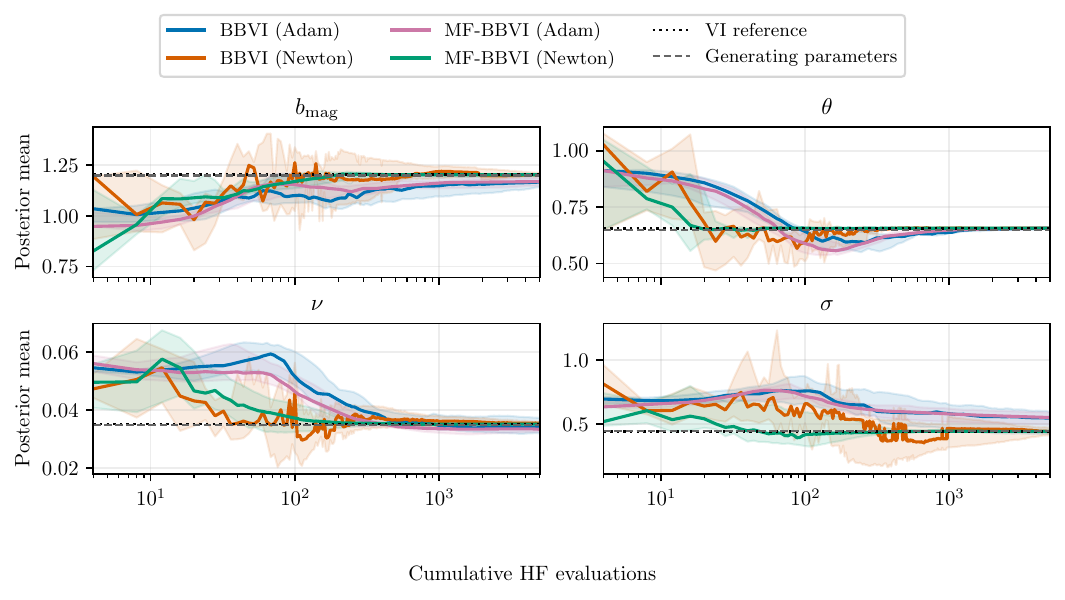}
\includegraphics[width=0.92\linewidth]{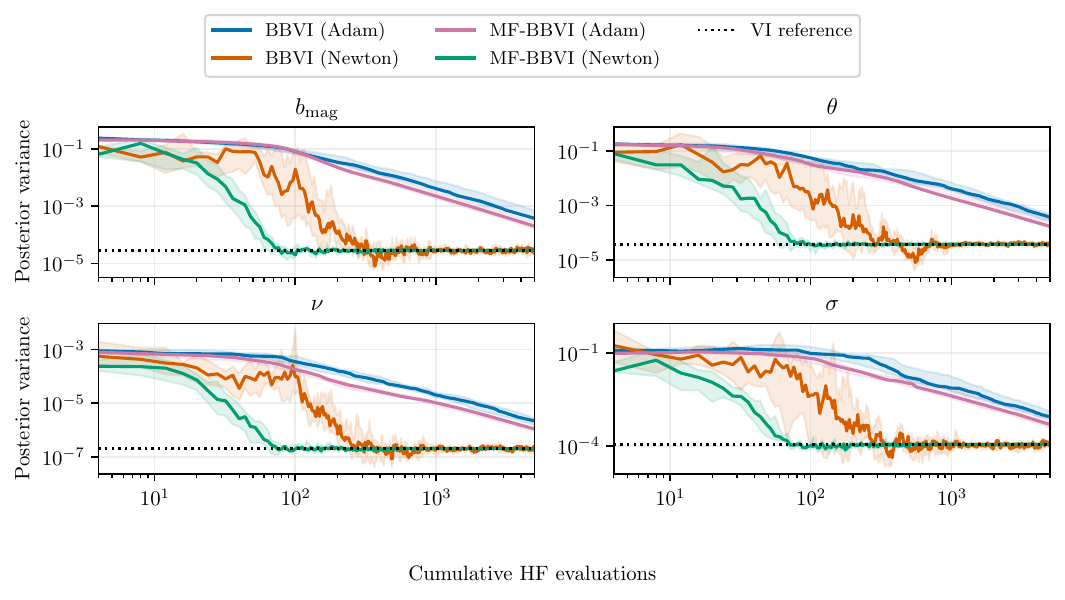}
\caption{CDR method comparison: physical posterior means (top) and variances (bottom). Curves show medians and interquartile ranges. Black dotted lines mark the VI reference moments; gray dashed lines mark the generating parameters.}%
\label{fig:cdr_moments}
\end{figure}

\begin{figure}
\centering
\includegraphics[width=\linewidth]{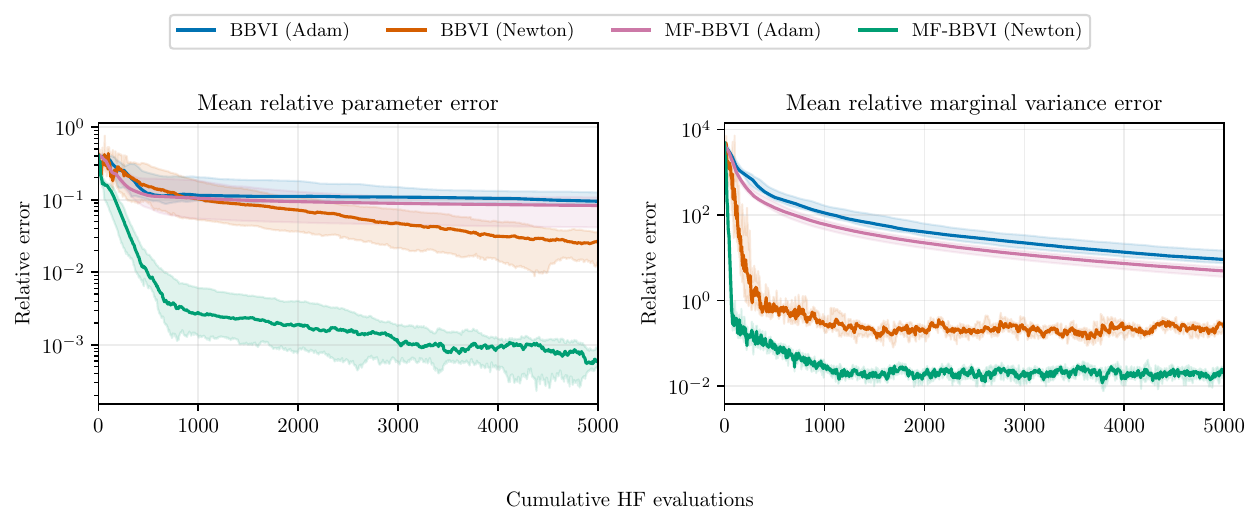}
\caption{CDR method comparison: mean relative parameter error (left) and mean relative marginal variance error (right) with respect to the numerical VI reference. Curves and bands show medians and interquartile ranges.}
\label{fig:cdr_relative_moment_errors}
\end{figure}

Figure~\ref{fig:cdr_moments} shows the physical mean and variance histories, while Figure~\ref{fig:cdr_relative_moment_errors} presents relative errors against the numerical VI reference solutions. The errors are averaged across parameters before summarizing the ten runs. We again observe that MF-BBVI (Newton) provides the best performance, and accurately approximates the posterior within 1000 FOM evaluations. BBVI (Newton) is the next best performing method, followed by MF-BBVI (Adam), and BBVI (Adam). MF-BBVI (Newton) also gives the smallest terminal median variance error, followed by BBVI (Newton), MF-BBVI (Adam), and BBVI (Adam). The MF estimators improve both error measures for both optimizers.

\subsubsection{Sample-size study}
\begin{figure}
\centering
\includegraphics[width=\linewidth]{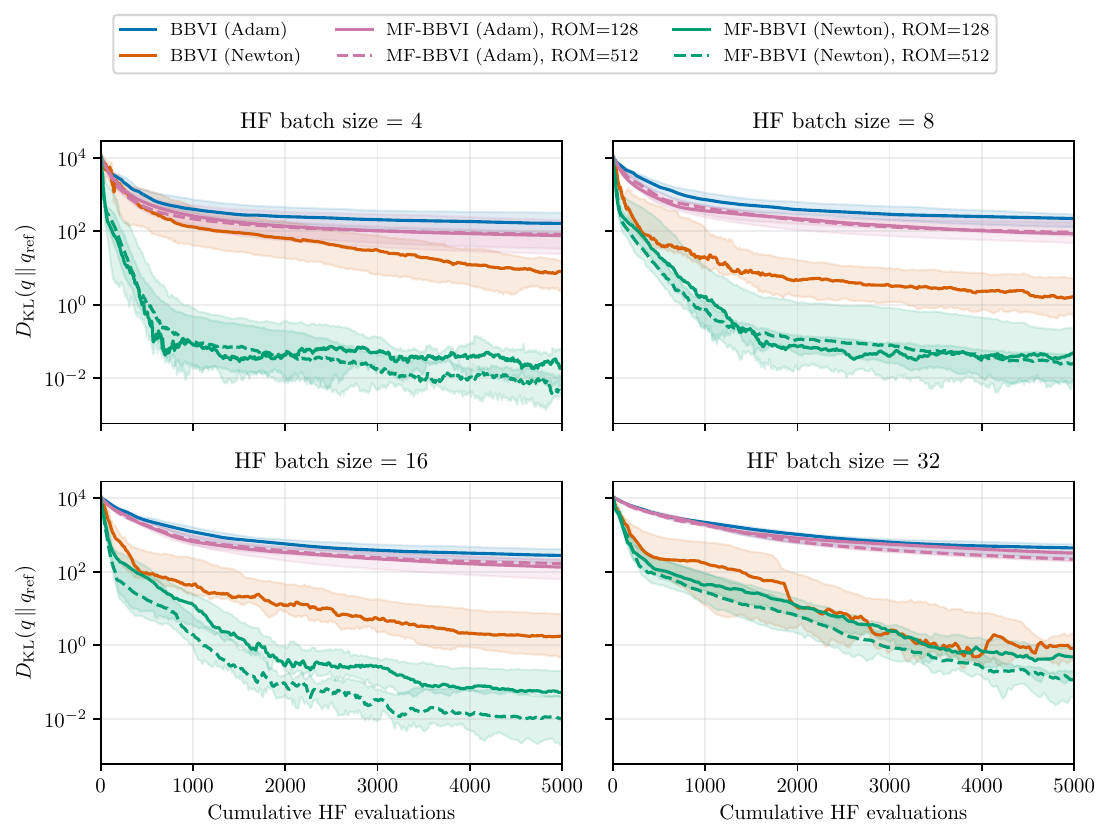}
\caption{CDR sample-size study: KL to the numerical VI reference solutions. Subfigures show a fixed HF batch size; solid and dashed multifidelity curves use 128 and 512 additional ROM samples. Curves and bands show medians and interquartile ranges.}
\label{fig:cdr_sample_size_ablation}
\end{figure}

Next, we vary the HF batch size over $\{4,8,16,32\}$ and, for the MF methods, the additional LF sample count over $\{128,512\}$. We retain the 5000-evaluation budget. Figure~\ref{fig:cdr_sample_size_ablation} shows KL convergence, and Table~\ref{tab:cdr_sample_size} reports the terminal medians. MF-BBVI (Newton) has the smallest terminal median for every tested configuration. Increasing the additional ROM sample count from 128 to 512 reduces its median in every configuration. The dependence on batch size varies across methods, since increasing the batch size results in more work per optimization step. The MF estimators improve the terminal median KL for both optimizers in every configuration considered. Both Newton methods outperform their Adam counterparts for every HF batch size.

\begin{table}[!htbp]
\centering
\small
\caption{CDR sample-size study: median terminal KL across 10 seeds. ROM denotes additional low-fidelity samples for the multifidelity methods. Single-fidelity values are repeated for comparison. Green and red mark the smallest and largest KL in each row.}
\label{tab:cdr_sample_size}
\begin{tabular}{rr|rrrr}
\hline
 & & BBVI & BBVI & MF-BBVI & MF-BBVI \\
HF batch & ROM & (Adam) & (Newton) & (Adam) & (Newton) \\
\hline
4 & 128 & \textcolor{red}{$164$} & $7.79$ & $76.4$ & \textcolor{bestkl}{$0.0221$} \\
4 & 512 & \textcolor{red}{$164$} & $7.79$ & $82$ & \textcolor{bestkl}{$5.36\times10^{-3}$} \\
8 & 128 & \textcolor{red}{$223$} & $1.66$ & $85.6$ & \textcolor{bestkl}{$0.0491$} \\
8 & 512 & \textcolor{red}{$223$} & $1.66$ & $93.2$ & \textcolor{bestkl}{$0.0252$} \\
16 & 128 & \textcolor{red}{$279$} & $1.79$ & $134$ & \textcolor{bestkl}{$0.0521$} \\
16 & 512 & \textcolor{red}{$279$} & $1.79$ & $168$ & \textcolor{bestkl}{$0.01$} \\
32 & 128 & \textcolor{red}{$447$} & $0.825$ & $323$ & \textcolor{bestkl}{$0.483$} \\
32 & 512 & \textcolor{red}{$447$} & $0.825$ & $220$ & \textcolor{bestkl}{$0.115$} \\
\hline
\end{tabular}
\end{table}

\subsubsection{Multifidelity studies}
\begin{figure}
\centering
\includegraphics[width=\linewidth]{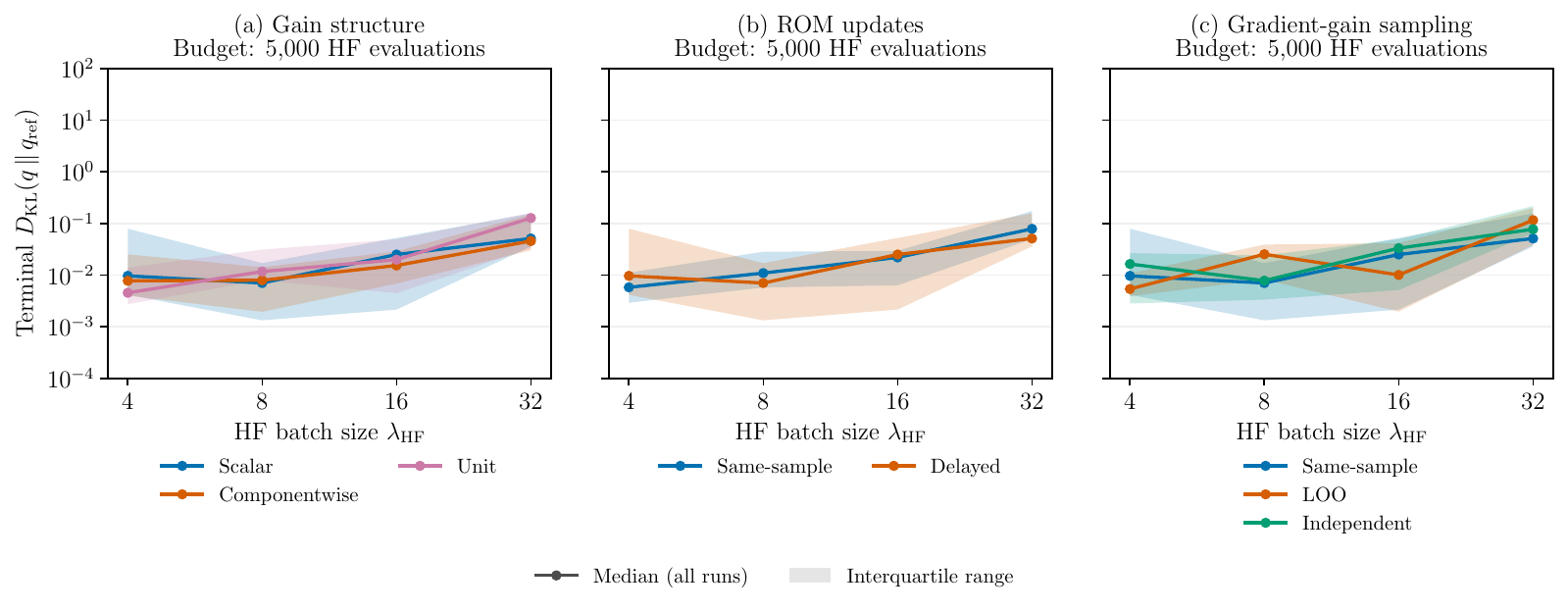}
\caption{CDR multifidelity studies: terminal KL to the numerical VI reference solution versus HF batch size for (a) gain structure, (b) ROM-update timing, and (c) gradient-gain sampling. Lines show medians, bands show interquartile ranges. Independent gradient-gain fitting in (c) incurs additional HF evaluations that are not accounted for.}
\label{fig:cdr_mf_ablation_summary}
\end{figure}

We next examine gain structure, the manner in which the LF model is updated, and sampling strategies for computing the gradient gain. The study is restricted to MF-BBVI (Newton). We use $\nSamplesHF\in\{4,8,16,32\}$, 512 additional LF samples, a budget of 5000 HF evaluations, and ten realizations per configuration. Figure~\ref{fig:cdr_mf_ablation_summary}(a--c) summarizes the terminal KL.%

For gain structure, we compare the scalar and componentwise gains with fixed unit-magnitude gains for both gradient and Hessian estimators. We employ same-sample gain estimation along with delayed LF updates. In Figure~\ref{fig:cdr_mf_ablation_summary}(a), we observe similar terminal medians at each HF batch size. Overall, we again observe that gains with more expressivity do not necessarily improve convergence due to finite sampling error. Interestingly, for this example we find that unit magnitude gains perform almost as well as both the scalar and componentwise gains. 

For LF model adaptation, Figure~\ref{fig:cdr_mf_ablation_summary}(b) compares same-sample and delayed updates with same-sample scalar gradient and Hessian gains. Same-sample updates yield a lower median at HF batch sizes of 4 and 16, while delayed updates yield a lower median at 8 and 32. Again, the results suggest that the bias incurred by the LF model update has a minimal result on the overall optimization problem.

Lastly, for the sampling strategy for computing the gradient gain, we compare same-sample, LOO, and independent scalar gains (Figure~\ref{fig:cdr_mf_ablation_summary}(c)). We employ delayed LF model adaptation and same-sample Hessian gains; the LOO and independent gradient gains therefore address both sources of bias discussed in Table~\ref{tab:adaptive_correlation_strategies}. Across all HF batch sizes, we again observe minimal statistical difference between the methods. The results again suggest that the bias incurred from fitting the gain has a minimal impact on the overall optimization. However, as independent sampling requires $\nSamplesHF$ additional HF samples at each iteration, the same-sample and LOO updates are clearly advantageous in terms of cost.  

\FloatBarrier

\subsection{Premixed \texorpdfstring{$\mathrm{H}_2$}{H2}-air flame}
\begin{figure}
\centering
\begin{tikzpicture}[font=\large, line cap=round, line join=round]

\def\W{9}   %
\def\H{4.5} %

\draw[line width=2.2pt] (0,0) rectangle (\W,\H);

\node at (1.2,3.85) {$\Gamma_{1}$};
\node at (1.2,2.25) {$\Gamma_{2}$};
\node at (1.2,0.65) {$\Gamma_{3}$};
\node at (4.7,0.55) {$\Gamma_{4}$};
\node at (8.2,2.35) {$\Gamma_{5}$};
\node at (4.5,3.85) {$\Gamma_{6}$};

\tikzset{
  dim/.style={line width=0.9pt},
  dimArrow/.style={<->,>=Stealth,dim},
  dimTick/.style={dim}
}

\draw[dim] (-0.45,0) -- (-0.45,\H);
\draw[dimTick] (-0.55,0) -- (-0.35,0);
\draw[dimTick] (-0.55,1.5) -- (-0.35,1.5);
\draw[dimTick] (-0.55,3.0) -- (-0.35,3.0);
\draw[dimTick] (-0.55,\H) -- (-0.35,\H);

\draw[dimArrow] (-1.05,3.0) -- (-1.05,\H)
  node[midway,left=3pt] {$3~\mathrm{mm}$};
\draw[dimArrow] (-1.05,1.5) -- (-1.05,3.0)
  node[midway,left=3pt] {$3~\mathrm{mm}$};
\draw[dimArrow] (-1.05,0) -- (-1.05,1.5)
  node[midway,left=3pt] {$3~\mathrm{mm}$};

\node[rotate=90] at (-3.0,2.25) {$x_{2}$};

\draw[dim] (\W+0.45,0) -- (\W+0.45,\H);
\draw[dimTick] (\W+0.35,0) -- (\W+0.55,0);
\draw[dimTick] (\W+0.35,\H) -- (\W+0.55,\H);
\draw[dimArrow] (\W+1.05,0) -- (\W+1.05,\H)
  node[midway,right=3pt] {$9~\mathrm{mm}$};

\draw[dim] (0,-0.45) -- (\W,-0.45);
\draw[dimTick] (0,-0.35) -- (0,-0.55);
\draw[dimTick] (\W,-0.35) -- (\W,-0.55);
\draw[dimArrow] (0,-1.05) -- (\W,-1.05)
  node[midway,below=3pt] {$18~\mathrm{mm}$};

\node at (\W/2,-1.95) {$x_{1}$};

\end{tikzpicture}
\caption{Premixed $\HTwo$ air-flame example. Schematic of geometry and boundary conditions.}
\label{fig:air_flame_schematic}
\end{figure}
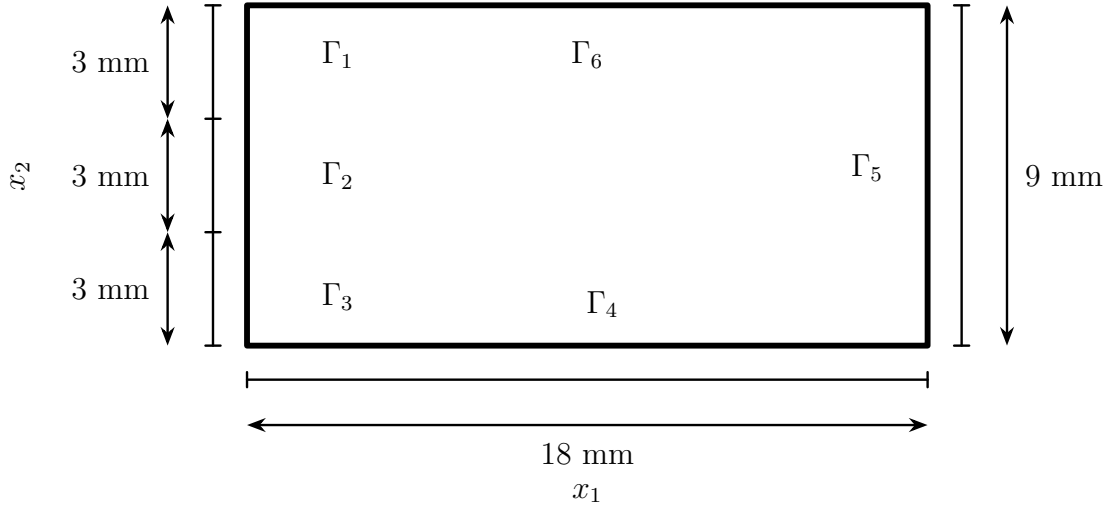

We next consider a two-dimensional premixed \(\mathrm{H}_2\)-air advection-diffusion-reaction system~\cite{LeCa20,GaWaZh20,BuWi10}. The governing equations are
\begin{equation}\label{eq:air-flame}
\begin{aligned}
\frac{\partial }{\partial t} \boldsymbol{u}(y,t,\params) + \boldsymbol{\beta} \cdot \nabla \boldsymbol{u} (y,t,\params) &= \kappa \Delta \boldsymbol{u} (y,t,\params) + \boldsymbol{R}(\boldsymbol{u}, \params), & \quad & y \in \Omega, \\
\boldsymbol{u} (y,t,\params) &= \boldsymbol{u_d}(y), & \quad & y \in \Gamma_d, \\
\boldsymbol{\nabla}  \boldsymbol u (y,t,\params) \cdot \boldsymbol{n} &= 0, & \quad & y \in \Gamma_n, \\
\boldsymbol{u}(y,0,\params) &= \boldsymbol{u_0}(y), & \quad & y \in \Omega. \\
\end{aligned}
\end{equation}
The spatial domain is $\Omega \defeq \left[0,L_x\right] \times \left[0,L_y \right]$ with $L_x = 1.8$ $\mathrm{cm}$ and $L_y = 0.90$ $\mathrm{cm}$, and the temporal domain is $t \in [0, 0.06 \mathrm{s}]$.
The state $\boldsymbol u \defeq \left[ \mathrm{Y}_{\mathrm{H}_2},\mathrm{Y}_{\mathrm{O}_2},\mathrm{Y}_{\mathrm{H}_2\mathrm{O}}, \Theta \right]$ contains the mass fractions of $\mathrm{H}_2$, $\mathrm{O}_2$, and $\mathrm{H}_{2}\mathrm{O}$, and the temperature, respectively. The velocity $\boldsymbol \beta\in\RR{2}$ has units $\mathrm{cm}\cdot \mathrm{s}^{-1}$, and the diffusivity $\kappa$ has units $\mathrm{cm}^{2}\cdot\mathrm{s}^{-1}$. The nonlinear reaction source term is
\begin{equation}
\begin{aligned}
    \boldsymbol R_i &= -\nu_i \left( \frac{W_i}{\rho} \right) \left( \frac{\rho \mathrm{Y}_{\mathrm{H}_2}}{\mathrm{W}_{\mathrm{H}_2}}  \right)^{\nu_{\mathrm{H}_2}} \left( \frac{\rho \mathrm{Y}_{\mathrm{O}_2}}{\mathrm{W}_{\mathrm{O}_2}}  \right)^{\mathrm{\nu}_{\mathrm{O}_2}} A \exp \left( -\frac{E}{R_{\mathrm{gas}} \Theta} \right), \\
\boldsymbol R_T &= Q R_{\HTwoO}.
\end{aligned}
\end{equation}
Here, $i\in \{\HTwo,\OTwo,\HTwoO\}$ indexes the species,  
$\nu_{\HTwo},\nu_{\OTwo},\nu_{\HTwoO} = 2,1,-2$ are the stoichiometric coefficients, $W_{\HTwo}, W_{\OTwo}, W_{\HTwoO} = 2.016, 31.9, 18$ are the molecular weights with units $\mathrm{g} \cdot \mathrm{mol}^{-1}$, $\rho = 1.39 \times 10^{-3}$ $\mathrm{g} \cdot \mathrm{cm}^{-3}$ is the density, $R_{\mathrm{gas}} = 8.314$ $\mathrm{J} \cdot  \mathrm{mol}^{-1} \cdot \mathrm{K}^{-1}$ is the universal gas constant, and $Q = 9800 \; \mathrm{K}$ is the heat of reaction. The pre-exponential factor is set to $A=4.465 \times 10^{12}$. 
As depicted in Figure~\ref{fig:air_flame_schematic}, the domain boundary $\partial \Omega$ is decomposed into six segments, $\partial \Omega = \cup_{i=1}^6 \Gamma_i$, with a Dirichlet part $\Gamma_d = \cup_{i=1}^3 \Gamma_i$ and a Neumann part $\Gamma_n = \cup_{i=4}^6 \Gamma_i$. The initial interior state is $\boldsymbol u(y,0,\params) = \left[0,0,0,300 \right]^\intercal$. The boundary conditions on $\Gamma_d$ are
\[\boldsymbol u_d(y) = \begin{cases}
\left[0,0,0,300 \right]^\intercal \qquad &y \in \Gamma_1 \cup \Gamma_3, \\
\left[0.0282,0.2259,0,950 \right]^\intercal \qquad &y  \in \Gamma_2. \\
\end{cases}
\] 
We discretize Eq.~\eqref{eq:air-flame} on a $(n_x,n_y)=(64,32)$ nodal grid. We use second-order upwinding in $x$, first-order upwinding in $y$, and second-order centered differences for diffusion. Time integration uses the Crank--Nicolson method with $\Delta t=10^{-3}\,\mathrm{s}$.

\begin{figure}
\centering
\includegraphics[width=\linewidth]{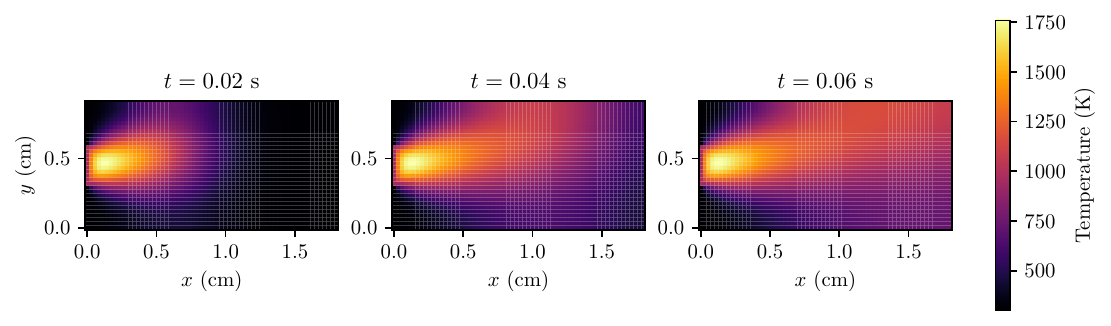}
\caption{H$_2$-air flame: temperature at the generating parameters at three times.}
\label{fig:flame_solution}
\end{figure}

We infer the four parameters \(\params=(\kappa,E/10^3,\beta_1,\beta_2)\), where \(E/10^3\) is the activation energy in units of \(10^3\,\mathrm{J}\,\mathrm{mol}^{-1}\). We set the parameter bounds to \([0.5,4]\), \([4,12]\), \([20,60]\), and \([1,20]\), respectively. Synthetic observations are generated at
\[
\params^\star=(2,8,40,7).
\]
As observations, we collect temperature data every fourth grid point in each spatial direction and every fourth stored snapshot in time. Snapshots are stored every ten time steps. We add synthetic noise to the data according to%
\begin{equation}
\mathbf y_{\rm obs}=\mathbf y^\star+\boldsymbol\varepsilon,\qquad
\varepsilon_i\overset{\rm iid}{\sim}\mathcal N(0,\sigma_{\rm obs}^2),\qquad
\sigma_{\rm obs}=0.05\frac{\|\mathbf y^\star\|_2}{\sqrt{N_{\rm obs}}},
\end{equation}
where \(\mathbf y^\star\in\mathbb R^{N_{\rm obs}}\) is the noiseless observation vector.
The Gaussian prior and the initial variational distribution are set to have the same mean \((2.525,9.2,46,9.85)\) and standard deviations \((0.37875,1.38,6.9,1.4775)\). All methods use the same bounded sigmoid parameter transformation and a diagonal-Gaussian variational family in the transformed coordinates.
Figure~\ref{fig:flame_solution} shows the temperature field at the generating parameters.

As a ``truth'' comparison, we again compare against a high-fidelity VI reference solution. We use the saved distribution after 225 Newton VI updates with 512 samples per iteration, corresponding to 115,712 HF evaluations including the initial batch. The reference is initialized at the previously computed MAP point, with standard deviations equal to 5\% of its component magnitudes, and uses fixed steps without a line search.

This premixed air-flame example incurs a higher computational cost than the previous two exemplars. As the previous two experiments suggested that the gain structure, ROM update timing, and gain sampling strategy have a relatively small impact on the overall optimization, we do not repeat these studies here. All MF methods use the scalar gain, delayed ROM updates, and LOO gain sampling. 

\subsubsection{Method comparison}
We now compare the various methods with a fixed computational budget. Each method uses four HF samples per batch and a budget of 1000 HF evaluations; the multifidelity methods use 512 additional ROM samples. We repeat each method four times to account for stochastic optimization. 

\begin{figure}
\centering
\includegraphics[width=0.92\linewidth]{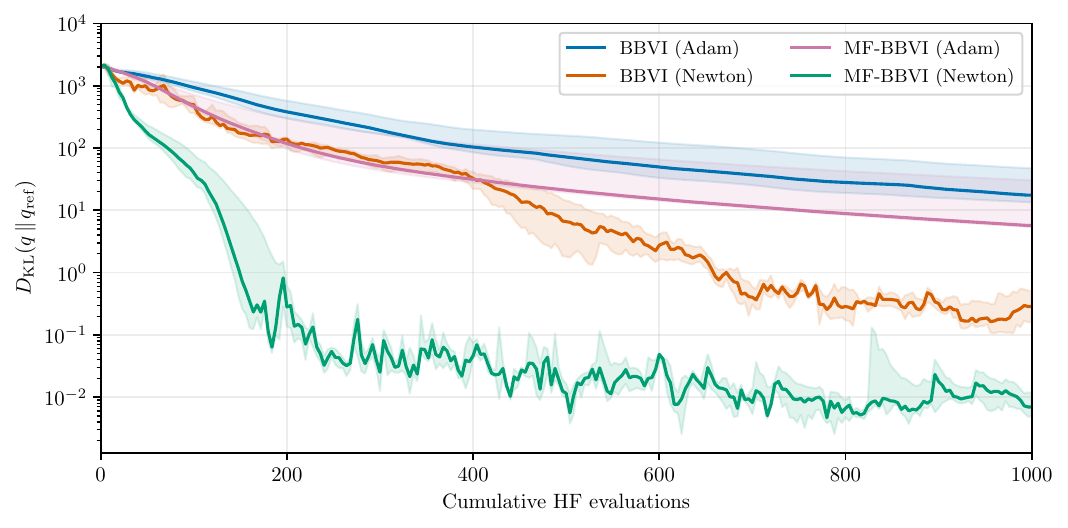}
\caption{H$_2$-air flame method comparison: KL to the numerical VI reference versus cumulative HF evaluations. Curves and bands show medians and interquartile ranges.}
\label{fig:air_flame_kl}
\end{figure}

Figure~\ref{fig:air_flame_kl} shows KL convergence, and Table~\ref{tab:flame_method_ablation} reports terminal medians across four seeds. MF-BBVI (Newton) attains the smallest median KL, $6.97 \times 10^{-3}$. BBVI (Newton) is again the next best performing method, with a terminal median KL of $0.285$. BBVI (Adam) is again the worst performing method, with a terminal median KL of $17.5$; MF-BBVI (Adam) again improves upon this with a terminal median KL of $5.70$. Overall, multifidelity sampling again improves both optimizers, and the Newton optimizer again outperforms Adam.

\begin{table}[!htbp]
\centering
\caption{H$_2$-air method comparison: median KL at the last recorded iteration across 4 seeds. Green and red mark the smallest and largest median KL, respectively.}
\label{tab:flame_method_ablation}
\begin{tabular}{lr}
\hline
Method & Median KL \\
\hline
BBVI (Adam) & \textcolor{red}{$17.5$} \\
BBVI (Newton) & $0.285$ \\
MF-BBVI (Adam) & $5.7$ \\
MF-BBVI (Newton) & \textcolor{bestkl}{$6.97\times10^{-3}$} \\
\hline
\end{tabular}
\end{table}

\begin{figure}
\centering
\includegraphics[width=0.92\linewidth]{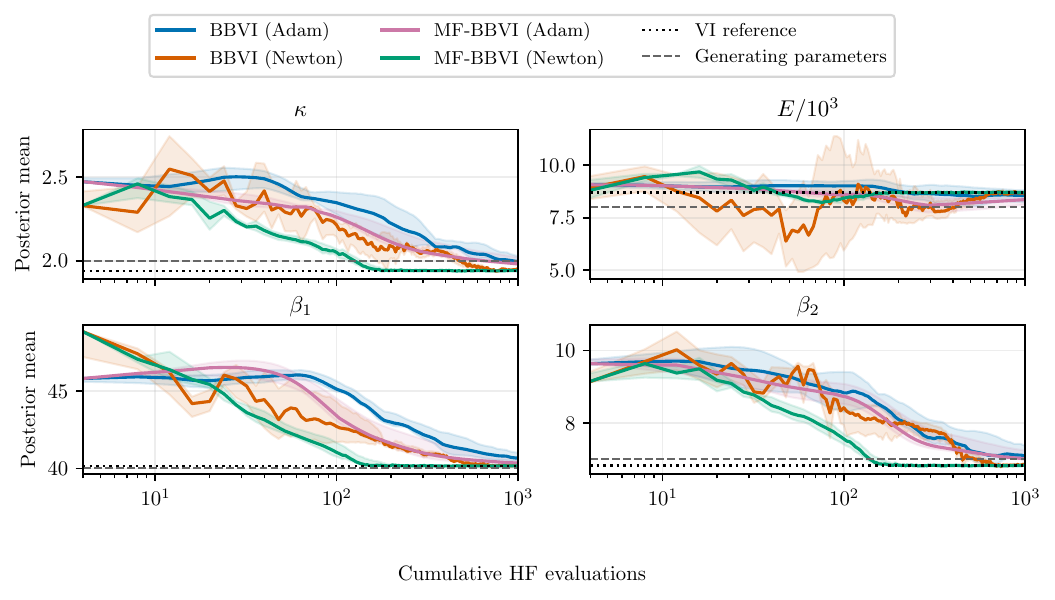}
\includegraphics[width=0.92\linewidth]{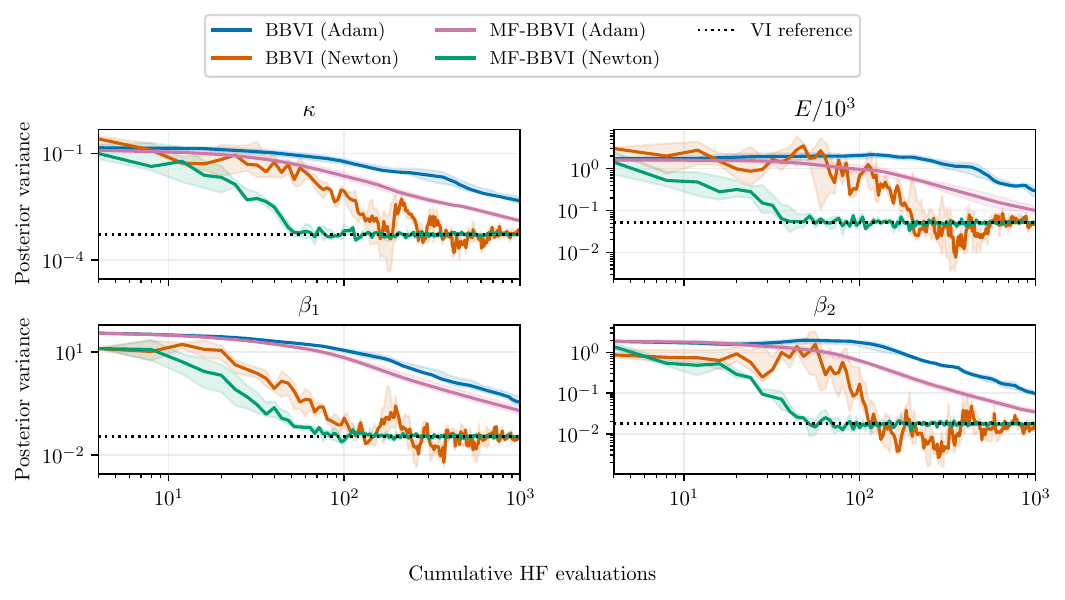}
\caption{H$_2$-air flame method comparison: physical posterior means (top) and marginal variances (bottom). Curves show medians and interquartile ranges. Black dotted lines mark the VI reference moments; gray dashed lines in the mean panels mark the generating parameters.}
\label{fig:flame_moments}
\end{figure}

\begin{figure}
\centering
\includegraphics[width=\linewidth]{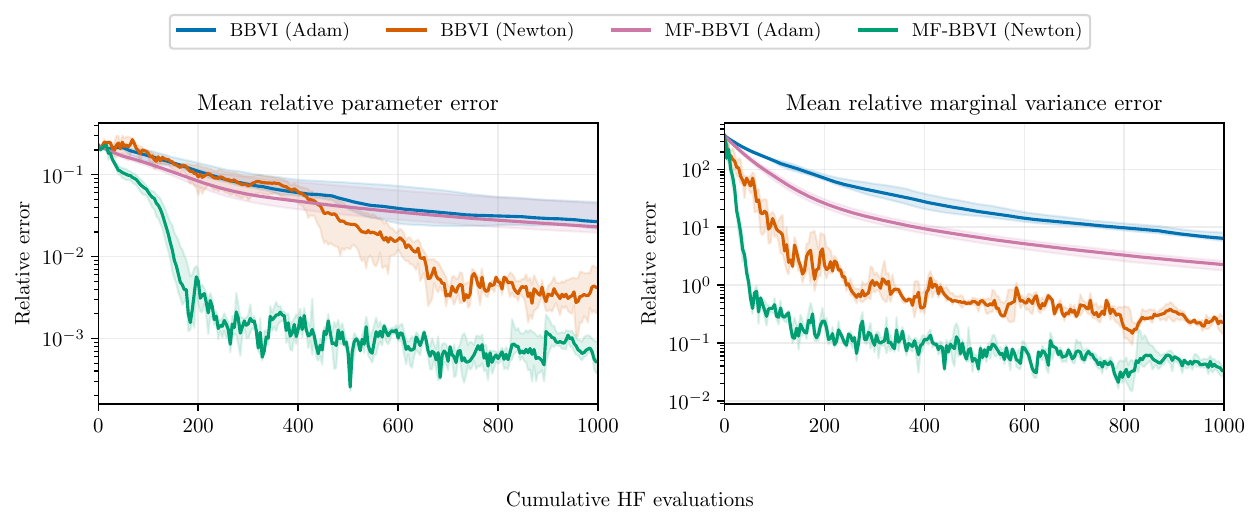}
\caption{H$_2$-air flame method comparison: mean relative parameter error (left) and mean relative marginal variance error (right) with respect to the numerical VI reference. Curves and bands show medians and interquartile ranges.}
\label{fig:flame_relative_moment_errors}
\end{figure}

Figure~\ref{fig:flame_moments} shows the physical mean and variance histories, while Figure~\ref{fig:flame_relative_moment_errors} summarizes agreement with the numerical VI reference solution. We observe that MF-BBVI (Newton) well-approximates the reference posterior for all parameters within two hundred iterations. BBVI (Newton) is the next best performing method, followed by MF-BBVI (Adam). MF sampling improves both optimizers, while the Newton-type optimizer outperforms Adam.  

\subsubsection{Sample-size study}
\begin{figure}
\centering
\includegraphics[width=\linewidth]{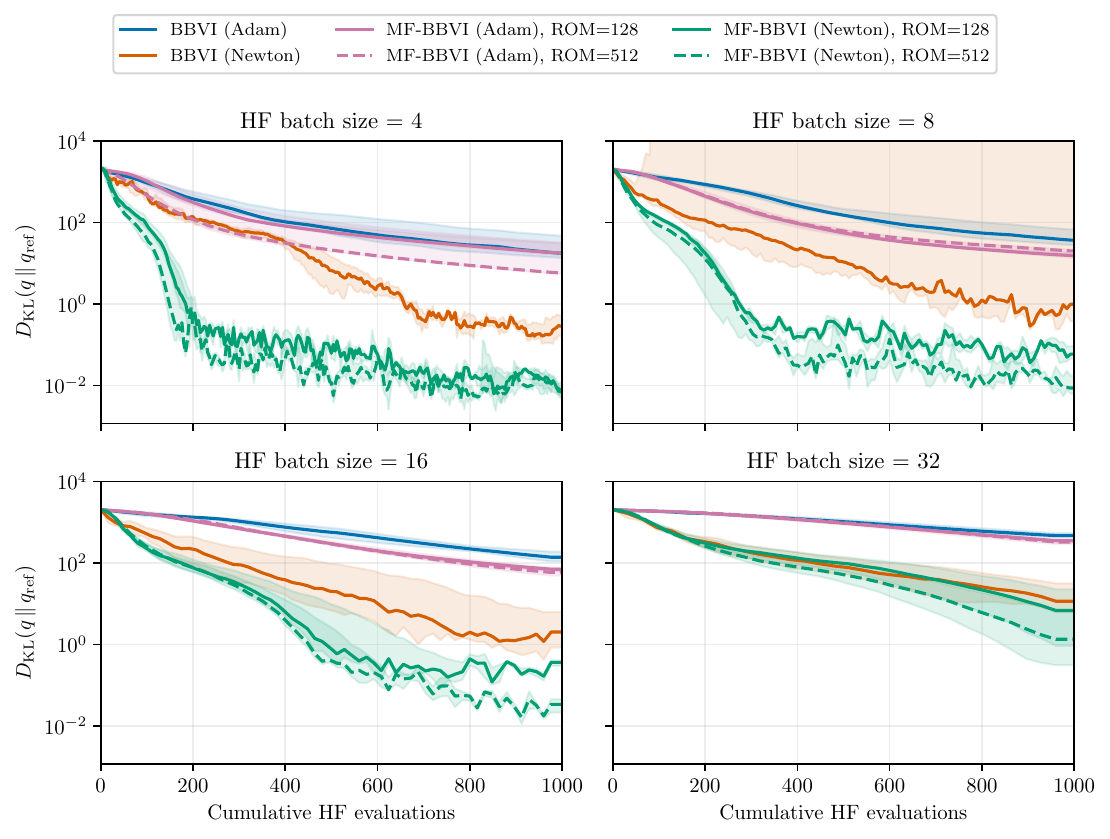}
\caption{H$_2$-air flame sample-size study: KL to the numerical VI reference versus cumulative HF evaluations. Curves and bands show medians and interquartile ranges.}
\label{fig:flame_sample_size_ablation}
\end{figure}

We now vary the HF batch size over $\{4,8,16,32\}$ and the additional ROM sample count over $\{128,512\}$. We retain the 1000-evaluation budget and the optimizer and surrogate settings used in the method comparison. Figure~\ref{fig:flame_sample_size_ablation} shows KL convergence to the same VI reference solution, and Table~\ref{tab:flame_sample_size} reports terminal medians across four seeds.

MF-BBVI (Newton) has the smallest terminal median for every configuration, giving, on average, an order of magnitude improvement compared to BBVI (Newton). With the exception of a single Adam configuration, we observe that MF sampling improves every configuration. We additionally observe that increasing the ROM sample size improves every configuration, and that Newton outperforms Adam. We again observe a non-monotonic improvement as the HF batch size is increased for a fixed computational budget. 

\begin{table}[!htbp]
\centering
\small
\caption{H$_2$-air sample-size study: median terminal KL across 4 seeds. ROM denotes additional low-fidelity samples for the multifidelity methods. Single-fidelity values are repeated for comparison. Green and red mark the smallest and largest KL in each row.}
\label{tab:flame_sample_size}
\begin{tabular}{rr|rrrr}
\hline
 & & BBVI & BBVI & MF-BBVI & MF-BBVI \\
HF batch & ROM & (Adam) & (Newton) & (Adam) & (Newton) \\
\hline
4 & 128 & $17.5$ & $0.285$ & \textcolor{red}{$17.6$} & \textcolor{bestkl}{$7.72\times10^{-3}$} \\
4 & 512 & \textcolor{red}{$17.5$} & $0.285$ & $5.7$ & \textcolor{bestkl}{$6.97\times10^{-3}$} \\
8 & 128 & \textcolor{red}{$37$} & $0.988$ & $15.4$ & \textcolor{bestkl}{$0.0586$} \\
8 & 512 & \textcolor{red}{$37$} & $0.988$ & $20$ & \textcolor{bestkl}{$8.64\times10^{-3}$} \\
16 & 128 & \textcolor{red}{$138$} & $2.02$ & $70.4$ & \textcolor{bestkl}{$0.363$} \\
16 & 512 & \textcolor{red}{$138$} & $2.02$ & $58.8$ & \textcolor{bestkl}{$0.0338$} \\
32 & 128 & \textcolor{red}{$468$} & $11.5$ & $352$ & \textcolor{bestkl}{$6.82$} \\
32 & 512 & \textcolor{red}{$468$} & $11.5$ & $328$ & \textcolor{bestkl}{$1.34$} \\
\hline
\end{tabular}
\end{table}

\FloatBarrier
\subsection{Shallow-water bathymetry inverse problem}
\label{sec:swe_bathymetry}
\begin{figure}
\centering
\includegraphics[width=\linewidth]{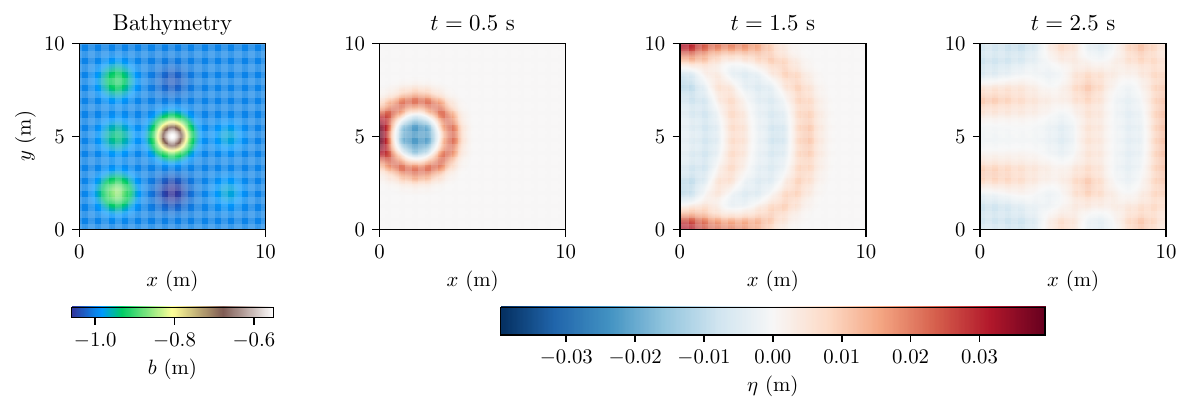}
\caption{Shallow-water problem: generating bathymetry and free-surface elevation at three times.}
\label{fig:swe_solution}
\end{figure}

Finally, we consider a two-dimensional inverse problem governed by the shallow-water equations. The governing equations are given on the domain $\Omega=[0,10]\times[0,10]\,\mathrm{m}^2$ as
\begin{equation}
 \frac{\partial h}{\partial t}+\nabla\cdot(h\mathbf u)=0,
 \qquad
 \frac{\partial(h\mathbf u)}{\partial t}+\nabla\cdot\left(h\mathbf u\otimes\mathbf u+\frac12gh^2I\right)=-gh\nabla b,
\end{equation}
where $h$ is the water depth, $\mathbf u$ is the depth-averaged horizontal velocity, $g=9.81\,\mathrm{m}\,\mathrm{s}^{-2}$, and $b$ is the bathymetry. The free-surface elevation is $\eta=h+b$. The unknown bathymetry is represented by nine Gaussian basis functions,
\begin{equation}
 b(x,y;\params)=-1+\sum_{i=1}^{9}\theta_i\exp\left[-\frac{(x-c_{x,i})^2+(y-c_{y,i})^2}{2\ell_b^2}\right],
 \qquad \ell_b=0.6\,\mathrm{m},
\end{equation}
with centers on the $3\times3$ tensor grid $c_x,c_y\in\{2,5,8\}\,\mathrm{m}$, ordered first in $x$ and then in $y$. The generating coefficients are
\begin{equation}
 \params^\star=(0.15,-0.06,0.03,0.06,0.45,0.03,0.09,-0.03,0)^\top\,\mathrm{m}.
\end{equation}
The initial velocity is zero, and the initial free surface is
\begin{equation}
 \eta(x,y,0)=0.15\exp\left[-\frac{(x-2)^2+(y-5)^2}{2(0.45)^2}\right]\,\mathrm{m}.
\end{equation}
Reflective slip-wall conditions are imposed on all boundaries. Figure~\ref{fig:swe_solution} shows the bathymetry and free-surface evolution.

The HF model uses a $200\times200$ Cartesian finite-volume mesh, hydrostatic reconstruction, Rusanov fluxes, and second-order SSP Runge--Kutta time integration with $\Delta t=2.5\times10^{-3}\,\mathrm{s}$ over $t\in[0,2.5]$ s. We observe $\eta$ at every twentieth cell center in both directions, giving 100 spatial locations, and at 21 times $t=0.5,0.6,\ldots,2.5$ s; the resulting observation map has 2100 components. We add independent Gaussian noise to the  observations with a standard deviation of $5\times10^{-4}\,\mathrm{m}$. The prior is $\mathcal N(0,0.15^2I)$, which is then restricted to the coefficient bounds $[-0.8,0.8]$ m.  We again compare to a high-fidelity reference ``truth'' solution. We obtain this solution via a VI Newton run employing $\nSamplesHF=512$ samples per iteration for 102 optimization iterations. The reference run employs a zero mean initial condition with a standard deviation of $0.15$.  

\subsubsection{Method comparison}
\begin{figure}
\centering
\includegraphics[width=0.92\linewidth]{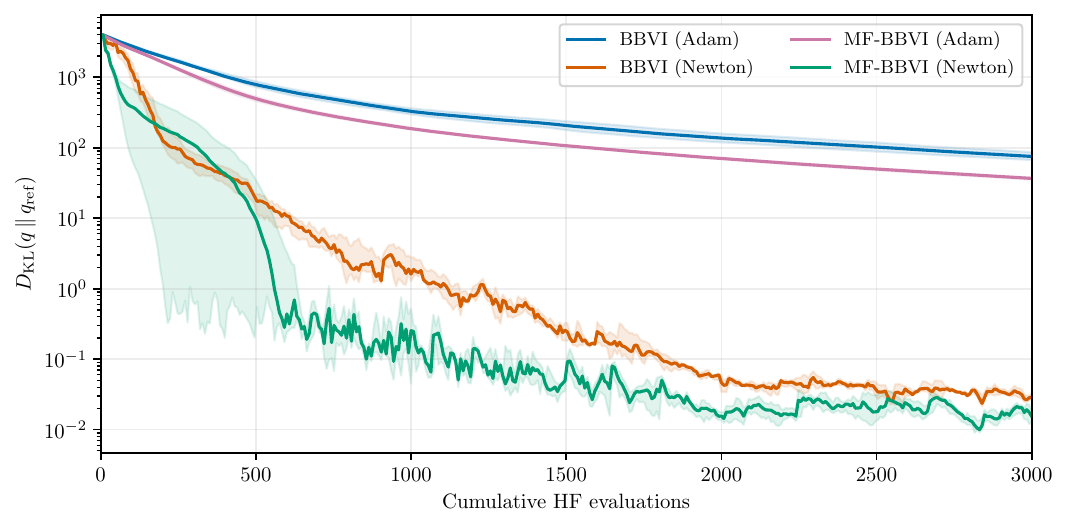}
\caption{Shallow-water method comparison: KL to the numerical VI reference versus cumulative HF evaluations. Curves show medians across four seeds, with interquartile bands.}
\label{fig:swe_kl}
\end{figure}

\begin{figure}
\centering
\includegraphics[width=0.95\linewidth]{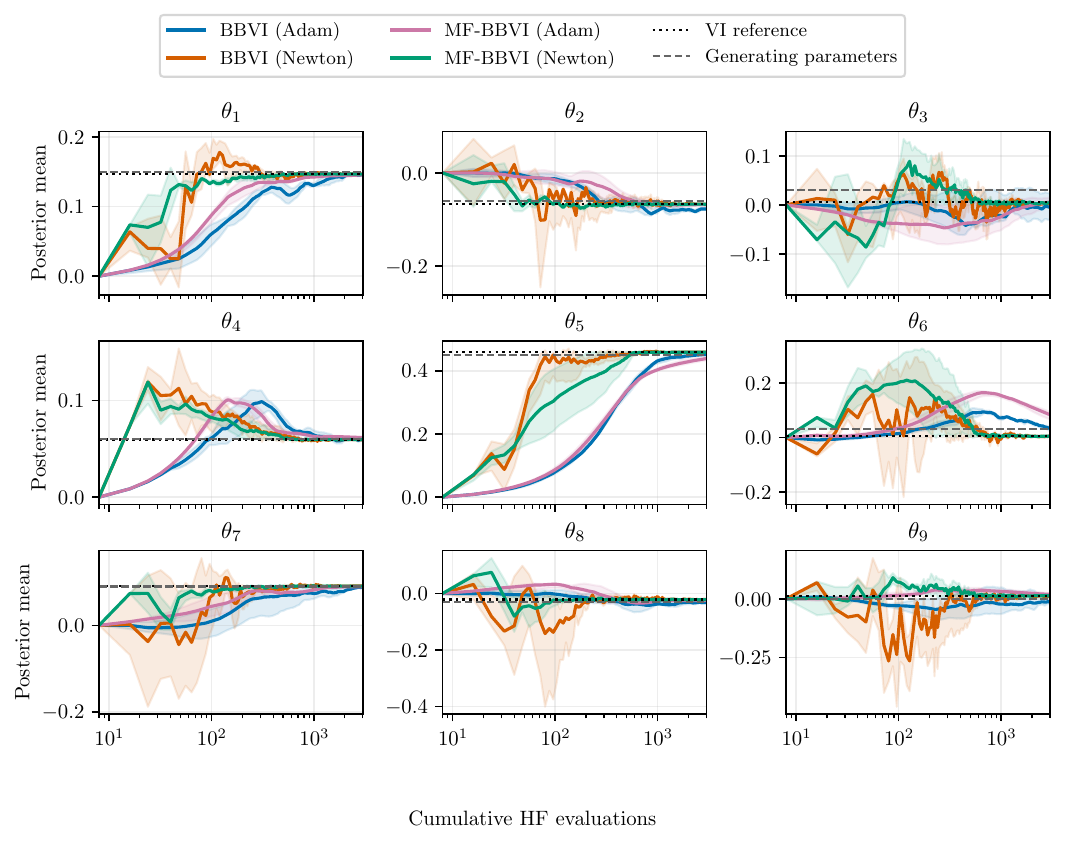}
\caption{Shallow-water method comparison: physical posterior means. Curves show medians and interquartile ranges. Black dotted lines mark the VI reference means; gray dashed lines mark the generating parameters.}
\label{fig:swe_moments}
\end{figure}

\begin{figure}
\centering
\includegraphics[width=0.95\linewidth]{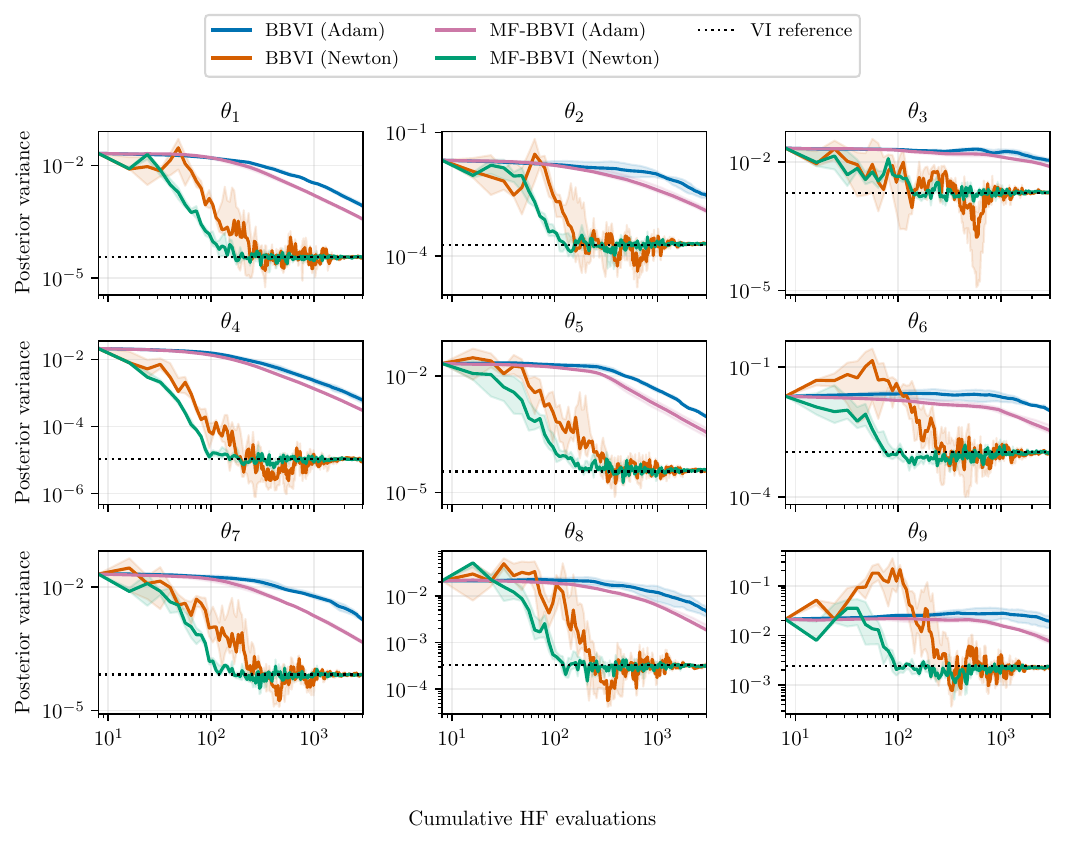}
\caption{Shallow-water method comparison: physical marginal variances. Curves show medians and interquartile ranges. Black dotted lines mark the VI reference variances.}
\label{fig:swe_moments_variance}
\end{figure}

Each method uses 8 HF samples per batch and a nominal budget of 3000 HF evaluations, with four realizations per method. The multifidelity methods use 512 additional ROM samples. Figure~\ref{fig:swe_kl} shows KL convergence, and Table~\ref{tab:swe_method_ablation} reports terminal medians. Once again, we observe that MF-BBVI (Newton) has the smallest median KL, followed by BBVI (Newton). MF-BBVI (Adam) again outperforms BBVI (Adam), but neither are as accurate as the Newton-based solutions. Thus, multifidelity sampling again lowers the median KL for both optimizers, and Newton optimization outperforms Adam. Figures~\ref{fig:swe_moments} and~\ref{fig:swe_moments_variance} show the physical mean and marginal variance histories for all nine coefficients. Both Newton methods closely approximate the reference means and marginal variances, with MF-BBVI (Newton) giving the smaller aggregate errors. The Adam methods retain larger mean errors and substantially overestimate the marginal variances within the specified computational budget.

\begin{table}[!htbp]
\centering
\caption{SWE method comparison: median KL at the last recorded iteration across 4 seeds. Green and red mark the smallest and largest median KL, respectively.}
\label{tab:swe_method_ablation}
\begin{tabular}{lr}
\hline
Method & Median KL \\
\hline
BBVI (Adam) & \textcolor{red}{$75.3$} \\
BBVI (Newton) & $0.0289$ \\
MF-BBVI (Adam) & $36.7$ \\
MF-BBVI (Newton) & \textcolor{bestkl}{$0.0154$} \\
\hline
\end{tabular}
\end{table}

\begin{figure}
\centering
\includegraphics[width=\linewidth]{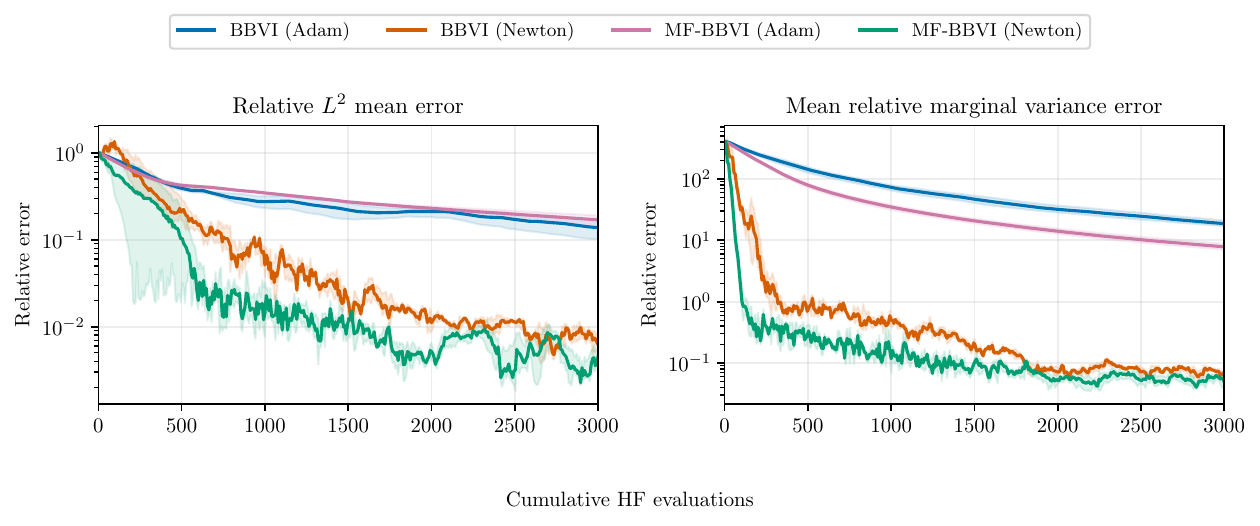}
\caption{Shallow-water method comparison: relative $L^2$ error of the posterior mean (left) and mean componentwise relative marginal variance error (right) with respect to the numerical VI reference. Curves show medians across four seeds, with interquartile bands.}
\label{fig:swe_relative_moment_errors}
\end{figure}

Figure~\ref{fig:swe_relative_moment_errors} summarizes agreement with the numerical VI reference. We measure the posterior mean error by $\|\boldsymbol\mu^{\mathrm{phys}}_k-\boldsymbol\mu^{\mathrm{phys}}_{\mathrm{ref}}\|_2/\|\boldsymbol\mu^{\mathrm{phys}}_{\mathrm{ref}}\|_2$, avoiding division by individual near-zero reference coefficients. Relative marginal variance errors are averaged across parameters within each run; both metrics are then summarized across the four runs. The terminal median mean error is $0.448\%$ for MF-BBVI (Newton), compared with $0.627\%$ for BBVI (Newton), $13.8\%$ for BBVI (Adam), and $17.0\%$ for MF-BBVI (Adam). The corresponding median relative variance errors are $4.92\%$, $7.09\%$, $1871\%$, and $784\%$, respectively. Both Newton methods therefore closely recover the reference moments, with multifidelity sampling improving both error measures. For Adam, multifidelity sampling reduces the variance error and KL despite a larger mean error, illustrating that improved KL need not imply improved agreement of every moment.

\section{Conclusions}
\label{sec:conclusions}
We presented a multifidelity black-box variational inference framework for PDE-constrained Bayesian inverse problems in settings where PDE gradients are unavailable. The approach combines score-function gradient and Hessian estimators with approximate control variates utilizing adaptive reduced-order models to reduce the cost of variational inference. We further outlined a PDE-derivative-free, regularized diagonal-Newton optimization strategy based on the gradient and Hessian estimators. We investigated several strategies for adapting the ROM and estimating the approximate control variate gain, and examined the tradeoff between low-variance and unbiased estimators. Numerical results across four benchmark problems with up to nine parameters demonstrated that the proposed multifidelity approach yields substantial improvements over single-fidelity score-function BBVI using the same optimizers.

The present contribution introduces various promising directions for future work. First, we believe that more sample-efficient schemes can be developed by combining the developed multifidelity estimators with adaptive sampling and sample-reuse approaches~\cite{rei2025blackboxvariationalinference}. Second, we believe that further improvements could be realized via an inference approach with ROM substeps, wherein the ROM is used to perform multiple optimization steps before sampling the high-fidelity model. Third, we believe there are numerous opportunities to improve the Newton-like optimization scheme outlined here. We employ a fairly simple regularized Newton optimization strategy. The development of a stochastic trust region approach with adaptive regularization is one such promising future direction. Lastly, the present work primarily considered the mean-field variational family. The extension of the approach to richer variational families, including multivariate Gaussian and neural-network-based variational densities, will enable more accurate representation of correlated and non-Gaussian posteriors. 

\section*{Acknowledgements}
Sandia National Laboratories is a multi-mission
laboratory managed and operated by National Technology and Engineering
 Solutions of Sandia, LLC (NTESS), a wholly owned subsidiary
of Honeywell International, Inc., for the U.S. Department of Energy's
National Nuclear Security Administration (DOE/NNSA) under Contract
 DE-NA0003525. This written work is authored by an employee of
NTESS. The employee, not NTESS, owns the right, title, and interest
in and to the written work and is responsible for its contents. Any subjective
 views or opinions that might be expressed in the written work do
not necessarily represent the views of the U.S. Government. The publisher
 acknowledges that the U.S. Government retains a non-exclusive,
paid-up, irrevocable, world-wide license to publish or reproduce the
published form of this written work or allow others to do so, for U.S.
Government purposes. The DOE will provide public access to results of
federally sponsored research in accordance with the DOE Public Access
Plan.

\section*{Declaration of generative AI and AI-assisted technologies in the writing process}
During the preparation of this work, the author used generative AI tools for manuscript formatting, assistance with table generation, general manuscript editing and review, and code documentation. After using these tools, the author reviewed and edited the content as needed and takes full responsibility for the content of the publication.

\appendix

\section{Derivation of optimal control-variate gains}
\label{sec:appendix_optimal_gains}
This appendix derives the optimal control-variate and approximate control variate gains. We refer the reader to Ref.~\cite{petersen2012matrixcookbook} for several useful matrix-calculus identities. 

\paragraph{Standard control variates}
We consider the control-variate estimator
\[
X_{\mathrm{cv}}(\mathbf{A})=X+\mathbf{A}(Y-\mathbb{E}[Y]).
\]
Its covariance is
\[
\mathrm{Cov}[X_{\mathrm{cv}}(\mathbf{A})]
=
\mathbf{\Sigma}_{XX}+\mathbf{A}\mathbf{\Sigma}_{YY}\mathbf{A}^\top+\mathbf{A}\mathbf{\Sigma}_{YX}+\mathbf{\Sigma}_{XY}\mathbf{A}^\top,
\]
where
\[
\mathbf{\Sigma}_{XX}:=\mathrm{Cov}(X),\qquad
\mathbf{\Sigma}_{YY}:=\mathrm{Cov}(Y),\qquad
\mathbf{\Sigma}_{XY}:=\mathrm{Cov}(X,Y),\qquad
\mathbf{\Sigma}_{YX}:=\mathbf{\Sigma}_{XY}^\top.
\]
We define the objective to be the trace of the covariance,
\[
\mathcal{J}_{\mathrm{cv}}(\mathbf{A}):=\mathrm{tr}(\mathrm{Cov}[X_{\mathrm{cv}}(\mathbf{A})]).
\]
The above can be expanded as
\[
\mathcal{J}_{\mathrm{cv}}(\mathbf{A})
=
\mathrm{tr}(\mathbf{\Sigma}_{XX})
+\mathrm{tr}(\mathbf{A}\mathbf{\Sigma}_{YY}\mathbf{A}^\top)
+2\,\mathrm{tr}(\mathbf{A}\mathbf{\Sigma}_{YX}),
\]
using \(\mathrm{tr}(\mathbf{A}\mathbf{\Sigma}_{YX})=\mathrm{tr}(\mathbf{\Sigma}_{XY}\mathbf{A}^\top)\).

\paragraph{Full matrix gain}
For a full matrix gain, we make use of \(\nabla_{\mathbf{A}}\,\mathrm{tr}(\mathbf{A}\mathbf{\Sigma}_{YY}\mathbf{A}^\top)=2\mathbf{A}\mathbf{\Sigma}_{YY}\) and \(\nabla_{\mathbf{A}}\,\mathrm{tr}(\mathbf{A}\mathbf{\Sigma}_{YX})=\mathbf{\Sigma}_{YX}^\top=\mathbf{\Sigma}_{XY}\) to write,
\[
\nabla_{\mathbf{A}} \mathcal{J}_{\mathrm{cv}}(\mathbf{A})=2\mathbf{A}\mathbf{\Sigma}_{YY}+2\mathbf{\Sigma}_{XY}.
\]
Setting this derivative to zero gives
\[
\mathbf{A}^\star=-\,\mathbf{\Sigma}_{XY}\mathbf{\Sigma}_{YY}^{-1},
\]
which is~\eqref{eq:generic_control_variate_optimal_matrix}.

\paragraph{Diagonal gain}
In the case where we restrict to a component-wise control variate gain, we set \(\mathbf{A}=\mathbf{D}=\mathrm{diag}(a_1,\ldots,a_d)\),
\[
\mathcal{J}_{\mathrm{cv}}(\mathbf{D})
=
\mathrm{tr}(\mathbf{\Sigma}_{XX})
+\sum_{i=1}^{d}\left(a_i^2(\mathbf{\Sigma}_{YY})_{ii}+2a_i(\mathbf{\Sigma}_{XY})_{ii}\right).
\]
Each component decouples, and setting the gradient to be zero yields
\[
\frac{\partial \mathcal{J}_{\mathrm{cv}}}{\partial a_i}
=
2a_i(\mathbf{\Sigma}_{YY})_{ii}+2(\mathbf{\Sigma}_{XY})_{ii}=0
\quad\Rightarrow\quad
a_i^\star=-\,\frac{(\mathbf{\Sigma}_{XY})_{ii}}{(\mathbf{\Sigma}_{YY})_{ii}},
\]
which is the componentwise formula in~\eqref{eq:compwise}.

\paragraph{Approximate control variates}
We now turn to the multifidelity Monte Carlo approximate-control-variate estimator. For simplicity we define 
\[
\beta := \left(\frac{1}{\nSamplesHF}-\frac{1}{\nSamplesLF}\right)>0.
\]
The MFMC ACV estimator is
\[
\widehat X_{\mathrm{acv}}(\mathbf{A})
=
\left\langle X\right\rangle_{\mathcal{S}^{\mathrm{HF}}}
+
\mathbf{A}\left(
\left\langle Y\right\rangle_{\mathcal{S}^{\mathrm{HF}}}
-\left\langle Y\right\rangle_{\mathcal{S}^{\mathrm{LF}}}
\right).
\]
With this convention, the ACV gains have the same sign as the additive control-variate gains in Section~\ref{sec:vi}.
We define the means as
\[
\bar X_{\mathrm{HF}}:=\left\langle X\right\rangle_{\mathcal{S}^{\mathrm{HF}}},\qquad
\bar Y_{\mathrm{HF}}:=\left\langle Y\right\rangle_{\mathcal{S}^{\mathrm{HF}}},\qquad
\bar Y_{\mathrm{LF}}:=\left\langle Y\right\rangle_{\mathcal{S}^{\mathrm{LF}}}.
\]
Since $\sampleSetHighFidelity\subseteq\sampleSetLowFidelity$ and the samples are i.i.d., we may write
$\bar Y_{\mathrm{LF}}=\frac{\nSamplesHF}{\nSamplesLF}\bar Y_{\mathrm{HF}}+\left(1-\frac{\nSamplesHF}{\nSamplesLF}\right)\bar Y_{\mathrm{e}}$,
where $\bar Y_{\mathrm{e}}$ is the average of $Y$ over $\sampleSetLowFidelity\setminus\sampleSetHighFidelity$ and is independent of $\bar X_{\mathrm{HF}}$ and $\bar Y_{\mathrm{HF}}$. Hence
\begin{align*}
\mathrm{Cov}[\bar X_{\mathrm{HF}},\bar Y_{\mathrm{HF}}-\bar Y_{\mathrm{LF}}]
&=\left(1-\frac{\nSamplesHF}{\nSamplesLF}\right)\frac{1}{\nSamplesHF}\mathbf{\Sigma}_{XY}
=\beta\,\mathbf{\Sigma}_{XY},
\\
\mathrm{Cov}[\bar Y_{\mathrm{HF}}-\bar Y_{\mathrm{LF}}]
&=\left(1-\frac{\nSamplesHF}{\nSamplesLF}\right)^2\left(\frac{1}{\nSamplesHF}+\frac{1}{\nSamplesLF-\nSamplesHF}\right)\mathbf{\Sigma}_{YY}
=\beta\,\mathbf{\Sigma}_{YY},
\end{align*}
and $\mathrm{Cov}[\bar X_{\mathrm{HF}}]=\mathbf{\Sigma}_{XX}/\nSamplesHF$. The covariance of the estimator is therefore
\begin{align*}
\mathrm{Cov}\!\left[\widehat X_{\mathrm{acv}}(\mathbf{A})\right]
&=
\mathrm{Cov}[\bar X_{\mathrm{HF}}]
+\mathbf{A}\,\mathrm{Cov}[\bar Y_{\mathrm{HF}}-\bar Y_{\mathrm{LF}}]\,\mathbf{A}^\top\\
&\quad
+\mathbf{A}\,\mathrm{Cov}[\bar Y_{\mathrm{HF}}-\bar Y_{\mathrm{LF}},\bar X_{\mathrm{HF}}]
+\mathrm{Cov}[\bar X_{\mathrm{HF}},\bar Y_{\mathrm{HF}}-\bar Y_{\mathrm{LF}}]\mathbf{A}^\top\\
&=
\frac{1}{\nSamplesHF}\mathbf{\Sigma}_{XX}
+\beta\left(
\mathbf{A}\mathbf{\Sigma}_{YY}\mathbf{A}^\top
+\mathbf{A}\mathbf{\Sigma}_{YX}
+\mathbf{\Sigma}_{XY}\mathbf{A}^\top
\right),
\end{align*}
which is~\eqref{eq:mf_covariance}.

We again define our objective as the trace of the estimator covariance,
\[
\mathcal{J}(\mathbf{A})
=
\frac{1}{\nSamplesHF}\mathrm{tr}(\mathbf{\Sigma}_{XX})
+
\beta\Big(
\mathrm{tr}(\mathbf{A}\mathbf{\Sigma}_{YY}\mathbf{A}^\top)+2\mathrm{tr}(\mathbf{A}\mathbf{\Sigma}_{YX})
\Big),
\]
using \(\mathrm{tr}(\mathbf{A}\mathbf{\Sigma}_{YX})=\mathrm{tr}(\mathbf{\Sigma}_{XY}\mathbf{A}^\top)\).

\paragraph{Full matrix gain}
Using \(\nabla_{\mathbf{A}}\,\mathrm{tr}(\mathbf{A}\mathbf{\Sigma}_{YY}\mathbf{A}^\top)=2\mathbf{A}\mathbf{\Sigma}_{YY}\) and \(\nabla_{\mathbf{A}}\,\mathrm{tr}(\mathbf{A}\mathbf{\Sigma}_{YX})=\mathbf{\Sigma}_{YX}^\top=\mathbf{\Sigma}_{XY}\),
\[
\nabla_{\mathbf{A}} \mathcal{J}(\mathbf{A})=2\beta\left(\mathbf{A}\mathbf{\Sigma}_{YY}+\mathbf{\Sigma}_{XY}\right).
\]
Setting the gradient to be zero and solving gives 
\[
\mathbf{A}^\star=-\,\mathbf{\Sigma}_{XY}\mathbf{\Sigma}_{YY}^{-1}.
\]

\paragraph{Diagonal gain}
For the diagonal case \(\mathbf{A}=\mathbf{D}=\mathrm{diag}(a_1,\ldots,a_{\nParams})\), the objective simplifies to
\[
\mathcal{J}(\mathbf{D})=\frac{1}{\nSamplesHF}\mathrm{tr}(\mathbf{\Sigma}_{XX})
+\beta\sum_{i=1}^{\nParams}\left(a_i^2(\mathbf{\Sigma}_{YY})_{ii}+2a_i(\mathbf{\Sigma}_{XY})_{ii}\right).
\]
Each component decouples, so
\[
\frac{\partial \mathcal{J}}{\partial a_i}
=
2\beta\left(a_i(\mathbf{\Sigma}_{YY})_{ii}+(\mathbf{\Sigma}_{XY})_{ii}\right)=0
\quad\Rightarrow\quad
a_i^\star=-\,\frac{(\mathbf{\Sigma}_{XY})_{ii}}{(\mathbf{\Sigma}_{YY})_{ii}}.
\]
\paragraph{Scalar gain}
The following derivation applies to one fitted vector block. In VI, it is
applied separately to the mean and log-standard-deviation gradient
blocks, and to the vectorized mean--mean and log-standard-deviation--log-standard-deviation Hessian blocks, as specified in
Eqs.~\eqref{eq:scalar_gradient_gain_blocks} and~\eqref{eq:scalar_hessian_gain_blocks}.
Restricting \(\mathbf{A}=\alpha \mathbf{I}\), we obtain
\[
\mathcal{J}(\alpha)
=
\frac{1}{\nSamplesHF}\mathrm{tr}(\mathbf{\Sigma}_{XX})
+\beta\left(\alpha^2\mathrm{tr}(\mathbf{\Sigma}_{YY})+2\alpha\,\mathrm{tr}(\mathbf{\Sigma}_{XY})\right).
\]
Setting the gradient to zero and solving gives
\[
\frac{d\mathcal{J}}{d\alpha}
=
2\beta\left(\alpha\,\mathrm{tr}(\mathbf{\Sigma}_{YY})+\mathrm{tr}(\mathbf{\Sigma}_{XY})\right)=0
\quad\Rightarrow\quad
\alpha^\star=-\,\frac{\mathrm{tr}(\mathbf{\Sigma}_{XY})}{\mathrm{tr}(\mathbf{\Sigma}_{YY})}.
\]

\bibliographystyle{elsarticle-num}
\bibliography{refs}

@article{RiPaBl26,
	author = {Rizzi, Francesco and Parish, Eric J. and Blonigan, Patrick J. and Tencer, John T. and Wentland, Christopher R. and Kurzawski, Andrew J. and Van Heyningen, R. Loek and Schilly, Caleb and Brunini, Victor},
	date = {2026/04/27},
	doi = {10.1007/s00158-026-04300-2},
	id = {Rizzi2026},
	isbn = {1615-1488},
	journal = {Structural and Multidisciplinary Optimization},
	number = {5},
	pages = {112},
	title = {Projection-based model reduction using Pressio: application to hypersonic aerothermodynamics},
	url = {https://doi.org/10.1007/s00158-026-04300-2},
	volume = {69},
	year = {2026}}

@article{Br22,
	author = {Douglas J. Brinkerhoff},
	doi = {https://doi.org/10.1016/j.jcp.2022.111095},
	issn = {0021-9991},
	journal = {Journal of Computational Physics},
	pages = {111095},
	title = {Variational inference at glacier scale},
	url = {https://www.sciencedirect.com/science/article/pii/S0021999122001577},
	volume = {459},
	year = {2022}}

@misc{petersen2012matrixcookbook,
  title        = {The Matrix Cookbook},
  author       = {Petersen, Kaare Brandt and Pedersen, Michael Syskind},
  year         = {2012},
  note         = {Version 20121115},
  url          = {https://www.math.uwaterloo.ca/~hwolkowi/matrixcookbook.pdf}
}

@inproceedings{LiLiJi18,
author = {Li, Ximing and Li, Changchun and Chi, Jinjin and Ouyang, Jihong},
title = {Variance reduction in black-box variational inference by adaptive importance sampling},
year = {2018},
isbn = {9780999241127},
publisher = {AAAI Press},
booktitle = {Proceedings of the 27th International Joint Conference on Artificial Intelligence},
pages = {2404–2410},
numpages = {7},
location = {Stockholm, Sweden},
series = {IJCAI'18}
}

@article{AgRaKo25,
	author = {Agrawal, Atul and Ravi, Kislaya and Koutsourelakis, Phaedon-Stelios and Bungartz, Hans-Joachim},
	doi = {10.1088/2632-2153/ad8e2b},
	journal = {Machine Learning: Science and Technology},
	month = {jan},
	number = {1},
	pages = {015024},
	publisher = {IOP Publishing},
	title = {{Stochastic black-box optimization using multi-fidelity score function estimator}},
	url = {https://doi.org/10.1088/2632-2153/ad8e2b},
	volume = {6},
	year = {2025}}

@misc{ng2024pathwisegradientvariancereduction,
      title={Pathwise Gradient Variance Reduction with Control Variates in Variational Inference}, 
      author={Kenyon Ng and Susan Wei},
      year={2024},
      eprint={2410.05753},
      archivePrefix={arXiv},
      primaryClass={stat.ML},
      url={https://arxiv.org/abs/2410.05753}, 
}

@misc{liu2019steinvariationalgradientdescent,
      title={Stein Variational Gradient Descent: A General Purpose Bayesian Inference Algorithm}, 
      author={Qiang Liu and Dilin Wang},
      year={2019},
      eprint={1608.04471},
      archivePrefix={arXiv},
      primaryClass={stat.ML},
      url={https://arxiv.org/abs/1608.04471}, 
}

@book{Brooks_2011,
   title={Handbook of Markov Chain Monte Carlo},
   ISBN={9780429138508},
   url={http://dx.doi.org/10.1201/b10905},
   DOI={10.1201/b10905},
   publisher={Chapman and Hall/CRC},
   author={Brooks, Steve and Gelman, Andrew and Jones, Galin and Meng, Xiao-Li},
   year={2011},
   month=may }

@misc{kingma2022autoencodingvariationalbayes,
      title={Auto-Encoding Variational Bayes}, 
      author={Diederik P Kingma and Max Welling},
      year={2022},
      eprint={1312.6114},
      archivePrefix={arXiv},
      primaryClass={stat.ML},
      url={https://arxiv.org/abs/1312.6114}, 
}

@misc{rezende2014stochasticbackpropagationapproximateinference,
      title={Stochastic Backpropagation and Approximate Inference in Deep Generative Models}, 
      author={Danilo Jimenez Rezende and Shakir Mohamed and Daan Wierstra},
      year={2014},
      eprint={1401.4082},
      archivePrefix={arXiv},
      primaryClass={stat.ML},
      url={https://arxiv.org/abs/1401.4082}, 
}

@misc{schulman2016gradientestimationusingstochastic,
      title={Gradient Estimation Using Stochastic Computation Graphs}, 
      author={John Schulman and Nicolas Heess and Theophane Weber and Pieter Abbeel},
      year={2016},
      eprint={1506.05254},
      archivePrefix={arXiv},
      primaryClass={cs.LG},
      url={https://arxiv.org/abs/1506.05254}, 
}

@article{WaJo08,
author = {Wainwright, Martin J. and Jordan, Michael I.},
title = {{Graphical Models, Exponential Families, and Variational Inference}},
year = {2008},
issue_date = {January 2008},
publisher = {Now Publishers Inc.},
address = {Hanover, MA, USA},
volume = {1},
number = {1–2},
issn = {1935-8237},
url = {https://doi.org/10.1561/2200000001},
doi = {10.1561/2200000001},
journal = {Found. Trends Mach. Learn.},
month = jan,
pages = {1–305},
numpages = {305}
}

@article{BlKuMc17,
   title={Variational Inference: A Review for Statisticians},
   volume={112},
   ISSN={1537-274X},
   url={http://dx.doi.org/10.1080/01621459.2017.1285773},
   DOI={10.1080/01621459.2017.1285773},
   number={518},
   journal={Journal of the American Statistical Association},
   publisher={Informa UK Limited},
   author={Blei, David M. and Kucukelbir, Alp and McAuliffe, Jon D.},
   year={2017},
   month=apr, pages={859–877} }

@article{Fu06,
  title = {{Gradient estimation}},
  author = {Fu, Michael C},
  journal = {Handbooks in operations research and management science},
  volume = {13},
  pages = {575--616},
  year = {2006},
  publisher = {Elsevier}
}

@article{Wi92,
	author = {Williams, Ronald J. },
	date = {1992/05/01},
	doi = {10.1007/BF00992696},
	id = {Williams1992},
	isbn = {1573-0565},
	journal = {Machine Learning},
	number = {3},
	pages = {229--256},
	title = {{Simple statistical gradient-following algorithms for connectionist reinforcement learning}},
	url = {https://doi.org/10.1007/BF00992696},
	volume = {8},
	year = {1992}}

@inbook{Jo98,
	address = {Dordrecht},
	author = {Jordan, Michael I. and Ghahramani, Zoubin and Jaakkola, Tommi S. and Saul, Lawrence K.},
	booktitle = {Learning in Graphical Models},
	doi = {10.1007/978-94-011-5014-9_5},
	editor = {Jordan, Michael I.},
	isbn = {978-94-011-5014-9},
	pages = {105--161},
	publisher = {Springer Netherlands},
	title = {{An Introduction to Variational Methods for Graphical Models}},
	url = {https://doi.org/10.1007/978-94-011-5014-9_5},
	year = {1998}}

@article{GORODETSKY2020109257,
	author = {Alex A. Gorodetsky and Gianluca Geraci and Michael S. Eldred and John D. Jakeman},
	doi = {https://doi.org/10.1016/j.jcp.2020.109257},
	issn = {0021-9991},
	journal = {Journal of Computational Physics},
	pages = {109257},
	title = {{A generalized approximate control variate framework for multifidelity uncertainty quantification}},
	url = {https://www.sciencedirect.com/science/article/pii/S0021999120300310},
	volume = {408},
	year = {2020}}

@article{ZhChZh24,
author = {Wu, Zhizhang and Zhang, Cheng and Zhang, Zhiwen},
title = {{Reduced-order model-based variational inference with normalizing flows for Bayesian elliptic inverse problems}},
year = {2024},
issue_date = {May 2024},
publisher = {Elsevier Science Publishers B. V.},
address = {NLD},
volume = {441},
number = {C},
issn = {0377-0427},
url = {https://doi.org/10.1016/j.cam.2023.115659},
doi = {10.1016/j.cam.2023.115659},
journal = {J. Comput. Appl. Math.},
month = may,
numpages = {17}
}

@article{rubinstein1985efficiency,
	Author = {Rubinstein, R. Y. and Marcus, R.},
	journal = {Oper. Res.},
	fjournal = {Operations Research},
	Number = {3},
	Pages = {661--677},
	Publisher = {INFORMS},
	Title = {Efficiency of multivariate control variates in {M}onte {C}arlo simulation},
	Volume = {33},
	Year = {1985}
}

@article{PeWiGu16,
	author = {Peherstorfer, Benjamin and Willcox, Karen and Gunzburger, Max},
	doi = {10.1137/15M1046472},
	eprint = {https://doi.org/10.1137/15M1046472},
	journal = {SIAM Journal on Scientific Computing},
	number = {5},
	pages = {A3163-A3194},
	title = {{Optimal Model Management for Multifidelity Monte Carlo Estimation}},
	url = {https://doi.org/10.1137/15M1046472},
	volume = {38},
	year = {2016}}

@inproceedings{Ranganath2014BBVI,
  title     = {Black Box Variational Inference},
  author    = {Ranganath, Rajesh and Gerrish, Sean and Blei, David M.},
  booktitle = {Proceedings of the Seventeenth International Conference on Artificial Intelligence and Statistics},
  series    = {Proceedings of Machine Learning Research},
  volume    = {33},
  pages     = {814--822},
  publisher = {PMLR},
  year      = {2014}
}

@incollection{ZhNaZh21,
	author = {Xin Zhang and Muhammad Atif Nawaz and Xuebin Zhao and Andrew Curtis},
	booktitle = {{Inversion of Geophysical Data}},
	doi = {https://doi.org/10.1016/bs.agph.2021.06.003},
	editor = {Cedric Schmelzbach},
	issn = {0065-2687},
	pages = {73-140},
	publisher = {Elsevier},
	series = {Advances in Geophysics},
	title = {{Chapter Two - An introduction to variational inference in geophysical inverse problems}},
	url = {https://www.sciencedirect.com/science/article/pii/S0065268721000030},
	volume = {62},
	year = {2021}}

@article{AgKo24,
	author = {Atul Agrawal and Phaedon-Stelios Koutsourelakis},
	doi = {https://doi.org/10.1016/j.jcp.2024.112982},
	issn = {0021-9991},
	journal = {Journal of Computational Physics},
	pages = {112982},
	title = {{A probabilistic, data-driven closure model for RANS simulations with aleatoric, model uncertainty}},
	url = {https://www.sciencedirect.com/science/article/pii/S0021999124002316},
	volume = {508},
	year = {2024}}

@inproceedings{Wierstra2008NES,
  author    = {Wierstra, Daan and Schaul, Tom and Peters, Jan and Schmidhuber, J{\"u}rgen},
  title     = {{Natural Evolution Strategies}},
  booktitle = {{Proceedings of the IEEE Congress on Evolutionary Computation (CEC)}},
  year      = {2008},
  pages     = {3381--3387},
  doi       = {10.1109/CEC.2008.4631255}
}

@article{Wierstra2014NES,
  author  = {Wierstra, Daan and Schaul, Tom and Peters, Jan and Schmidhuber, J{\"u}rgen},
  title   = {Natural Evolution Strategies},
  journal = {Journal of Machine Learning Research},
  volume  = {15},
  pages   = {949--980},
  year    = {2014}
}

@article{Amari1998NaturalGradient,
  author  = {Amari, Shun-ichi},
  title   = {Natural Gradient Works Efficiently in Learning},
  journal = {Neural Computation},
  volume  = {10},
  number  = {2},
  pages   = {251--276},
  year    = {1998},
  doi     = {10.1162/089976698300017746}
}

@misc{paisley2012variationalbayesianinferencestochastic,
      title={Variational Bayesian Inference with Stochastic Search}, 
      author={John Paisley and David Blei and Michael Jordan},
      year={2012},
      eprint={1206.6430},
      archivePrefix={arXiv},
      primaryClass={cs.LG},
      url={https://arxiv.org/abs/1206.6430}, 
}

@article{CaMaYo15,
	author = {Cui, Tiangang and Marzouk, Youssef M. and Willcox, Karen E.},
	doi = {https://doi.org/10.1002/nme.4748},
	eprint = {https://onlinelibrary.wiley.com/doi/pdf/10.1002/nme.4748},
	journal = {International Journal for Numerical Methods in Engineering},
	number = {5},
	pages = {966-990},
	title = {{Data-driven model reduction for the Bayesian solution of inverse problems}},
	url = {https://onlinelibrary.wiley.com/doi/abs/10.1002/nme.4748},
	volume = {102},
	year = {2015}}

@article{ChFo05,
	author = {J. Andr{\'e}s Christen and Colin Fox},
	doi = {10.1198/106186005X76983},
	eprint = {https://doi.org/10.1198/106186005X76983},
	journal = {Journal of Computational and Graphical Statistics},
	number = {4},
	pages = {795--810},
	publisher = {Taylor \& Francis},
	title = {{Markov chain Monte Carlo Using an Approximation}},
	url = {https://doi.org/10.1198/106186005X76983},
	volume = {14},
	year = {2005}}

@article{DoYa22,
	author = {Geoff Donoghue and Masayuki Yano},
	doi = {https://doi.org/10.1016/j.cma.2022.115282},
	issn = {0045-7825},
	journal = {Computer Methods in Applied Mechanics and Engineering},
	pages = {115282},
	title = {{A multi-fidelity ensemble Kalman filter with hyperreduced reduced-order models}},
	url = {https://www.sciencedirect.com/science/article/pii/S0045782522004005},
	volume = {398},
	year = {2022}}

@article{popov2021multifidelity,
  title = {{A multifidelity ensemble Kalman filter with reduced order control variates}},
  author = {Popov, A. A. and Mou, C. and Iliescu, T. and Sandu, A.},
  journal = {SIAM J. Sci. Comput.},
  fjournal = {SIAM Journal on Scientific Computing},
  volume = {43},
  number = {2},
  pages = {A1134--A1162},
  year = {2021},
  publisher = {SIAM}
}

@article{chada2023review,
  title = {A Review of the {EnKF} for Parameter Estimation},
  author = {Chada, N. K.},
  journal = {Inverse Problems: Recent Advances and Applications},
  year = {2023},
  publisher = {IntechOpen}
}

@article{yan2019adaptive,
  title = {{An adaptive multifidelity PC-based ensemble Kalman inversion for inverse problems}},
  author = {Yan, L. and Zhou, T.},
  journal = {Int. J. Uncertain. Quantif.},
  fjournal = {International Journal for Uncertainty Quantification},
  volume = {9},
  number = {3},
  year = {2019},
  publisher = {Begel House Inc.}
}

@article{gao2021bi,
  title = {A bi-fidelity ensemble Kalman method for PDE-constrained inverse problems in computational mechanics},
  author = {Gao, H. and Wang, J.-X.},
  journal = {Comput. Mech.},
  fjournal = {Computational Mechanics},
  volume = {67},
  number = {4},
  pages = {1115--1131},
  year = {2021},
  publisher = {Springer}
}

@book{evensen2009data,
  title = {{Data assimilation: the ensemble Kalman filter}},
  author = {Evensen, G.},
  year = {2009},
  publisher = {Springer}
}

@book{kalnay2003atmospheric,
  title = {{Atmospheric modeling, data assimilation, and predictability}},
  author = {Kalnay, E.},
  year = {2003},
  publisher = {Cambridge University Press}
}

@article{hoel2016multilevel,
  title = {{Multilevel ensemble Kalman filtering}},
  author = {Hoel, H. and Law, K. J. H. and Tempone, R.},
  journal = {SIAM J. Numer. Anal.},
  fjournal = {SIAM Journal on Numerical Analysis},
  volume = {54},
  number = {3},
  pages = {1813--1839},
  year = {2016},
  publisher = {SIAM}
}

@article{rei2025blackboxvariationalinference,
      title={A Black Box Variational Inference Scheme for Inverse Problems with Demanding Physics-Based Models}, 
      author={G. Robalo Rei and C. P. Schmidt and J. Nitzler and M. Dinkel and W. A. Wall},
      year={2026},
      journal={Advanced Modeling and Simulation in Engineering Sciences},
      volume={13},
      number={19},
      doi={10.1186/s40323-026-00335-8},
      eprint={2510.25038},
      archivePrefix={arXiv},
      primaryClass={cs.CE},
      url={https://doi.org/10.1186/s40323-026-00335-8},
}

@article{Gl90,
  title = {{Likelihood ratio gradient estimation for stochastic systems}},
  author = {Glynn, Peter W},
  journal = {Communications of the ACM},
  volume = {33},
  number = {10},
  pages = {75--84},
  year = {1990},
  publisher = {ACM}}

@article{WaZa05,
	author = {Jingbo Wang and Nicholas Zabaras},
	doi = {https://doi.org/10.1016/j.ijheatmasstransfer.2004.08.009},
	issn = {0017-9310},
	journal = {International Journal of Heat and Mass Transfer},
	number = {1},
	pages = {15-29},
	title = {{Using Bayesian statistics in the estimation of heat source in radiation}},
	url = {https://www.sciencedirect.com/science/article/pii/S0017931004003709},
	volume = {48},
	year = {2005}}

@article{LiWiGh10,
	author = {Lieberman, Chad and Willcox, Karen and Ghattas, Omar},
	doi = {10.1137/090775622},
	eprint = {https://doi.org/10.1137/090775622},
	journal = {SIAM Journal on Scientific Computing},
	number = {5},
	pages = {2523-2542},
	title = {{Parameter and State Model Reduction for Large-Scale Statistical Inverse Problems}},
	url = {https://doi.org/10.1137/090775622},
	volume = {32},
	year = {2010}}

@article{GaFiWi09,
	author = {Galbally, D. and Fidkowski, K. and Willcox, K. and Ghattas, O.},
	doi = {https://doi.org/10.1002/nme.2746},
	eprint = {https://onlinelibrary.wiley.com/doi/pdf/10.1002/nme.2746},
	journal = {International Journal for Numerical Methods in Engineering},
	number = {12},
	pages = {1581-1608},
	title = {{Non-linear model reduction for uncertainty quantification in large-scale inverse problems}},
	url = {https://onlinelibrary.wiley.com/doi/abs/10.1002/nme.2746},
	volume = {81},
	year = {2010}}

@book{wasserman2004all,
  title={All of Statistics: A Concise Course in Statistical Inference},
  author={Wasserman, Larry},
  year={2004},
  publisher={Springer}
}

@inbook{BuWi10,
    author = {Marcelo Buffoni and Karen Willcox},
    booktitle = {40th Fluid Dynamics Conference and Exhibit},
    doi = {10.2514/6.2010-5008},
    eprint = {https://arc.aiaa.org/doi/pdf/10.2514/6.2010-5008},
    title = {{Projection-based} model reduction for reacting flows},
    url = {https://arc.aiaa.org/doi/abs/10.2514/6.2010-5008}}

@article{GaWaZh20,
    author = {Han Gao and Jian-Xun Wang and Matthew J. Zahr},
    doi = {https://doi.org/10.1016/j.physd.2020.132614},
    issn = {0167-2789},
    journal = {Physica D: Nonlinear Phenomena},
    pages = {132614},
    title = {{Non-intrusive} model reduction of large-scale, nonlinear dynamical systems using deep learning},
    url = {https://www.sciencedirect.com/science/article/pii/S0167278919305573},
    volume = {412},
    year = {2020}}

@article{LeCa20,
    author = {Kookjin Lee and Kevin T. Carlberg},
    doi = {https://doi.org/10.1016/j.jcp.2019.108973},
    issn = {0021-9991},
    journal = {Journal of Computational Physics},
    pages = {108973},
    title = {{Model} reduction of dynamical systems on nonlinear manifolds using deep convolutional autoencoders},
    url = {https://www.sciencedirect.com/science/article/pii/S0021999119306783},
    volume = {404},
    year = {2020}}

@misc{dayta2024variancecontrolblackbox,
      title={Variance Control for Black Box Variational Inference Using The James-Stein Estimator}, 
      author={Dominic B. Dayta},
      year={2024},
      eprint={2405.05485},
      archivePrefix={arXiv},
      primaryClass={cs.LG},
      url={https://arxiv.org/abs/2405.05485}, 
}

@misc{
kool2019buy,
title={Buy 4 {REINFORCE} Samples, Get a Baseline for Free!},
author={Wouter Kool and Herke van Hoof and Max Welling},
year={2019},
url={https://openreview.net/forum?id=r1lgTGL5DE}
}

@article{FuSa21,
  author  = {Fujisawa, Masahiro and Sato, Issei},
  title   = {{Multilevel Monte Carlo Variational Inference}},
  journal = {Journal of Machine Learning Research},
  volume  = {22},
  number  = {278},
  pages   = {1--44},
  year    = {2021},
  url     = {https://jmlr.org/papers/v22/20-653.html}
}

@misc{LiYaChHeSu25,
  author        = {Li, Kaiyu and Yang, Yiming and Cheng, Xiaoyuan and He, Yi and Sun, Zhuo},
  title         = {{Multilevel Control Functional}},
  year          = {2025},
  eprint        = {2305.12996},
  archivePrefix = {arXiv},
  primaryClass  = {stat.ME},
  url           = {https://arxiv.org/abs/2305.12996}
}

@article{ChChHuHuWa25,
  author  = {Che, Baojun and Chen, Yifan and Huan, Zhenghao and Huang, Daniel Zhengyu and Wang, Weijie},
  title   = {{Stable Derivative Free Gaussian Mixture Variational Inference for Bayesian Inverse Problems}},
  journal = {SIAM Journal on Scientific Computing},
  volume  = {47},
  number  = {5},
  pages   = {A2583--A2608},
  year    = {2025},
  doi     = {10.1137/25M1722548},
  url     = {https://doi.org/10.1137/25M1722548}
}

@article{MoRoFiMn20,
  author  = {Mohamed, Shakir and Rosca, Mihaela and Figurnov, Michael and Mnih, Andriy},
  title   = {{Monte Carlo Gradient Estimation in Machine Learning}},
  journal = {Journal of Machine Learning Research},
  volume  = {21},
  number  = {132},
  pages   = {1--62},
  year    = {2020},
  url     = {https://jmlr.org/papers/v21/19-346.html}
}

@inproceedings{RiBoNuRuAk20,
  author    = {Richter, Lorenz and Boustati, Ayman and N{\"u}sken, Nikolas and Ruiz, Francisco J. R. and Akyildiz, {\"O}mer Deniz},
  title     = {{VarGrad: A Low-Variance Gradient Estimator for Variational Inference}},
  booktitle = {Advances in Neural Information Processing Systems},
  volume    = {33},
  pages     = {13481--13492},
  year      = {2020},
  url       = {https://proceedings.neurips.cc/paper/2020/hash/9c22c0b51b3202246463e986c7e205df-Abstract.html}
}

@inproceedings{FaWaBeKwHe15,
  author    = {Fan, Kai and Wang, Ziteng and Beck, Jeffrey and Kwok, James T. and Heller, Katherine A.},
  title     = {{Fast Second-Order Stochastic Backpropagation for Variational Inference}},
  booktitle = {Advances in Neural Information Processing Systems},
  volume    = {28},
  pages     = {1387--1395},
  year      = {2015},
  url       = {https://proceedings.neurips.cc/paper/2015/hash/fc3cf452d3da8402bebb765225ce8c0e-Abstract.html}
}

@inproceedings{ReJoMc17,
  author        = {Regier, Jeffrey and Jordan, Michael I. and McAuliffe, Jon},
  title         = {{Fast Black-box Variational Inference through Stochastic Trust-Region Optimization}},
  year          = {2017},
  booktitle     = {Advances in Neural Information Processing Systems},
  volume        = {30},
  eprint        = {1706.02375},
  archivePrefix = {arXiv},
  primaryClass  = {stat.ML},
  url           = {https://proceedings.neurips.cc/paper_files/paper/2017/hash/9a1756fd0c741126d7bbd4b692ccbd91-Abstract.html}
}

@inproceedings{FoFaAlRoXiWh18,
  author    = {Foerster, Jakob and Farquhar, Gregory and Al-Shedivat, Maruan and Rockt{\"a}schel, Tim and Xing, Eric P. and Whiteson, Shimon},
  title     = {{DiCE: The Infinitely Differentiable Monte Carlo Estimator}},
  booktitle = {Proceedings of the 35th International Conference on Machine Learning},
  series    = {Proceedings of Machine Learning Research},
  volume    = {80},
  pages     = {1529--1538},
  year      = {2018},
  url       = {https://proceedings.mlr.press/v80/foerster18a.html}
}

@inproceedings{CoZhLiChCa21,
  author    = {Cong, Yulai and Zhao, Miaoyun and Li, Jianqiao and Chen, Junya and Carin, Lawrence},
  title     = {{GO Hessian for Expectation-Based Objectives}},
  booktitle = {Proceedings of the AAAI Conference on Artificial Intelligence},
  volume    = {35},
  pages     = {12060--12068},
  year      = {2021},
  doi       = {10.1609/aaai.v35i13.17432},
  url       = {https://doi.org/10.1609/aaai.v35i13.17432}
}

@article{MaShIo26,
  author = {Magris, Martin and Shabani, Mostafa and Iosifidis, Alexandros},
  title = {Quasi black-box variational inference for {Bayesian} learning},
  journal = {Pattern Recognition},
  volume = {177},
  pages = {113219},
  year = {2026},
  doi = {10.1016/j.patcog.2026.113219},
  url = {https://doi.org/10.1016/j.patcog.2026.113219}
}

@inproceedings{JiSuSu21,
  author = {Ji, Geng and Sujono, Debora and Sudderth, Erik B.},
  title = {Marginalized Stochastic Natural Gradients for Black-Box Variational Inference},
  booktitle = {Proceedings of the 38th International Conference on Machine Learning},
  series = {Proceedings of Machine Learning Research},
  volume = {139},
  pages = {4870--4881},
  year = {2021},
  url = {https://proceedings.mlr.press/v139/ji21b.html}
}

@inproceedings{TrCoMe18,
  author = {Trusheim, Felix and Condurache, Alexandru Paul and Mertins, Alfred},
  title = {Boosting Black-Box Variational Inference by Incorporating the Natural Gradient},
  booktitle = {2018 24th International Conference on Pattern Recognition (ICPR)},
  year = {2018},
  url = {https://www.isip.uni-luebeck.de/fileadmin/files/Trusheim2018.pdf}
}

@misc{ReMc16,
  author = {Regier, Jeffrey and McAuliffe, Jon},
  title = {Second-order stochastic variational inference},
  year = {2016},
  note = {Author manuscript},
  url = {https://sites.lsa.umich.edu/regier/wp-content/uploads/sites/1349/2024/09/regier2016second.pdf}
}

@article{LiOw21,
  author = {Liu, Sifan and Owen, Art B.},
  title = {Quasi-{Monte Carlo} Quasi-{Newton} in Variational {Bayes}},
  journal = {Journal of Machine Learning Research},
  volume = {22},
  number = {243},
  pages = {1--23},
  year = {2021},
  url = {https://jmlr.org/papers/v22/21-0498.html}
}

@article{Ru89,
  author = {Rubinstein, Reuven Y.},
  title = {Sensitivity Analysis and Performance Extrapolation for Computer Simulation Models},
  journal = {Operations Research},
  volume = {37},
  number = {1},
  pages = {72--81},
  year = {1989},
  doi = {10.1287/opre.37.1.72},
  url = {https://doi.org/10.1287/opre.37.1.72}
}

@book{Gl04,
  author = {Glasserman, Paul},
  title = {{Monte Carlo} Methods in Financial Engineering},
  publisher = {Springer},
  year = {2004},
  doi = {10.1007/978-0-387-21617-1},
  url = {https://doi.org/10.1007/978-0-387-21617-1}
}

@article{FuLeBa16,
  author = {Furmston, Thomas and Lever, Guy and Barber, David},
  title = {Approximate {Newton} Methods for Policy Search in {Markov} Decision Processes},
  journal = {Journal of Machine Learning Research},
  volume = {17},
  number = {226},
  pages = {1--51},
  year = {2016},
  url = {https://jmlr.org/papers/v17/15-414.html}
}

@misc{MaMoPrBh23,
  author = {Maniyar, Mizhaan Prajit and Mondal, Akash and {Prashanth L. A.} and Bhatnagar, Shalabh},
  title = {A Cubic-regularized Policy {Newton} Algorithm for Reinforcement Learning},
  year = {2023},
  eprint = {2304.10951},
  archivePrefix = {arXiv},
  url = {https://arxiv.org/abs/2304.10951}
}

\end{document}